\documentclass{amsart}

\usepackage{amsmath,amsthm,amsfonts,amssymb,bm,graphicx,mathrsfs, dsfont, mathtools}

\usepackage{comment}
\usepackage
{hyperref}

\hypersetup{colorlinks=true,citecolor=blue,linkcolor=blue,urlcolor=blue,
pdfstartview=FitH }

\theoremstyle{plain}
\newtheorem{lemma}{Lemma}
\newtheorem{theorem}[lemma]{Theorem}

\newtheorem{prop}[lemma]{Proposition}

\numberwithin{lemma}{section}
\numberwithin{equation}{section}

\theoremstyle{definition}

\theoremstyle{remark}

\renewcommand{\pmod}[1]{\,(\textup{mod}\,#1)}
\mathtoolsset{showonlyrefs=true}
\newcommand{\N}{\mathfrak{N}}
\newcommand{\bd}{\bold{d}}

\begin{document}
\author{Valentin Blomer}
\author{Lasse Grimmelt}
\author{Junxian Li}
\author{Simon L. Rydin Myerson}
 
\address{Mathematisches Institut, Endenicher Allee 60, 53115 Bonn}
\email{blomer@math.uni-bonn.de}
\address{Mathematical Institute, Andrew Wiles Building. University of Oxford, Oxford OX2 6GG}
\email{grimmelt@maths.ox.ac.uk}
\address{Mathematisches Institut, Endenicher Allee 60, 53115 Bonn}
\email{jli135@math.uni-bonn.de}
\address{Mathematics Institute, Zeeman Building. University of Warwick, Coventry CV4 7AL}
\email{simon.rydin-myerson@warwick.ac.uk}

\title{Additive problems with almost prime squares}

\thanks{The first and third author were supported in part by DFG grants BL 915/2-2 and BL 915/5-1 and Germany's Excellence Strategy grant EXC-2047/1 - 390685813. The second author received funding from the European Research Council (ERC) under the 
	European Union's Horizon 2020 research and innovation programme, grant 
	agreement no. 851318. The third and fourth author acknowledge support of the Max Planck Institute for Mathematics. The fourth author was supported by DFG project number 255083470 and by a Leverhulme Early Career Fellowship, and received support from the European Reseach Council (ERC) under the European Union's Horizon 2020 research and innovation program (grant ID 648329).
}

\begin{abstract}  We show that every sufficiently large integer is a sum of a prime and two almost prime squares, and also a sum of a smooth number and two almost prime squares.  The number of such representations is of the expected order of magnitude. We likewise treat representations of shifted primes \(p-1\) as sums of two almost prime squares. The methods involve a combination of analytic, automorphic and algebraic arguments to handle representations by restricted binary quadratic forms with a high degree of uniformity.
	
\end{abstract}


\setcounter{tocdepth}{2}  \maketitle 

\maketitle

\section{Introduction} 

\subsection{Statement of results}

In his first paper, Hooley \cite{Ho} gave a proof, conditional upon the Generalized Riemann Hypothesis,  of a conjecture made by Hardy and Littlewood some 30 years before \cite{HL}: the number of representations of an integer $n$ in the form 
\begin{equation}\label{1}
   n=p+x_1^2+x_2^2
\end{equation}  
    in primes $p$ and non-zero integers $x_1, x_2$ as $n \rightarrow \infty$ is asymptotic to
\begin{equation}\label{12}
\pi \, \mathfrak S(n){\rm Li}(n),\quad \text{ where } \mathfrak S(n)= \prod_{p}\Big(1+\frac{\chi_{-4}(p)}{p(p-1)}\Big)\prod_{p\mid n}\frac{(p-1)(p-\chi_{-4}(p))}{p^2-p+\chi_{-4}(p)}. 
\end{equation}
In particular every sufficiently large number is representable in this form. The Generalized Riemann Hypothesis (GRH) was removed a few years later by Linnik  \cite{Li}   who used his dispersion method to obtain the first unconditional proof of the Hardy-Littlewood problem. The proof could be greatly streamlined and simplified once the Bombieri-Vinogradov theorem became available, cf.\ e.g.\ \cite{HoBook}. 
The best error term $O_A(n (\log n)^{-A})$ in the asymptotic formula was recently established by  Assing, Blomer and Li \cite{ABL} by invoking in addition the theory of automorphic forms. 

In his last paper, Hooley \cite{Ho2} returned to this topic and considered a refined version. He showed that every sufficiently large number is representable as a sum of a prime and two squares of \emph{square-free} numbers and gave a lower bound for the number of such representations of the expected order of magnitude. At first sight this appears to be easy -- square-free numbers can be detected by zero-dimensional sieve, so only a little bit of inclusion-exclusion is required. On second thought, a major issue emerges: to sieve out even a fixed number of squares of primes, one needs an asymptotic formula for \eqref{1} with $d_1 \mid x_1$, $d_2 \mid x_2$. However, 
 the quadratic form $d_1^2 x_1^2 + d_2^2 x_2^2$ fails to have class number one, in fact very soon it even fails to have one class per genus, and hence there is no convolution formula for the number of ways it represents a given integer. Without such a formula, all approaches to the Hardy-Littlewood problem and its variations become very problematic. Hooley finds an ingenious (but rather ad hoc) way to circumvent this issue altogether for the particular problem at hand. 

In this paper, we continue on this avenue. We introduce a general sieve process into the problem with the goal of showing that every sufficiently large integer is a sum of a prime and two squares of \emph{almost primes}. The density of such representations requires us to have a good error term (a saving of at least $(\log x)^{3+\epsilon}$) in solutions to \eqref{1}, which ultimately requires the study of primes in arithmetic progressions to large moduli that go  beyond the range that GRH could handle individually. More importantly, in stark contrast to Hooley's result, the application of a sieve to make $x_1, x_2$ almost prime requires the coefficients $d_1$ and $d_2$ to be as large as some fixed power of $n$. In addition to the fact that the quadratic forms $d_1^2x_1^2 + d_2^2x_2^2$ fail to have class number one, this large range of uniformity for individual $d_1,d_2$ seems to require inevitably a count of solutions of \eqref{1} with a power-saving error term, which is only obtainable under some form of GRH.  

We will give more details on the combination of algebraic, automorphic and analytic methods involved in a moment, but emphasize at this point that we are nevertheless able to solve this problem unconditionally. 
  
\begin{theorem}\label{thm1} There exists a constant $C > 0$ such that every sufficiently large integer  $n\equiv 1, 3\bmod 6$ can be represented in the form \eqref{1}, where \(p\) is a prime and  $x_1, x_2$ are integers all of whose prime factors are greater than $n^{1/C}$. The number of such representations is $\gg \mathfrak S(n) n (\log n)^{-3}$ with $\mathfrak{S}(n)$ as in \eqref{12}. In particular, every sufficiently large integer $n\equiv 1, 3 \mod 6$ can be written as the sum of a prime and two squares of integers with no more than $C$ prime factors. 
\end{theorem}

Theorem \ref{thm1} gives the expected order of magnitude given the restrictions on $x_1,x_2$. We have made no effort to compute and optimize the value of $C$, although an upper bound of size $C \leq 10^6$, say, can easily be obtained by careful book-keeping.
Note that the congruence condition on $n$ is necessary to guarantee solutions to \eqref{1} with $(x_1x_2, 6)=1$ and \(p\) a prime greater than 3. At the cost of slightly more work it is possible to remove the condition modulo 6 if \(x_1,x_2\) are permitted to have the small primes  \(2\) and \(3\) as factors.
Going further, our method can be generalized to deal with general fixed quadratic forms $F(x_1, x_2)$ of fundamental discriminant, instead of $x_1^2 + x_2^2$, cf.\  the remark after Lemma \ref{lem1}.

The technique used in the proof of Theorem \ref{thm1} also gives the following result. 

\begin{theorem}\label{thm1'}
There exists a constant $C>0$ such that the number of solutions to the equation
\begin{equation}\label{FI-equation}
p = x_1^2 + x_2^2 + 1
\end{equation}
in primes $p \leq x$ and $x_1, x_2$ all of whose prime factors are greater  than $p^{1/C}$ is $\gg x (\log x)^{-3}$. In particular, 
there are infinitely many  primes  shifted by one that can be written as two squares of almost primes. 
\end{theorem}

This should be compared with a beautiful recent result of Friedlander-Iwaniec \cite{FI2} who  considered \eqref{FI-equation} for a prime $x_1$ and an almost prime $x_2$, but without the additive shift. The multiplicativity in Gau\ss ian primes plays a crucial rule in their proof. Adding a  shift parameter ruins the multiplicative structure and so it is no surprise that our proof takes a different route, see Section \ref{sec12}. Note that Theorem \ref{thm1'} is slightly easier than Theorem \ref{thm1} as the shift parameter is fixed (although of course one can prove a similar result with a more general shift $f$ uniformly for $|f|\ll x$).

From a multiplicative point of view, complementary to primes are smooth numbers (``entiers friables'' in French) all of whose prime factors are small. The smallness is measured in terms of a parameter $y$, and following the usual notation we let $\Psi(x, y) = \#S(x, y)$ where $S(x, y)$ denotes the set of positive integers $m \leq x$ such that $p\mid m$ implies $p \leq y$ for all primes $p$. Smooth numbers can appear quite frequently (e.g $\Psi(x,y)\gg x$ if $\log x \asymp \log y$), yet additive problems with smooth numbers are often surprisingly difficult, see e.g.\ \cite{Ba, BBD, KMS}. The present case is no exception. 
Here we consider the equation
\begin{equation}\label{1.2}
n=m+x_1^2+x_2^2
\end{equation}  
 where $x_1, x_2$ are almost primes and $m$ is $y$-smooth with $y \geq \mathcal{L}(n) := (\log n)^D$ for some sufficiently large constant $D$. We remark that the smoothing parameter \(y\) could be much smaller than those that are admissible in \cite{Ba, BBD, KMS}, which makes the problem quite delicate due to the sparseness of smooth numbers when $y$ becomes such small. We have that 
 \begin{align}
 \Psi(n,y)\asymp n u^{-(1+o(1))u}, \quad u=\frac{\log n}{\log y}.
 \end{align}
  In particular, we have 
 $$\Psi(n,(\log n)^D)=n^{1-1/D+o(1)}.$$ To capture smooth numbers with such small density, we need to define the quantity
 \begin{align}\label{Fdef}
\mathcal F(n, y)=\prod_{p\mid n}\Big(1-\frac{\chi_{-4}(p)}{p}\Big)\prod_{\substack{p\mid n\\ p\leq y}}\Big(1+\frac{\chi_{-4}(p)}{p^{\alpha}}\Big)
\end{align} 
where  $\alpha=\alpha(n,y)$ is the unique positive solution to the equation 
\begin{align}\label{alphadef}
\sum_{p\leq y}\frac{\log p}{p^{\alpha}-1}=\log n.
\end{align} 

\begin{theorem}\label{thm2} 
There exist  constants $C , D> 0$ 
 such that for any function $g$ with 
$(\log n)^D\leq g(n)\leq n$,  every sufficiently large  integer $n$ can be represented in the form \eqref{1.2} with $m \in S(n, g(n))$ and integers $x_1, x_2$ all of whose prime factors are greater than \(n^{1/C}\). The number of such representations is $\gg  \mathcal F(n, g(n))\Psi(n,g(n))(\log n)^{-2}$. 
 \end{theorem}

In fact, in the representations we construct, the smooth number $m$ will in addition be a square-free number times a power of 2. For orientation we remark that $$\alpha=1-\frac{\log (u\log u)+O(1)}{\log y}\gg 1-1/D+o(1) \asymp 1$$ for $y \geq \mathcal{L}(n)$.  The arithmetic factor could vary quite a lot when $y$ is small as
 $$\exp(-(\log n)^{1/D+o(1)})\ll\mathcal F(n, (\log n)^D)\ll \exp((\log n)^{1/D+o(1)})$$ and whilst when  $\log y\geq (\log_2n)^2$, we have $\mathcal F(n, y)\asymp 1$.

We mention one possible path for further work, namely Lagrange's four square theorem with restricted variables. That is, one seeks to solve $N=x_1^2+x_2^2+x_3^2+x_4^2$ with $x_i$ restricted to thin subsets of the integers such as primes, almost primes, or smooth numbers. See for instance   Br\"udern-Fouvry \cite{BF}, Heath-Brown-Tolev \cite{HBT} and Blomer-Br\"udern-Dietmann \cite{BBD}. It is feasible that the techniques in this paper can be helpful to study related problems.

\subsection{Methods and ideas}\label{sec12}  While the statements of Theorems \ref{thm1} and \ref{thm2} are purely arithmetic, the proofs combine a variety of different methods from different fields. We need the algebraic structure theory of orders in quadratic number fields, the full machinery of automorphic forms including some higher rank $L$-functions, and of course the toolbox of analytic number theory. 

Starting with the latter, the obvious approach to almost primes is the application of a sieve. Thus, for the proof of both theorems we need to analyze quite accurately the quantity
\begin{equation}\label{2}
\sum_{\ell} r_{d_1, d_2}(n-\ell), \quad r_{d_1, d_2}(n) = \#\{(x, y) \in \Bbb{Z}^2 \mid d_1^2 x_1^2 + d_2^2 x_2^2 = n\}
\end{equation}
for square-free integers $d_1, d_2$ where $\ell$ runs either over primes or over smooth numbers. As mentioned above, in general there is no closed convolution formula for $r_{d_1, d_2}$, since there is typically more than one class per genus. 
For fixed $d_1, d_2$ one can lift the necessary ingredients to analyze \eqref{2} from \cite{Iw} or \cite{Vi}, but for almost primes we need strong uniformity where $d_1, d_2$ can be as large as some small power of $n$.  

A robust substitute that may also be useful in other situations is as follows. Let $F$ be a primitive binary quadratic form of   discriminant $D  = D_0f^2 < 0$ (with $D_0$ fundamental) and let $\mathcal{O}_D \subseteq \Bbb{Q}(\sqrt{D})$ be the (not necessarily maximal) order of discriminant $D$. The set of equivalence classes of primitive forms of discriminant $D$ is in bijection with the proper (i.e.\ invertible) ideal classes $\mathcal{C}_D$ of $\mathcal{O}_D$ (see \cite[Section 7]{Co}), and we denote by $C_F \in \mathcal{C}_D$ the class corresponding to $F$. Then for $(n, f) = 1$  we have the exact formula
\begin{equation}\label{3}
\#\{(x_1, x_2) \in \Bbb{Z}^2 \mid F(x_1, x_2) = n\} = \frac{w}{\# \widehat{\mathcal{C}}_D}\sum_{\chi \in \widehat{\mathcal{C}}_D} \bar{\chi}(C_F) \sum_{\substack{N\mathfrak{a} = n\\ \mathfrak{a} \subseteq \mathcal{O}_K}}\chi(\mathfrak{a})
\end{equation}
where $w$ is the number of units in $\mathcal{O}_D$. In our application, 
the underlying quadratic field is fixed: for all $d_1, d_2$ in \eqref{2} we obtain the Gau{\ss}ian number field $\Bbb{Q}(i)$, but the discriminant of the order increases with $d_1, d_2$, in particular the order is highly non-maximal. This will cause  technical challenges and involved computations of the corresponding main terms.

The key point is now to split the outer sum into two parts: class group characters $\chi$ of order at most two are genus characters and can be computed explicitly by some kind of convolution identity. For a character of order three or more, the inner sum is the $n$-th Hecke eigenvalue of a cusp form of weight 1 and level $|D|$. In this way, the theory of automorphic forms enters, and \eqref{3} can be seen as a version of Siegel's mass formula. 

The upshot is that -- at least morally -- we need to understand a sum as in \eqref{2} in two cases: when the arithmetic function in question is a Fourier coefficients of some cusp form and when the arithmetic function is a Dirichlet convolution of two characters, in other words the Fourier coefficient of an Eisenstein series. In particular, this is yet another example where the originally purely arithmetic problem of representing $n$ in the form \eqref{1}, \eqref{FI-equation} or \eqref{1.2} is  intimately linked to the theory of automorphic forms. Indeed, here we need to invoke the full power of the spectral theory of automorphic forms. 

As an aside we remark that there is a different way to build a bridge to automorphic forms. Instead of applying the algebraic identity \eqref{3}, one could model $r_{d_1, d_2}$ by a suitably truncated  expression coming from a formal application of the circle method. More rigorously, one could apply a $\delta$-symbol method \`a la Heath-Brown \cite{HB} and Duke-Friedlander-Iwaniec \cite{DFI}, followed by Poisson summation. This leads to sum of Kloosterman sums (for congruence subgroups of levels roughly $d_1^2d_2^2$) which after an application of the Kuznetsov formula returns a similar expression of cuspidal and non-cuspidal Fourier coefficients. This can be seen as a double Kloosterman's refinement of the circle method. An application of \eqref{3}, however, seems more direct.

To understand \eqref{2} when $r_{d_1, d_2}$ is replaced by a cuspidal Fourier coefficient, 
we refine methods and results established in  \cite[Theorem 1.3]{ABL}  to obtain an unconditional power saving for \emph{individual} pair $(d_1,d_2)$ with the required degree of uniformity. However, the Eisenstein contribution requires much more (or at least rather different) work. As mentioned above, if we were going to handle each pair $(d_1,d_2)$ \emph{separately} for the Eisenstein contribution, we would ultimately require a zero-free strip for the relevant Dirichlet $L$-functions. To avoid this use of GRH, we instead incorporate the \emph{sum} over $d_1,d_2$ from the sieve weights into the analysis of primes (or smooth numbers) in long arithmetic progressions (see Proposition \ref{TypeII} below). The corresponding generalizations and extensions of \cite[Theorem 2.1]{ABL} are given in  Propositions \ref{PrimeAp} and \ref{SmoothAp} below. This is of independent interest and should find applications elsewhere.

The proof of Theorem \ref{thm2} features additional difficulties. The arithmetic of smooth numbers is rather subtle, and it turns out  that leading term in a suitable asymptotic formula of \eqref{2} for $\ell \in S(n, g(n))$ is not multiplicative  in $d_1, d_2$. Multiplicativity  does not fail by a lot, but we obtain secondary terms involving derivatives of Euler products. This makes the application of a sieve perhaps not impossible, but at least seriously problematic. There are some ways to handle non-multiplicative main terms without developing the theory in full generality (see \cite[p.\ 37]{FI1}). But the fact that we have two variables to sieve makes the monotonicity principle, in particular \cite[Corollary 5.4/5.5]{FI1}, not applicable.  We therefore take a different route and incorporate the sieve weights directly into the analysis. A prototype of this idea can be found in the proof of \cite[Theorem 4]{FT96}, but such a strategy in the context of additive problems seems to be new. There are a number of  interesting but more technical details that we will discuss as they arise.  \\

\textbf{Roadmap for the rest of the paper:} In Section \ref{algebraic}, we give the explicit formula for $r_{d_1,d_2}$ (see Lemma \ref{lem1}), which consists of a cuspidal part and an Eisenstein part. To deal with the cuspidal part, in Section \ref{cuspidal} we prove general bounds for linear forms in Hecke eigenvalues over prime numbers (Theorem \ref{prime}) and smooth numbers  (Theorem \ref{smooth})  assuming a bilinear estimate (Proposition \ref{bilinear}) and a $\tau_3$-type estimate (Proposition \ref{d3}) with Hecke eigenvalues, whose proofs are given in Section \ref{bid3}. We then give a general averaged bilinear estimates in Section \ref{averageBilinear} which will be used to handle the error terms when applying sieve methods. The proofs of Theorems \ref{thm1}/\ref{thm1'} and \ref{thm2} are given in Sections \ref{prime-section} and \ref{smooth-section} respectively. The most involved part is the rather involved analysis of the main term in the sieving process, which covers Subsections \ref{52} -- \ref{53} in the prime number case and  Subsections \ref{part1} -- \ref{part4} in the smooth number case. \\

\textbf{Acknowledgement:} We would like to thank the referee for a careful reading of the manuscript.

\section{Algebraic considerations}\label{algebraic}

We recall the classical theory on binary quadratic forms. For $f \in \Bbb{N}$ let $\mathcal{O}_{f}  = \Bbb{Z} + f\Bbb{Z}[i] \subseteq \Bbb{Z}[i]$ be the  unique order of discriminant $-4f^2$ (the notation is slightly different  than in Section \ref{sec12}). We have $\mathcal{O}_f^{\ast} = \{\pm 1\}$ unless $f = 1$ in which case $\mathcal{O}_f^{\ast} = \{\pm 1, \pm i\}$. We write $w_f = \#\mathcal{O}_f^{\ast}$. Its class group $\mathcal{C}_f$  is the set of  proper (= invertible) fractional $\mathcal{O}_f$-ideals $I_f$ modulo principal ideals $P_f$. This is a finite group of cardinality \cite[Theorem 7.24]{Co}
 \begin{equation}\label{classnumber}
 h_f = f\prod_{p \mid f} \Big(1 - \frac{\chi_{-4}(p)}{p}\Big) \begin{cases} 1/2, & f > 1,\\ 1, & f = 1.\end{cases}
 \end{equation}
 We have a sequence of isomorphisms. Let $\tilde{\mathcal{C}}_f$ be the set of primitive positive binary quadratic forms with discriminant $-4f^2$ modulo the action of ${\rm SL}_2(\Bbb{Z})$. Let $I^*_f$ be the set of fractional $\mathcal{O}_f$-ideals coprime to $f$ (these are automatically proper \cite[Lemma 7.18(ii)]{Co}) and $P^*_f$ the subgroup of its principal ideals. Let $I^{\ast}$ be the set of fractional $\mathcal{O}_1 = \Bbb{Z}[i]$-ideals coprime to $f$ and  $P^{\ast}$ the set of principal $\Bbb{Z}[i]$-ideals with a generator $\alpha \in \mathcal{O}_f$ coprime to $f$. Then we have (\cite[ Theorem 7.7(ii), Exercise 7.17(a), Propositions 7.19,  7.20, 7.22]{Co} with the convention that $\Im (\beta/\alpha)>0$)
\begin{equation}\label{iso}
\begin{array}{ccccccc}
\tilde{\mathcal{C}}_f &\cong & \mathcal{C}_f = I_f/P_f& \cong &I^*_f/P^*_f &\cong& I^{\ast}/P^{\ast}\\
\frac{N(\alpha x - \beta y)}{N\mathfrak{a}} & \leftarrow & \mathfrak{a} = [\alpha, \beta]& \leftarrow& \mathfrak{a}& \rightarrow & \mathfrak{a}\Bbb{Z}[i]\\
 ax^2 + bxy + cy^2 &\rightarrow & [a, \frac{1}{2} (-b + \sqrt{-4f^2})] &&\mathfrak{a} \cap \mathcal{O}_f&\leftarrow& \mathfrak{a}
\end{array}
\end{equation}
The last isomorphism is norm-preserving (\cite[Proposition 7.20]{Co}). 
For a positive integer $m$ and a primitive binary form $F$ with corresponding ideal  class $C_F$, the first isomorphism induces a bijection (\cite[proof of Theorem 7.7(iii)]{Co} or \cite[Theorem 5 in Chapter 2, Section 7.6]{BS}) between
\begin{equation}\label{21}
\{(x_1, x_2) \mid  F(x_1, x_2)  = m \}/\mathcal{O}_f^{\ast} \longleftrightarrow  \{\mathfrak{a} \text{ integral } \mathcal{O}_f\text{-ideal in the class } C_F \mid N\mathfrak{a} = m\}.
\end{equation}
If $(m, f) = 1$, then the last isomorphism in \eqref{iso} shows
\begin{equation}\label{22}
 \sum_{\substack{N\mathfrak{a} = m\\ \mathfrak{a} \subseteq \mathcal{O}_f}} 1 = \sum_{\substack{N\mathfrak{a} = m\\ \mathfrak{a} \subseteq\Bbb{Z}[i]}} 1 = \sum_{d\mid m} \chi_{-4}(d) = r(m)
 \end{equation}
where $r(m) = \textbf{1} \ast \chi_{-4}$ is the usual sums-of-two-squares function.

The elements of order 2 in $\mathcal{C}_f$ are well-understood  by genus theory. In particular, for square-free $f$ the real characters $\chi$ of $\widehat{\mathcal{C}}_f$ are parametrized by odd divisors $D \mid f$,  
and for an ideal $\mathfrak{a} \in I_f^{\ast}$ (i.e.\ coprime to $f$) and a real character associated with the divisor $D$ we have 
\begin{equation}\label{23}
\chi(\mathfrak{a}) = \chi_{D^{\ast}}(N\mathfrak{a}), \quad D^{\ast} := \chi_{-4}(D) D \in \{\pm D\}.
\end{equation}
We write
\begin{equation}\label{genera}
G_{f} = \tau(f/(f, 2))\end{equation}
for the number of genera.

From now on we focus on the form $F_{a,b}=a^2x_1^2+b^2x_2^2$. We write $C_{a, b} \in \mathcal{C}_{ab}$ to denote the class corresponding to the form $F_{a,b}$. We will give an explicit formula for $r_{d_1, d_2}$ defined in \eqref{2}. 
We first notice that we can compute the real characters of $\tilde{\mathcal{C}}_{ab}$ explicitly for the class $C_{a,b}$. 

\begin{lemma}\label{lem21} Let $a, b \in \Bbb{N}$ be coprime and square-free. Let $F_{a, b}$ be as above and let $\chi$ be a  real character of $\tilde{\mathcal{C}}_{ab}$. Then $\chi(C_{a, b}) = 1$. 
\end{lemma}

 \textbf{Proof.} 
 Then  $F_{a, b}$ is equivalent to the form $(a^2 + b^2)x_1^2 + 2b^2x_1x_2 + b^2 x_2^2$, which under the isomorphism \eqref{iso} corresponds to an ideal of norm $a^2 +b^2$ coprime to $ab$, so if $\chi$ belongs to the  divisor $D\mid ab$, then by \eqref{23} we have
 $$\chi(C_{a, b})  = \chi_{D^{\ast}}(a^2+b^2) = 1.$$
 
 \emph{Remark:} Alternatively one can explicitly compute the class group structure from the
exact sequence \cite[(7.25), (7.27), Exercise 7.30]{Co}
 $$1 \longrightarrow \{\pm 1\} \overset{
 }{\longrightarrow} (\Bbb{Z}/f\Bbb{Z})^{\ast} \times \{\pm 1, \pm i\} \longrightarrow (\Bbb{Z}[i]/f\Bbb{Z}[i])^{\ast} \longrightarrow \mathcal{C}_f \longrightarrow 1.$$ 
 For instance, if $f$ is even (and squarefree), one can conclude by a computation based on the Chinese remainder theorem and the fact that the multiplicative group of a finite field is cyclic that $\text{rk}_2(\mathcal{C}_f) = \text{rk}_4(\mathcal{C}_f)$, so that the 2-part of $\mathcal{C}_f$ is a direct product of copies of $\Bbb{Z}/2^{k_j}\Bbb{Z}$ with $k_j \geq 2$ and hence  all real characters are trivial on elements of order 2 (which are squares) and so a fortiori on diagonal forms. \\

From Lemma \ref{lem21} we can give the explicit formula for $r_{d_1,d_2}$ in \eqref{2}. We recall the notation \eqref{classnumber} and \eqref{genera}. 
\begin{lemma}\label{lem1} Let $d_1, d_2$ be two square-free numbers, and write $\delta=(d_1,d_2), d_1' = d_1/\delta$, and $d_2' = d_2/\delta$. Let $\mathcal{L} \subseteq \Bbb{N}$ be a finite set and define
$$\mathcal{S} :=  \sum_{l \in \mathcal{L}} r_{d_1, d_2}(l).$$
Then $\mathcal{S} = \mathcal{S}_{\leq 2}+ \mathcal{S}_{\geq 3}$ where we have 
 \begin{equation}\label{nl2}
 \mathcal{S}_{\leq 2} =\sum_{\delta_1\delta_1' =  d_1'}\sum_{\delta_2 \delta_2'  = d_2'}   \frac{w_{\delta_1' \delta_2'}}{h_{\delta_1'\delta_2'}} G_{\delta_1'\delta_2'}\sum_{w\in \mathcal G_{\delta_1'\delta_2'}} \sum_{\substack{l \equiv 0\, (\text{{\rm mod }} (\delta\delta_1\delta_2)^2)\\ l\equiv w\pmod{\delta_1'\delta_2'}\\ l\in \mathcal{L}} } 
 r\Big(\frac{l}{\delta^2 \delta_1^2 \delta_2^2}\Big).
 \end{equation}
with
\begin{equation}\label{residues}
\mathcal{G}_{\delta_1'\delta_2'} =\big\{w\in (\Bbb{Z}/\delta_1'\delta_2'\Bbb{Z})^{\ast}: \chi_{p^*}(w)=1\text{ for all } 2\nmid p\mid \delta_1'\delta_2'\big\},
\end{equation}
and
\begin{equation}\label{sg3}
\mathcal{S}_{\geq 3} = \sum_{\delta_1\delta_1' =  d_1'}\sum_{\delta_2 \delta_2'  = d_2'}   \frac{w_{\delta_1' \delta_2'}}{h_{\delta_1'\delta_2'}} \sum_{\substack{\chi \in \widehat{\mathcal{C}}_{\delta_1'\delta_2'}\\ \text{{\rm ord}}\chi \geq 3}} \bar{\chi}(C_{\delta_1',\delta'_2}) \sum_{\substack{l \equiv 0\, (\text{{\rm mod }} (\delta\delta_1\delta_2)^2)\\ (l, \delta_1'\delta_2') = 1\\l\in \mathcal{L}}}  \lambda_{\chi}\Big( \frac{l}{\delta^2 \delta_1^2 \delta_2^2}\Big)
\end{equation}
with a Hecke eigenvalue $\lambda_{\chi}$ of some holomorphic cuspidal newform of weight $1$, character $\chi_{-4}$ and level dividing $4d_1^2d_2^2$.

  \end{lemma}

\textbf{Proof.} Since $\delta, d_1', d_2'$ are pairwise coprime and square-free, we have
 \begin{equation}\label{firststep}
 \begin{split}
\mathcal{S} :=  \sum_{l \in \mathcal{L}} r_{d_1, d_2}(l) &= \sum_{\delta^2 \mid l \in \mathcal{L} } r_{d'_1, d'_2}\Big(\frac{l}{\delta^2}\Big) = \sum_{\delta_1\delta_1' =  d_1'}\sum_{\delta_2 \delta_2'  = d_2'} \sum_{\substack{(\delta\delta_1\delta_2)^2 \mid l \in \mathcal{L} \\ (l, \delta_1'\delta_2') = 1}} r_{ \delta'_1, \delta_2'}\Big(\frac{l}{\delta^2\delta_1^2\delta_2^2}\Big).
  \end{split}
 \end{equation}
Note that $(\delta_1', \delta_2') = 1$, so the form $(\delta_1')^2 x_1^2+ (\delta_2')^2 x_2^2$ is primitive and $f = \delta_1'\delta_2'$ is square-free. Moreover, the argument of $r_{ \delta'_1, \delta_2'}$ is (nonzero and) coprime to $\delta_1'\delta_2'$.  
We have by \eqref{21} and orthogonality of characters
$$r_{a, b}(m) = \frac{w_{ab}}{h_{ab}} \sum_{\chi \in \widehat{\mathcal{C}}_{ab}} \bar{\chi}(C_{a,b})\sigma(\chi, m), \quad \sigma(\chi, m) = \sum_{N\mathfrak{a} = m} \chi(\mathfrak{a}).$$
If $\chi$ is real and belongs to the odd divisor $D \mid \delta_1'\delta_2'$, 
then by \eqref{23} and \eqref{22} for $(m, \delta_1'\delta_2') = 1$ we have 
$$\sigma(\chi, m) = \chi_{D^{\ast}}(m)\#\{N\mathfrak{a} = m\} =  \chi_{D^{\ast}}(m)r(m).$$
Combining this with Lemma \ref{lem21}, we obtain
 \begin{equation}
 \mathcal{S}_{\leq 2} =\sum_{\delta_1\delta_1' =  d_1'}\sum_{\delta_2 \delta_2'  = d_2'}   \frac{w_{\delta_1' \delta_2'}}{h_{\delta_1'\delta_2'}} \sum_{2\nmid D \mid \delta_1'\delta_2'}  \sum_{\substack{l \equiv 0\, (\text{{\rm mod }} (\delta\delta_1\delta_2)^2)\\ (l, \delta_1'\delta_2') = 1\\ l\in \mathcal{L}} } \chi_{D^{\ast}}\Big( \frac{l}{\delta^2 \delta_1^2 \delta_2^2}    \Big) r\Big(\frac{l}{\delta^2 \delta_1^2 \delta_2^2}\Big),
 \end{equation}
 which yields \eqref{nl2} after summing over $2\nmid D\mid \delta_1'\delta_2'$.

 If $\text{ord}(\chi) \geq 3$, it remains to show that $\sigma(\chi, m) $ is a normalized Hecke eigenvalue of a weight 1 newform of level dividing $4d_1^2d_2^2$ and character $\chi_{-4}$. 
 This is a special case of automorphic induction for imaginary quadratic fields \cite{AC}, but most of it can be seen elementarily. We observe first that by splitting into ideal classes the function
$$f_{\chi}(z) = \sum_{\mathfrak{a} \subseteq \mathcal{O}_{\delta_1'\delta_2'}} \chi(\mathfrak{a}) e(N\mathfrak{a} \cdot z)$$ is a linear combination of theta series corresponding  to forms $F \in \tilde{\mathcal{C}}_{\delta_1'\delta_2'}$ and hence by \cite[Theorem 10.9]{Iw2} a modular form of weight 1, character $\chi_{-4}$ and level $4(\delta_1'\delta_2')^2$.  Using \eqref{iso}, one can check directly from the definition of the Hecke operators $T_m$ that
$$T_m f_{\chi}(z) = \sigma(\chi, m) f_{\chi}(z)$$
for $(m, \delta_1'\delta_2') = 1$. Hence   $f_{\chi}$ belongs to the space generated by some newform of level dividing  $4(\delta_1'\delta_1')^2$, and for $(m, \delta_1'\delta_2') = 1$ the numbers $\sigma(\chi, m)$ coincide with the Hecke eigenvalues. 

It remains to show that  $f_{\chi}$ is a cusp form for  $\text{ord}(\chi) \geq 3$. This is again a   classical fact (also contained in \cite{AC}) that can be proved in several ways: by the parametrization of Eisenstein series, a modular form is non-cuspidal if and only if the associated Dirichlet series factorizes into two Dirichlet $L$-functions, and this happens if and only if its Rankin-Selberg square $L$-function has a double pole at $s=1$. The Rankin-Selberg square of a class group character $\chi$ is, up to finitely many Euler factors, $\zeta(s) L(s, \chi^2)$ which has a double pole at $s=1$ if and only if $\text{ord}(\chi) \leq 2$. 

This completes the proof of the lemma.  \\

\textbf{Remark:} 
The formula \eqref{nl2} depends on Lemma \ref{lem21}. For general quadratic forms not covered by Lemma \ref{lem21} the summation over $D$ can be carried out in the same way. This yields an expression similar to \eqref{nl2}
with a set $\mathcal G_{\delta_1'\delta_2'}$ which may be different from the one defined in \eqref{residues} but which has the same cardinality, namely $\phi(\delta_1'\delta_2')/G_{\delta_1'\delta_2'}$.
This is a consequence of the fact that numbers coprime to the order can only be represented in one genus. More precisely, 
the Chinese Remainder Theorem induces a bijection $\mathcal{G}_{\delta_1'\delta_2'} \cong \prod_{p \mid \delta_1'\delta_2'} \mathcal{G}_p$ where $\mathcal{G}_p \subseteq (\Bbb{Z}/p\Bbb{Z})^{\ast}$ has cardinality $\phi(p)/2$ for odd $p$ and cardinality 1 for   $p=2$. This is all the information needed in the subsequent arguments. 
In this way, our main results generalize to other quadratic forms than sums of two squares, although some extra care may be  necessary at the last step of \eqref{firststep}. 

If we consider solutions to $n=p+x_1^2+x_2^2$ where  $x_1,x_2$ have no small prime factors other than a small square-free factor of $r$, then we will encounter binary quadratic forms of the shape $x^2+(r^2y)^2$, which can be studied with our method while paying special attention to the real characters of $\mathcal C_f$ with non square-free $f$ and the step in \eqref{firststep}.

\section{Additive problems in cusp forms}\label{cuspidal}
 
In this section we deal with the innermost sum  of \eqref{sg3} when the set $\mathcal{L}$ consists of numbers $n-\ell$ where $\ell$ is either prime or smooth. We generally consider a holomorphic or Maa{\ss} cuspidal Hecke newform $\phi$ of level $4\mid N$, character $\chi_{-4}$ and Hecke eigenvalues $\lambda(n)$ (which by self-adjointness of Hecke operators are real) whose archimedean parameter (weight or spectral parameter) we denote by $\mu$. 
We start with two auxiliary results.

\begin{prop}\label{bilinear} Let $ M, X, Z \geq 1$, $\Delta\in \mathbb{N},$ and  $f, \sigma \in \mathbb{Z} \setminus \{0\}$. Let $\alpha_m, \beta_n$ be sequences supported on $[M, 2M]$, $[X, 2X]$ respectively. Suppose that $|\alpha_m|\ll m^\varepsilon$ for any $\varepsilon>0$ and $\alpha_m\beta_n$ vanishes unless $\sigma mn+f$ is in some dyadic interval $[\Xi, 2\Xi]$ with $|\sigma|MX/Z \ll \Xi \ll |\sigma|MX + |f|$. Then  
		\begin{displaymath}
		\begin{split}
	\underset{ \Delta\mid\sigma mn +f}{\sum\sum}
		\alpha_m \beta_n \lambda\Big(\frac{\sigma mn + f}{\Delta}\Big) 
		&\ll P^A
		  \|\beta\|\big( M^{11/14}X^{1/2}  +    M^{5/4}X^{3/8}     \big)^{1+\varepsilon} 
		\end{split}
		\end{displaymath}
	for some absolute constant $A > 0$ and any $\varepsilon>0$, 	where $P =\Delta|\sigma|(1 + |\mu|)NZ (1 +  |f|/MX)$. 
\end{prop}

We remark that the trivial bound   is $\| \beta \| (MX^{1/2})^{1+\varepsilon}$, so Proposition \ref{bilinear} is non-trivial as soon as $M \ll X^{1/2 - \delta}$ for some $\delta > 0$ as in \cite[Theorem 1.2]{Pi} when $P$ is fixed.

\begin{prop}\label{d3} Let $\Delta, q,d\in \mathbb{N},r, f\in \Bbb{Z} \setminus \{0\}$, $Z, X_1, X_2, X_3 \geq 1/2$. 
		Suppose $(\Delta, q)=1$ and that at least one of $r, f$ is positive. Let $G$ be a function with support on  $\bigtimes [X_j, 2X_j]$ and $G^{(\nu_1, \nu_2, \nu_3)} \ll_{{\bm \nu}} Z^{\nu_1+\nu_2+\nu_3}X_1^{-\nu_1} X_2^{-\nu_2} X_3^{-\nu_3}$ for all $\bm \nu \in \Bbb{N}_0^3$.  Write $X =  |r|X_1X_2X_3$ and suppose that $G(k, l, m)$ vanishes unless $rklm + f$ is in some dyadic interval $[\Xi, 2\Xi]$ with $X/Z \ll \Xi \ll X+|f|$. Then 
		\begin{align}
		&\sum_{\substack{k,l,m\\rklm \equiv -f \pmod \Delta \\ rklm\equiv d \pmod q}} G(k, l, m)\lambda\Big(\frac{rklm + f}{\Delta}\Big) \ll  P^A X^{1-1/54+\varepsilon}
		\end{align}
	for some absolute constant $A$ and any $\varepsilon > 0$, where $P = qr\Delta (1 + |\mu|)NZ(1 + |f|/X)$. 
\end{prop}

The trivial bound is $(X/|r|\Delta q)^{1+\varepsilon}$, and so Proposition \ref{d3} is non-trivial uniformly for $P$ up to some small power of $X$. Better exponents could be obtained with more work, but these bounds suffice for our purpose. In principle this is a shifted convolution problem of $\tau_3$ against $\lambda$, cf.\ e.g.\ \cite{Mu, To}. 

We postpone the proofs of Propositions \ref{bilinear} and \ref{d3} to the next section and state the two main results of this section.


\begin{theorem}\label{prime} There exist  absolute constants $A, \eta > 0$ with the following property. 
	Let $\sigma, f\in \mathbb{ Z}\setminus\{0\}$, $X, Z \geq 1 $, $\Delta, q, d\in \mathbb{ N}$ with $(\sigma q, \Delta)=1$.
	Let $G$ be a smooth function supported on $[X, 2X]$ with $G^{(\nu)}\ll_\nu(X/Z)^{-\nu}$ for $\nu \in \Bbb{N}_0$. Suppose that $G(n)$ vanishes unless $\sigma n+f $ is in some dyadic interval $[\Xi, 2\Xi]$ with $|\sigma |X/Z\ll \Xi \ll|\sigma| X+|f|$. Then
	\begin{align}
	\sum_{\substack{\sigma n\equiv -f \pmod \Delta\\n\equiv d\pmod q}}\lambda\Big(\frac{\sigma n+f}{\Delta}\Big)\Lambda(n)G(n)\ll P^A X^{1-\eta}
		\end{align}
		where $P =  q|\sigma|\Delta(1 + |\mu|)NZ(1+|f|/|\sigma |X) $. 
\end{theorem}

\textbf{Proof.} Based on Propositions \ref{bilinear} and \ref{d3}, this could be proved similarly  as in \cite[Prop.\ 9.1]{ABL} which in turn is the argument of \cite[Section 9]{Pi}. 
The proof in \cite{Pi} uses Vaughan's identity to decompose $\Lambda$ iteratively three times. This is slightly cumbersome, but has the advantage that we only need $\alpha_m$ to be supported on square-free integers in the bilinear estimate in Proposition \ref{bilinear}, which could yield better error terms numerically. Here we give a shorter independent proof at the cost of slightly weaker error terms. Instead of using Vaughan's identity as in \cite{Pi}, we apply Heath-Brown's identity in the form
	\begin{align}\label{HeathBrownId}
	\Lambda(n)=\sum_{j=1}^4 (-1)^{j-1}\binom{j}{4}\sum_{m_1,\dots, m_j\leq X^{1/4}}\prod_{i=1}^j \mu(m_i)\sum_{m_1\cdots m_jn_1\cdots n_j=n}\log n_1
	\end{align}	
for $1 \leq n \leq X$. Next we apply a smooth partition of unity to the variables $m_i, n_i$. To be precise,
let $0<\mathcal Z\leq 1$ be a parameter to be chosen later and let  $\psi(x)$ be a smooth function supported on $[-1-\mathcal Z, 1+\mathcal Z]$ that equals $1$ on $[-1,1]$ and satisfies 
$\psi^{(j)}(t)\ll \mathcal Z^j$. Let $\mathcal D=\{(1+\mathcal Z)^m, m\geq 0\}$. Then we have the smooth partition of unity 
\begin{align}\label{partunity}
1=\sum_{N\in \mathcal D}\psi_N(n)
\end{align}
for any natural number $n$, where $\psi_N(n)=\psi(n/N)-\psi((1+\mathcal Z)n/N)$. 
Here we can choose $\mathcal Z=1$ (later parts of the paper we use the same argument with different choices of $\mathcal Z$) so that
\begin{align}
	\sum_{\substack{\sigma n \equiv -f\pmod {\Delta}\\n\equiv d\pmod q}}\lambda\Big(\frac{\sigma n+f}{\Delta}\Big)\Lambda(n)G(n)\ll  {(\log X)}^{8}\max_{\substack{M_1, \dots, M_j\\ N_1, \dots,N_j}} |S|
\end{align}
where
\begin{align}
S=\sum_{\substack{\sigma \prod_{i=1}^jm_i n_i\equiv -f\pmod \Delta\\ \prod_{i=1}^jm_in_i\equiv d\pmod q}}&\prod_{i=1}^j\mu(m_i)\psi_{M_i}(m_i)\log (n_1)\psi_{N_1}(n_1)\prod_{i=2}^j\psi_{N_i}(n_i)\\
&\times  \lambda\Big(\frac{\sigma \prod_{i=1}^jm_in_i+f}{\Delta}\Big)G\Big(\prod_{i=1}^jm_in_i\Big).
\end{align}	
Therefore it is enough to show that there exist some absolute constants $A,\eta>0$ such that 
\begin{align}\label{primereduction}
S\ll P^A X^{1-\eta}.
\end{align}
	
We distinguish two cases. If there exists some $i$ such that  $X^{\eta_0}\leq M_i \text{ (or $N_i$)}\leq X^{1/4+\eta_0}$ for some $\eta_0>0$, then we can split $m_i$ or $n_i$ into residue classes modulo $q$, use Mellin inversion to separate variables $m_i,n_i$ (see e.g.\ \cite[Section 7]{Pi}) and then apply Proposition \ref{bilinear} to obtain 
	\begin{align}
	S\ll  P^{A}X^\varepsilon \big( {X}^{1-3\eta_0/14}+ X^{31/32+3\eta_0/8}\big).
	\end{align}
	 Otherwise, we must have $$R:=\prod_{\substack{i: M_i\leq X^{\eta_0}, N_i\leq X^{\eta_0}}}M_iN_i\leq X^{8\eta_0}$$ and  $$\#\{1\leq i\leq j:N_i\geq X^{1/4+\eta_0}\}\leq 3, \quad 1\leq j \leq 4,$$ in which case we can apply Proposition \ref{d3} to obtain
	\begin{align}
	S \ll P^AX^{1-1/54+8A\eta_0 +\varepsilon}.
	\end{align} 
	Combining the bounds in these two cases and chossing $\eta_0$ suitably, we complete the proof of \eqref{primereduction}.\\
	
Our second main result treats the case of smooth numbers. As usual, $P^+(n)$ denotes the largest prime factor of $n$ with the convention $P^+(1) = 1$.

\begin{theorem}\label{smooth} 
	There exist   absolute constants $A, \eta>0$ with the following property. 
	Let $\sigma, f\in \mathbb{Z}\setminus\{0\}$, $X, Z, y \geq 1$, $ \Delta, q,d\in \mathbb{N}$ with $(\sigma q, \Delta)=1$. Let $G$ be a smooth function supported on $[X, 2X]$ with $G^{(\nu)}\ll_\nu(X/Z)^{-\nu}$. Suppose that $G(n)$ vanishes unless $\sigma n+f $ is in some dyadic interval $[\Xi, 2\Xi]$ with $|\sigma| X/Z\ll \Xi \ll |\sigma| X+|f|$. Then
	\begin{align}
	\sum_{\substack{\sigma n \equiv -f\pmod \Delta\\ n\equiv d\pmod q\\ P^+(n)\leq y}}\lambda\Big(\frac{\sigma n+f}{\Delta}\Big)G(n)\ll P^AX^{1-\eta}.
	\end{align}
	where $P = q |\sigma|\Delta (1+|\mu|)NZ(1 + |f|/|\sigma| X)$.
\end{theorem}

\textbf{Proof.}  The set-up follows along the lines in \cite[Lemme 2.1]{FT90} on averages of smooth numbers in arithmetic progressions to large moduli, which we also generalize (see Proposition \ref{SmoothAp}).  Throughout the proof, we use $A$ to denote some absolute positive constant, not necessarily the same at each occurrence. For $t \geq 1$ let $g_t(n)$ 
	be the characteristic function on  numbers $n$ with $P^+(n) \leq t$. The starting point is Buchstab's identity that we iterate. For $y \geq 1$ we have 
	\begin{align}\label{SmoothBuchstab}
	g_y(n)=1-\sum_{\substack{p\mid n\\ p>y}}g_p\Big(\frac{n}{p}\Big)=1+\sum_{j=1}^3(-1)^j \sum_{\substack{p_1\cdots p_j\mid n\\y<p_1<\cdots<p_j}}1+\sum_{\substack{p_1\cdots p_4\mid n\\ y<p_1<\cdots <p_4}}g_{p_4}\Big(\frac{n}{p_1
		\cdots p_4}\Big).
	\end{align}
	Thus 
	\begin{align}
	\sum_{\substack{ \sigma n \equiv -f\pmod \Delta \\n\equiv d\pmod q\\ P^+(n)\leq y}}\lambda\Big(\frac{\sigma n+f}{\Delta}\Big)G(n)=\sum_{\substack{ \sigma n \equiv -f\pmod \Delta \\n\equiv d\pmod q}}\lambda\Big(\frac{\sigma n+f}{\Delta}\Big)g_y(n)G(n)=S_0+\sum_{i=1}^3 (-1)^i S_i + \overline{S},
	\end{align}
	say, where 
	\begin{align}
	&S_0=\sum_{\substack{\sigma n \equiv -f\pmod \Delta\\ n\equiv d\pmod q}}\lambda\Big(\frac{\sigma n+f}{\Delta}\Big)G(n),\\
	&S_i=\sum_{\substack{ \sigma n \equiv -f\pmod \Delta\\ y<p_1< \cdots < p_i\\ p_1\cdots p_im\equiv d\pmod q}}\lambda\Big(\frac{\sigma p_1\cdots p_i m+ f}{\Delta}\Big)G(p_1\cdots p_i m),\\
	& \overline{S}=\sum_{\substack{\sigma n \equiv -f\pmod \Delta \\ y<p_1< \cdots < p_4\\ p_1\cdots p_4m\equiv d\pmod q}}\lambda\Big(\frac{\sigma p_1\cdots p_4 m+f}{\Delta}\Big)g_{p_4}(m)G(p_1\cdots p_4 m).
	\end{align}
	For $S_0$, we apply \cite[Corollary 7.6]{ABL} to obtain 
	\begin{align}\label{S0bound}
	S_0=\sum_{  k \equiv \bar{ \Delta}(f+\sigma d)\, (\text{mod } q) }\lambda(k)G\Big(\frac{\Delta k-f}{\sigma }\Big)\ll P^A.
	\end{align}
For $\overline{S}$, we localize $p_i$ in intervals $(P_j, P_{j+1}]$ where $P_j=y(1+\mathcal {Z})^j$ with some $\mathcal {Z}\leq 1$ to be chosen later and $0\leq j\leq J= 1+\lfloor \frac{\log (X/y)}{\log (1+\mathcal {Z})}\rfloor$. Then 
	\begin{align}\label{Spartition}
	\overline{S}&=\sum_{\substack{0\leq i_1\leq i_2\leq i_3\leq i_4\leq J\\ p_k\in(P_{i_k},P_{i_k+1}]\\ \sigma p_1\cdots p_4m \equiv -f\pmod \Delta \\ p_1\cdots p_4m\equiv d\pmod q}}\lambda\Big(\frac{\sigma p_1\cdots p_4 m+f}{\Delta}\Big)g_{p_4}(m)G(p_1\cdots p_4m)
	 +O\Big(P^A X^{1+\varepsilon} \mathcal{Z}^{1/2}\Big). 
	\end{align}
	To justify this,
	we see that $\overline{S}$ and the sum in \eqref{Spartition} differ only by integers $n$ that    lie in $[X(1+\mathcal {Z})^{-4},X(1+\mathcal {Z})^4]$ or have at least two prime factors $p_i,p_j$ with $p_i<p_j<(1+\mathcal {Z})p_i$.
	The contribution from these integers can be bounded by
	\begin{align}
	&\Big(\sum_{\substack{X(1+\mathcal {Z})^{-4}<n\leq X(1+\mathcal {Z})^4}}+\sum_{\substack{p_i<p_j<(1+\mathcal {Z})p_i\\ p_ip_j\mid n\\ n \ll X}}\Big)\Big|\lambda\Big(\frac{\sigma n+f}{\Delta}\Big)\Big|\\
	&\ll \Bigg(\Big(\sum_{\substack{X/(1+\mathcal {Z})^{-4}< n\leq X(1+\mathcal {Z})^4}}1 \Big)^{1/2}+\Big( \sum_{\substack{n\ll X\\ p_ip_j\mid n\\ p_i<p_j\leq (1+\mathcal {Z})p_i}} 1\Big)^{1/2}\Bigg)\Big(\sum_{\substack{n\ll X\\ \sigma n \equiv -f\pmod \Delta}}\Big|\lambda\Big(\frac{\sigma n+f}{\Delta}\Big)\Big|^2\Big)^{1/2}\\
	&\ll P^A X^{1+\varepsilon} \mathcal{Z}^{1/2}. 
	\end{align}
	Since $y<P_{i_1}\ll X^{1/4}$, we can apply Proposition \ref{bilinear} with $\bold m=p_1$ and $\bold n=|\sigma |p_2\cdots p_4 m$ to the main term in \eqref{Spartition} by splitting $\bold m,\bold n$ into residue classes modulo $q$ to obtain 
	\begin{equation}\label{bars}
	\overline{S}\ll \mathcal {Z}^{-4} X^{\varepsilon} P^A  \Big( \frac{X}{y^{3/14}}+ X^{31/32}\Big)+P^A X^{1+\varepsilon} \mathcal{Z}^{1/2}  \ll P^A X^{1+\varepsilon}    \max(y^{-1/42}, X^{-1/288})
	\end{equation}
	upon  choosing $\mathcal{Z} = \max(y^{-1/21}, X^{-1/144})$.	 

For $S_k$, $1 \leq k \leq 3$, we first replace the characteristic function on primes with the von Mangoldt function $\Lambda$. As above, we see that the contribution of higher prime powers is at most
$$\ll P^A X^{1+\varepsilon} y^{-1/2}$$
which can be absorbed in the existing bounds. Now similar to the proof of Theorem \ref{primereduction}, we use Heath-Brown's identity \eqref{HeathBrownId} to decompose the prime variables. Again we split all variables into short intervals using the smooth partition \eqref{partunity} for some parameter $\mathcal{Z}$ to be chosen later. We first bound 
	\begin{align}
	\begin{split}
	\tilde{S}=\sum_{m}\sum_{m_i } &\sum_{\substack{n_i\\ \sigma m m_1 \cdots m_r n_1 \cdots n_s \equiv -f \pmod \Delta \\ mn_1\cdots m_rn_1\cdots n_s\equiv d\pmod q}} \psi_M(m)\prod_{i=1}^r \mu(m_i)\psi_{M_i}(m_i) \prod_{j=1}^s\psi_N(n_j)\\
	&\lambda\Big(\frac{\sigma m m_1\cdots m_rn_1\cdots n_s+f}{\Delta}\Big)G(mm_1\cdots m_r n_1 \cdots n_s)
	\end{split}
	\end{align}
	where $1\leq r, s \leq 12$, $ M\prod_{i=1}^r M_i \prod_{j=1}^s N_j\asymp X$ and $M_i\ll X^{1/4}$. Again we distinguish two cases.  
	
	If there exists a subset $\mathcal I \subset \{M,M_i,N_j\}$ such that  $K:=\prod_{I\in \mathcal I}I$ such that $X^{\eta_0}< K \leq X^{1/4+\eta_0}$, then we can apply Proposition \ref{bilinear} to obtain 
	\begin{align}
	\tilde S\ll P^A X^{\varepsilon} \big( {X^{1-3\eta_0/14 }} + X^{31/32+3\eta_0/8}\big).
	\end{align}
	If such $\mathcal I$ does not exist, then  all the $M_i$'s and possibly some of the $M$, $N_j$'s can be combined into $R\ll X^{\eta_0}$ and the number of $M$,$N_j$'s of size $\gg X^{1/4+\eta_0}$ is less than, or equal to three. So it is enough to bound 
	\begin{align}
	\sum_{r\ll R|\sigma|} \Big|\sum_{\substack{k,l,m\\ rklm \equiv -f \pmod \Delta \\ rklm \equiv d\pmod q }}\alpha_r\lambda\Big(\frac{rklm+f}{\Delta}\Big)\psi_K(k)\psi_L(l)\psi_M(m) G(rklm)\Big|
	\end{align}
	which by Proposition \ref{d3} is 
	\begin{align}
	&\ll P^A  X^{1-1/54+A\eta_0 +\varepsilon}
		\end{align}
for some absolute constant $A > 0$. Then by the same argument as above we see that
$$S_i \ll P^A X^{\varepsilon} (\mathcal{Z}^{-24} X^{1-1/54+A\eta_0 +\varepsilon} + X\mathcal{Z}^{1/2}), \quad 1\leq i\leq 3$$
and hence by choosing $\mathcal{Z}=X^{\eta_1}$ for some suitable $\eta_1>0$ we obtain
\begin {equation}\label{si}
S_i \ll  P^A X^{1 - \eta_2 }, \quad 1\leq i\leq 3
\end{equation}
for some $ \eta_2 > 0$. Combining \eqref{S0bound}, \eqref{bars} and \eqref{si}, we have shown that there exist some absolute constants $A, \eta>0$ such that 
 \begin{align}\label{smoothylarge}
 \sum_{\substack{ \sigma n \equiv -f\pmod \Delta \\ n\equiv d\pmod q\\ P^+(n)\leq y}}\lambda\Big(\frac{\sigma n+f}{\Delta}\Big)G(n)\ll P^AX^\varepsilon\big(X^{1-\eta}+X y^{-\eta}\big).
 \end{align}

	On the other hand, when $y$ is very small we can do better using the flexible factorization of smooth numbers to create some bilinear structure so that Proposition \ref{bilinear} can be applied directly. Recall \cite[Lemme 3.1]{FT96}: for any $M \geq 1$, every $n\geq M$ with $P^+(n)\leq y$ has a unique representation in the form 
	\begin{align}\label{smoothfactor}
	n=lm, \quad P^+(l)\leq P^-(m), \quad M\leq m\leq yM.
	\end{align}
	We can separate $l$ and $m$ in the condition $P^+(l)\leq P^-(m)$ using 
	\cite[eq. (3.37)]{Dr1} with an acceptable error. Thus by the same argument as before, we can localize $m$   into short intervals and apply Proposition \ref{bilinear} to obtain 
	\begin{align}\label{smoothysmall}
	\sum_{\substack{ \sigma n \equiv -f\pmod \Delta \\ n\equiv d\pmod q\\ P^+(n)\leq y}}\lambda\Big(\frac{\sigma n+f}{\Delta}\Big)G(n)&\ll P^AX^\varepsilon\big(\mathcal {Z}^{-1}  (XM^{-3/14}+(M y)^{3/8}X^{7/8})+ X\mathcal Z^{1/2}\big)\\
	&\ll P^A X^{1-1/66+\varepsilon}y^{1/22}
	\end{align}
	by choosing $M=X^{7/33}y^{-7/11}, \mathcal {Z}=X^{1/33}y^{-1/11}$.
	Applying \eqref{smoothylarge} when $y\geq X^{\eta}$ and \eqref{smoothysmall} when $y\leq X^{\eta}$ for some suitable $\eta>0$ we obtain the theorem.\\
	
 \textbf{Remark.} Arguing similarly as in the proofs of Theorems \ref{prime} and \ref{smooth},  we can replace $\Lambda(n)$ or $\mathds{1}_{P^+(n)\leq y}$ by $\tau_k(n)$ (for fixed $k$) and obtain a power saving bound \emph{uniformly} in all parameters.  Moreover, the proof works also for a wide class of arithmetic functions which possess a similar combinatorial decomposition (cf. e.g. \cite{DT, FT22} as well as \cite[section 5.1]{MT}), such as generalized divisor functions $\tau_z(n)$ for any complex $z$ and the indicator function of norm forms of abelian extensions of $\mathbb{Q}$.  
 \\	
 

 \section{Proofs of Propositions \ref{bilinear} and \ref{d3}}\label{bid3}
 
 We still owe the proofs of  Propositions \ref{bilinear} and \ref{d3}. Recall that $P$ is the product of the ``unimportant'' parameters which has slightly different meanings in Proposition \ref{bilinear} and Proposition \ref{d3}.  For notational simplicity we write
 \begin{equation}\label{curly}
 	R \preccurlyeq S \quad : \Longleftrightarrow \quad R \ll X^{\varepsilon} P^A S
 \end{equation}
 for any $\varepsilon > 0$ and some   $A > 0$, not necessarily the same at each occurrence.

  \subsection{Proof of Proposition \ref{bilinear}} This is an extension of \cite[Proposition 8.1]{ABL} where we relax the condition $\alpha_m$ supported on square-free numbers at the cost of slightly weaker bounds. Let $\mathscr{B}(\alpha, \beta)$ denote the $m, n$-sum we want to bound. We use the notation \eqref{curly} with $P=\Delta|\sigma|(1+|\mu|)N Z(1+|f|/MX)$. Note that we can assume without loss of generality that $\sigma=\{\pm 1\}$ as otherwise we can set $\tilde{\beta}_n=\beta_n\mathds{1}_{\sigma \mid n}$. 
  When $\alpha_m$ is supported on square-free numbers, then a straightforward extension of \cite[Proposition 8.1]{ABL} (replacing the $F$ in the definition of set of moduli $\mathcal Q$ by $[N_1, D^2]\Delta_1$, with the notation $\Delta=\Delta_1\Delta_2$ and $\Delta_2=(\Delta, (m_3m_4)^\infty)$ to accommodate the new parameter $\Delta$) shows
  \begin{equation}\label{SquarefreeBilinear}
  \mathscr{B}(\alpha, \beta) \preccurlyeq  \| \beta\| \big((XM)^{1/2} + X^{1/4} M^{3/2} + X^{3/8} M^{5/4}\big). 
  \end{equation}
  To deduce the result for general sequences, we write
  \begin{align*}
  \mathscr{B}(\alpha, \beta)=  \underset{\Delta \mid \sigma l^2mn + f}{\sum\sum} \mu^2(m)\alpha_{ml^2}\beta_n\lambda\Big(\frac{\sigma ml^2n+f}{\Delta}\Big).
  \end{align*}
  We choose a parameter $L > 1$. By H\"older's inequality, we see that the contribution from $l\geq L$ gives
  \begin{align*}
  &\ll\Big(\sum_{L\leq l\ll M^{1/2}}\sum_{m\ll M/l^2}\sum_{n\asymp X}|\alpha_{ml^2}\beta_n|^{4/3}\Big)^{3/4} \Big(\sum_{L\leq l\ll M^{1/2}} \sum_{m\ll M/l^2}\sum_{n\asymp X}|\lambda\Big(\frac{\sigma ml^2n+f}{\Delta}\Big)|^4\Big)^{1/4}\\
  &\preccurlyeq  \Big (\sum_{L\leq l\ll M^{1/2}}\sum_{m\ll M/l^2} 1\Big)^{3/4} \Big(\sum_{n}|\beta_n|^2\Big)^{1/2}\Big(\sum_{n \asymp X}1\Big)^{1/4} (MX)^{1/4} \ll  \|\beta \| \frac{MX^{1/2}}{L^{3/4}}
  \end{align*}
  using the well-known bound for the fourth moment of Hecke eigenvalues
  \begin{equation}\label{L4bound}
  \sum_{n \leq Y} |\lambda(n)|^4  \leq Y^{1+\varepsilon} \sum_{n \leq Y} \frac{|\lambda(n)|^4 }{n^{1+\varepsilon}} \preccurlyeq Y.
  \end{equation}
 (Indeed,  the Dirichlet series $\sum_n \lambda_\phi(n)^4 n^{-s}$ equals  $L(s, \text{sym}^2 \phi \times \text{sym}^2 \phi) L(s, \text{sym}^2\phi)^2 \zeta(s)$ up to an Euler product that is absolutely convergent in $\Re s > 1/2 + 2\theta$ where $\theta \leq 7/64$ is an admissible exponent for the Ramanujan-Petersson conjecture. Hence if $C_\phi$ denotes the conductor of $\phi$, then we conclude from   standard bounds for $L$-functions at $s=1$ \cite{XLi}  the second bound in the above previous display.)

  On the other hand, for $l\leq L$ we combine $n$ and $l^2$ to one variable and apply \eqref{SquarefreeBilinear} with $\tilde{\beta}^{(l)}_n=\beta_{n/l^2}\mathds{1}_{l^2\mid n}$ to obtain
  \begin{align*}
  &\sum_{l\leq L}\sum_{m\asymp M/l^2}\sum_{n\asymp Xl^2}\mu^2(m)\alpha_{ml^2}\tilde{\beta}_{n}^{(l)}\lambda\Big(\frac{\sigma mn+f}{\Delta}\Big)\\
  &\preccurlyeq\sum_{l\leq L} \|\tilde{\beta}^{(l)}\|  \Big((MX)^{1/2}+ (Xl^2)^{3/4} (M/l^2)^{3/2}+ (M/l^2)^{5/4} (Xl^2)^{3/8}\Big)  \\
  &\ll  \|\beta\|\big( L  (MX)^{1/2} +X^{3/4} M^{3/2} +M^{5/4}X^{3/8}\big). 
  \end{align*}
  With $L = M^{2/7}$, we obtain
  $$\mathscr{B}(\alpha, \beta) \preccurlyeq \| \beta \| (X^{1/2} M^{11/14}  +X^{1/4} M^{3/2} +M^{5/4}X^{3/8}).$$
  Here we can drop the middle term, because if it dominates the last term, then $M > X^{1/2}$ in which case the claim is trivially true. This completes the proof.

 \subsection{Proof of Proposition \ref{d3}} This result has a precursor in  \cite[Theorem 2.5]{ABL}, but it turns out that here the extra congruence conditions modulo $q$ and $\Delta$  are somewhat subtle and not completely  straightforward to implement due to well-known difficulties in the Voronoi summation formula for general moduli. We therefore modify the proof, the biggest difference being that we replace the delta-symbol method (\cite[Lemma 7.3]{ABL}) of Duke-Friedlander-Iwaniec with Jutila's circle method.\footnote{As an aside we remark that we could also follow Munshi's strategy \cite{Mu}, but note that his Lemma 6 needs to be corrected by a factor $q_1$ when $q_1 = \tilde{q}_1$. }
 
 Let us denote by $\mathscr{S}$ the $k, l, m$-sum that we want to bound and recall that $P$ is the product of the ``unimportant'' parameters (keeping in mind that the definition of $P$ differs slightly from the previous proof).  Assume without loss of generality $X_1 \geq X_2 \geq X_3$. 

 By Voronoi summation (cf.\ \cite[Corollary 7.7]{ABL}) it is easy to see that
 \begin{equation}\label{bound1}
 \mathscr{S} \preccurlyeq (X/X_1)^{3/2}.
 \end{equation}
We can also write 
\begin{align*}
\mathscr S=\sum_{\substack{k,l,m\\ rklm\equiv d\, (\text{mod } q)}}\sum_{n}G(k,l,m)g(n)\lambda(n)\int_0^1 e\big(\alpha( rklm+f-\Delta n)\big) d\alpha.
\end{align*}
where the redundant function $g(n)$  has support on $\frac{1}{2}\Xi/\Delta\leq n\leq 3\Xi/\Delta$ and satisfies $g(n)=1$ for $\Xi/\Delta\leq n\leq 2\Xi/\Delta$ and $g^{(\nu)}(n)\ll (\Xi/\Delta)^{-\nu}$ for all $\nu\in \mathbb N_0$.
We now employ Jutila's circle method as in the corollary to \cite[Lemma 1]{Ju}, i.e.\ we insert the constant function into the $\alpha$-integral and approximate it by a step function of small rational intervals.  This device will simplify the argument especially with respect to the extra congruence condition on $q$. Let $C\gg X^{1/2}$ be a large parameter and $\mathcal{C}\subset [C, 2C]$ be the set of moduli which are coprime to $\Delta Nq$. Define $$L=\sum_{c\in\mathcal C}\phi(c)\gg C^{2}(\Delta NqC)^{-\varepsilon}$$ and let $\delta$ be a real number such that 
   $C^{-2}\ll \delta \ll C^{-1}.$ 
    Then we have (bounding the error term with Cauchy-Schwarz and Rankin-Selberg)
\begin{align*}
\mathscr S=\tilde{\mathscr S}+ O\Big( \frac{(X\Xi/\Delta)^{1/2}}{\delta^{1/2} C} (XC)^{\varepsilon} \Big)
\end{align*}
where
\begin{align*}
\tilde{\mathscr S}=\sum_{\substack{k,l,m\\ rklm\equiv d\, (\text{mod } q)}}\sum_{n}G(k,l,m)g(n)\lambda(n)\frac{1}{L}\sum_{c\in \mathcal{ C}}  \underset{a\, (\text{mod } c)}{\left.\sum\right.^{\ast}} \frac{1}{2\delta}\int_{-\delta}^\delta e\Big( \Big(\frac{a}{c}+\alpha\Big)(rklm+f-\Delta n)\Big)d\alpha.
\end{align*}
We now apply Voronoi summation \cite[Theorem A.4]{KMV} to the $n$-sum and note that by construction $(\Delta N, c) = 1$ (this flexibility is the main advantage of the present set-up). We obtain 
$$\sum_{n } g(n)e(\alpha \Delta n) \lambda(n)e\Big(\frac{a\Delta}{c}n\Big) = \frac{\xi}{c\sqrt{N}} \sum_{\pm} \chi_{-4}(\mp c) \sum_{n } \lambda(n)e\Big(\mp \frac{\overline{a\Delta N}}{c}n\Big) \mathcal{G}^{\pm}\Big(\frac{n}{c^2 N}\Big)$$
for $(a, c) = 1$ ,  some constant $\xi \in S^1$ depending only  on the cusp form $\phi$ and  functions $\mathcal{G}^{\pm}$
 which by \cite[Lemma 7.5]{ABL} satisfy uniformly in $|\alpha|\leq \delta$ the bound
  \begin{align*}
y^\nu\frac{ d^\nu\mathcal{G^{\pm}}(y)}{dy^\nu}\ll_{A, \nu} \frac{\Xi}{\Delta} \Big(1+\frac{1}{( y\Xi/\Delta)^{2\tau+\varepsilon}}\Big)\Big((1+|\mu|)\Big(1+\sqrt{ \frac{y\Xi}{\Delta}}\Big)\Big)^\nu\Big(1+\frac{y\Xi/\Delta}{(\delta\Xi/\Delta+1)^2}+\frac{y\Xi/\Delta}{(1+|\mu|)^2}\Big)^{-A}
 \end{align*}
where $\tau = 0$ if $\phi$ is holomorphic and $\tau = \Im t_{\phi}$ if $\phi$ is Maa{\ss} with spectral parameter $t_{\phi}$.  
We choose
 $$\delta = 1/X.$$
 With this choice we can truncate the $n$-sum at $\preccurlyeq C^2/X$ and obtain
 \begin{multline*}
   \mathscr{S} \preccurlyeq
    \frac{X^{3/2}}{C }+
     \frac{X}{L} \sup_{|\alpha|\leq \delta} \sum_{\pm} \sum_{c \in \mathcal{C}}  \sum_{n \preccurlyeq  C^2/X} \frac{|\lambda(n)|}{c} \Big(1 +  \frac{C^2}{nX}\Big)^{2\tau} 
   \\\Big| \sum_{k, l, m} G(k, l, m)e\Big(\Big(\alpha + \frac{b}{q}\Big)rklm\Big)S(rklm + f,  \pm \overline{\Delta} n, c)\Big| .
 \end{multline*}
 The factor $e(\alpha rklm)$ can be incorporated into $G$ at essentially no cost, since $\alpha rklm \ll \delta X = 1$. When  $X_1\ll C$, we can bound the $k, l, m$-sum using \cite[Lemma 7.2]{ABL} to obtain
 $$\preccurlyeq \frac{X(c, f, n)^{1/2}}{c^{1/2}} + \Big((c, f)^{1/2} c^{7/4} + X_1^{1/2} c^{3/2} + (X_2X_3)^{1/2} c^{5/4} + X_1(X_2X_3)^{1/2}c^{1/2}\Big)c_{\square}^{1/4}$$
 where $c_{\square}$ is the squarefull part of $c$.  We can now evaluate the sums over $n$ (cf.\ \cite[(7.1)]{ABL}) and $c$ (with Rankin's trick, cf.\ two displays after \cite[(8.3)]{ABL}) getting
  \begin{equation*}
 \begin{split}
   \mathscr{S} & \preccurlyeq   \frac{X }{C^{1/2}} + C^{7/4} + X_1^{1/2} C^{3/2} +X^{1/3} C^{5/4} + X_1^{1/2}X^{1/2}C^{1/2}     +\frac{X^{3/2}}{C }\\
   \end{split}
 \end{equation*}
(recall that $X_1 \geq X_2 \geq X_3$).   At this point we make the admissible choice $C = X^{1/2 + 1/54}$ getting 
  \begin{equation}\label{bound2}
 \begin{split}
   \mathscr{S}  \preccurlyeq   X^{53/54} + X_1^{1/2} X^{7/9}. 
   \end{split}
 \end{equation}
 Then the claim follows from \eqref{bound2} when $X^{1/3} \leq X_1 \leq X^{2/5}$ and from \eqref{bound1} when $X_1 \geq X^{2/5}$. 
 
\section{Averaged bilinear estimates}\label{averageBilinear}
The main result of this section is  an averaged bilinear estimate which will be used later to handle the error terms when applying sieve methods. Define
\begin{align}\label{urqd}
	\mathfrak u_R(n; q,a,d,w)&=\mathds{1}_{\substack{n\equiv a\pmod q\\ n\equiv w\pmod d}}-\frac{1}{\phi(q)\phi(d)}\sum_{\substack{ \chi \pmod q , \psi \pmod d\\ \operatorname{ cond}(\chi\psi)\leq R}}\chi(n\bar{a})\psi(n\overline{w}).
	\end{align}
Note that $\mathfrak u_R(n; q, a, 1, \ast)$ equals $\mathfrak u_R(n\bar{a}, q)$ as defined in \cite[Section 4]{ABL}. Here we only use $\psi \pmod d$ with $\operatorname{cond}(\psi)\leq R$ so that GRH can be avoided when estimating contributions from the second term in \eqref{urqd}. However, we could not treat the contribution from each $d$ \emph{individually} as in \cite[Proposition 4.1]{ABL}  and thus we make use of the \emph{sum} over $d$ as the following result indicates. 

\begin{prop}\label{TypeII}
	Let $M, N, Q, R\geq 1$, $a_1 \in \mathbb{Z}\backslash \{0\}$, write  $x=MN$. Let $c\in \Bbb{N}$, $c_0 \in \Bbb{Z}$ with $(c_0, c) = 1$. Let $\alpha_m$ and $\beta_n$ be two sequences supported in $m\in (M, 2M]$ and $n\in (N, 2N]$ such that for some $A \geq 1$ we have 
	$\alpha_m\leq \tau(m)^A, \beta_n\leq \tau(n)^A$,  
	and suppose $\lambda_d\ll \tau(d)^A$.

	Let $\eta > 0$ be any sufficiently small number. Then there exist $\delta = \delta(\eta) > 0$ and $D = D(\eta, A)$ with the following property.
	If
	\begin{equation*}
	x^\eta\leq N \leq x^{1/4+\eta},  \quad 
	Q\leq x^{1/2+\delta}, \quad 
	c, R\leq x^\delta, \quad 
	|a_1|\leq x^{1+\delta}, \quad |a_2|\leq x^\delta
	\end{equation*}
	then	
	$$ \sum_{d\leq x^\delta}\lambda_d\sup_{\substack{w\pmod d^\times\\a_3\mid d}}	\Big|\sum_{\substack{Q\leq q\leq 2Q\\ (q, a_1a_2d)=1\\ q\equiv c_0 \,(\text{{\rm mod }}  c)}}\sum_{\substack{m, n\\(n, a_2)=1}}\alpha_m\beta_{n}u_R(mn;q,a_1\overline{a_2a_3},d,w) \Big|\ll_{ \eta, A} cx(\log x)^{D}R^{-1/3}.$$
	\end{prop}
We remark that the condition $(n,a_2)=1$ is important in this proposition for a few reasons: it is used by \cite[around eq.\ (5.11)] {Dr2} to reduce to the case of $n$ being square-free as well as the modulo $1$ reduction of the exponential phase in \cite[(5.24)]{Dr2}. More importantly, it is required in the application of the Kuznetsov formula when estimating sum of Kloosterman sums (e.g. \cite[the condition $(q,s)=1$ in Proposition 4.12]{Dr2}).\\

\textbf{Proof.}  
Following the proof of \cite[Proposition 4.1]{ABL}, it is enough to prove for any smooth function $\gamma:\mathbb{R}_+\rightarrow [0,1]$ with 
$$\mathds{1}_{Q\leq q\leq 2Q}\leq \gamma(q)\leq \mathds{1}_{Q/2\leq q\leq 3Q/2}, \quad \|\gamma^{(j)}\|_\infty\ll_j Q^{-j+B\delta j} \text{ for some $B\geq 0$ and for all  $j\geq 0$}$$
that under the assumptions in Proposition \ref{TypeII} as well the additional condition that $\beta_n$ is supported on square-free integers we have  
\begin{equation}\label{bilinearsum}
\mathcal S:=\sum_{d\asymp D}\lambda_d\sup_{\substack{w\pmod d^\times\\a_3\mid d}}\Big|\sum_{\substack{q\leq Q\\ (q,a_1a_2d)=1\\q\equiv c_0 \pmod c}}\gamma(q)\operatornamewithlimits{\sum\sum}_{\substack{ m,n \\ (n,a_2)=1}}\alpha_m \beta_n \mathfrak u_R(mn;q,a_1\overline{a_2a_3},d,w) \Big|\ll \frac{cx(\log x)^{O(1)}}{R^{1/3}}
\end{equation}
for $D\ll x^\delta$. 
The proof of \eqref{bilinearsum} is similar to that in \cite[Proposition 4.2]{ABL} with some modifications. Here we have a sum over $d$ which will only influence the application of the large sieve inequalities, as we obtain power saving error terms in all other parts of the argument for each $d$. Note that we also have an additional condition $(q,d)=1$ compared to that in \cite[Proposition 4.2]{ABL}, which can be incorporated following the proof. The sup over $w\pmod d$ and $a_3\mid d$ cause no issue as the main expressions will be independent of $w, a_3$.  

To start with, we also assume $\beta_n$ is supported on $n\equiv b_0\pmod c$ as in  the second display of the proof of Proposition 4.2 of \cite{ABL}, which makes certain coprime conditions easier to verify when applying \cite[Theorem 2.1]{ABL} (and costs a factor $c$ in the final bound).  

We now apply the Cauchy-Schwarz inequality in $d, m$ with the bounds for $\lambda_d$ and $\alpha_m$.  After majorizing the summation over $m$ with the help of a suitable smooth weight $\alpha(m)$ we get
\begin{align}
\mathcal S&\ll(\log x)^{O(1)}(D M)^{1/2}\Big(\sum_{d\asymp D}\sup_{\substack{w\pmod d^\times\\a_3\mid d}}\sum_{\substack{m }}\alpha(m)\Big|\sum_{\substack{q\\(q,a_1a_2d)=1 \\ q\equiv c_0 \pmod c}}\gamma(q)\sum_{\substack{n \\ (n,a_2)=1}}\beta_n \mathfrak u_R(mn;q,a_1\overline{a_2a_3},d,w)\Big|^2\Big)^{1/2}\\
&=:(\log x)^{O(1)}(DM)^{1/2}\Big(\sum_{d\asymp D}\sup_{\substack{w\pmod d^\times\\a_3\mid d}} \mathcal S_1(w,a_3)-2\Re(\mathcal S_2(w,a_3))+\mathcal S_3(w,a_3)\Big)^{1/2},
\end{align}
where (if $(n,qd)\neq 1$ the above sum is empty)
\begin{align}
&\mathcal S_1(w,a_3)=\sum_{\substack{ (q_i, a_1a_2d)=1\\ q_i \equiv c_0\pmod c}}\gamma(q_1)\gamma(q_2)\sum_{\substack{n_i\asymp N\\ (n_i, q_ia_2d)=1}}\beta_{n_1}\overline{\beta_{n_2}}\sum_{\substack{  mn_i\equiv w\pmod d\\ mn_i \equiv a_1 \overline{a_2a_3}\pmod {q_i}}}\alpha(m),\\
& \mathcal S_2(w,a_3)=\sum_{\substack{ (q_i, a_1a_2d)=1\\ q_i \equiv c_0 \pmod c}}\frac{\gamma(q_1)\gamma(q_2)}{\phi(q_2)\phi(d)}\sum_{\substack{n_i\asymp N\\ (n_i, q_ia_2d)=1}}\beta_{n_1}\overline{\beta_{n_2}}\\
& \quad \quad \quad \quad \times\sum_{\substack{\chi_2 \pmod {q_2}\\ \psi_2\pmod {d}\\\operatorname{cond}(\chi_2\psi_2)\leq R }}\chi_2(n_2\overline{a_1}a_2a_3)\psi_2(\bar{w})\sum_{\substack{ mn_1\equiv a_1\overline{a_2a_3}\pmod {q_1}\\ mn_1\equiv w\pmod d}}\chi_2\psi_2(m)\alpha(m),
\intertext{and}
& \mathcal S_3(w,a_3)=\sum_{\substack{ (q_i, a_1a_2d)=1\\ q_i \equiv c_0 \pmod c}}\frac{\gamma(q_1)\gamma(q_2)}{\phi(q_1)\phi(q_2)\phi(d)^2}\sum_{\substack{n_i\asymp N\\ (n_i, q_ia_2d)=1}}\beta_{n_1}\overline{\beta_{n_2}}\\
& \quad \quad\quad \quad  \times \sum_{\substack{\chi_i \pmod {q_i}\\\psi_i\pmod d\\\operatorname{cond}(\chi_i\psi_i)\leq R }}\chi_1(n_1\overline{a_1}a_2a_3)\overline{\chi_2(n_2\overline{a_1}a_2a_3)}\psi_1\overline{\psi_2}(\bar{w})\sum_{\substack{m}}\chi_1\psi_1\overline{\chi_2\psi_2}(m)\alpha(m).
\end{align} 
It is then enough to show  (recall $MN = x$)
\begin{align}\label{dispersion}
\sum_{d\asymp D}\sup_{\substack{w\pmod d^\times\\a_3\mid d}}\Big(\mathcal S_1(w,a_3)-2\Re\mathcal S_2(w,a_3)+\mathcal S_3(w,a_3)\Big)\ll MN^2 D^{-1}R^{-2/3} (\log x)^{O(1)}.
\end{align}
We can evaluate $\mathcal{S}_3$ as in \cite[Section 5.3.1]{Dr2} so that the $m$-sum becomes 
\begin{align}
\frac{\hat{\alpha}(0)}{[q_1,q_2]d}\sum_{h\pmod W^\times}\chi_1\psi_1\overline{\chi_2\psi_2}(h)+O(W^\varepsilon R^{1/2})
\end{align}
where $W=[q_1,q_2]d$, 
and thus for $D, R\leq x^\delta$ and uniformly for $|a_2|\leq x^\delta$ we have 
\begin{align}
\mathcal S_3(w, a_3)=\hat{\alpha}(0)X_3+O(x^{\varepsilon} R^{5/2} N^2/D^2)
\end{align}
where 
\begin{align}
X_3=\sum_{\substack{ (q_i, a_1a_2d)=1\\ q_i \equiv c_0\pmod c}}\frac{\gamma(q_1)\gamma(q_2)}{[q_1,q_2]d\phi((q_1,q_2))\phi(d)}\sum_{\substack{\chi\pmod{(q_1,q_2)d}\\ \operatorname{cond}(\chi)\leq R}}\sum_{\substack{n_i\asymp N\\ (n_i, q_ia_2d)=1}}\beta_{n_1}\overline{\beta_{n_2}}\chi(n_1\overline{n_2}).
\end{align}
Note that $X_3$ is independent of $w,a_3$, and the error term is acceptable for \eqref{dispersion} (in fact much better than required).

Similarly we can evaluate $\mathcal S_2$ as in \cite[Section 5.3.2]{ Dr2}
so that for $D, R\leq x^\delta$ and uniformly for $|a_2|\leq x^\delta$ we have
\begin{align}
\mathcal S_2(w,a_3)=\hat{\alpha}(0)X_3+O(MN^2D^{-2}x^{-2\delta}).
\end{align}

We next evaluate $\mathcal S_1(w, a_3)$.
After Poisson summation in $m$, we get the expected main term as $\hat{\alpha}(0)X_1$, where
\begin{align}
X_1=\sum_{\substack{(q_i, a_1a_2d)=1\\ q_i\equiv c_0 \pmod c}}\frac{\gamma(q_1)\gamma(q_2)}{[q_1,q_2]d}\sum_{\substack{n_i \asymp N\\ (n_i, q_id)=1\\ n_1\equiv n_2\pmod{d(q_1,q_2)}}}\beta_{n_1}\overline{\beta_{n_2}},
\end{align}
which is also independent of $w, a_3$.
The rest of the evaluation of $\mathcal S_1$ is similar to that in  \cite[Proposition 4.2]{ABL} (with $a_2$ replaced by $a_2a_3$). The condition $(q_i,a_2d)=1$ is equivalent to $(q_i, a_2a_3)=1$ and $(q_i, \tilde{d})=1$ where $\tilde{d}=d/(d, (a_2a_3c)^\infty)$. Thus the additional condition $(q_i,\tilde{d})=1$ can be incorporated by following the proof of \cite[Proposition 4.2]{ABL} using M\"obius inversion. To be precise, we replace the condition $\delta_i\mid a_1$ by the condition that $\delta_i\mid a_1\tilde{d}$ so that the condition $(\delta_2, n_0a_2a_3c)=1$ is still satisfied when applying \cite[Theorem 2.3]{ABL} since $(n_i,d)=1$. In conclusion, we have for  $\delta=\delta(\eta),\kappa=\kappa(\eta)$ small enough, $D, R\leq x^\delta$ and uniformly for $|a_2|\leq x^\delta$ that
\begin{align}
\mathcal S_1(w,a_3)=\hat{\alpha}(0)X_1+O(MN^2 D^{-2} x^{-2\delta}).
\end{align}

It remains show that $$\mathcal X:=\sum_{d\asymp D}\sup_{\substack{w\pmod d^\times\\a_3\mid d}}|X_1-X_3| \ll N^2D^{-1}R^{-2/3}(\log x)^{O(1)}.$$ (Note that $X_1, X_3$ are independent of $w, a_3$.) The proof follows along the lines of \cite[Section 5.6]{Dr2}. We have
\begin{align}
\mathcal X&=\sum_{d\asymp D}\Big|\sum_{\substack{(q_i, a_1a_2d)=1\\ q_i \equiv c_0\pmod c}}\frac{\gamma(q_1)\gamma(q_2)}{[q_1,q_2]d\phi((q_1,q_2)d)\phi(d)}\sum_{\substack{ \chi \pmod {(q_1,q_2)d}\\ \operatorname{cond}(\chi)> R}}\sum_{\substack{n_i\asymp N\\ (n_i, q_ia_2d)=1}}\beta_{n_1}\overline{\beta_{n_2}}\chi(n_1\overline{n_2})\Big|\\
	&\ll \log_2 x\sum_{d\asymp D}\sum_{\substack{q_i\asymp Q\\ (q_i, a_1a_2d)=1}}\frac{1}{q_1q_2d^2}\sum_{\substack{ \chi \text{ primitive}\\ \operatorname{cond}(\chi)> R\\\operatorname{cond}(\chi)\mid (q_1,q_2)d}}\Big|\sum_{\substack{n_i\asymp N\\ (n_i, q_ia_2d)=1}}\beta_{n_1}\overline{\beta_{n_2}}\chi(n_1\overline{n_2})\Big|.
\end{align}	
Before we can estimate this with the large-sieve, we need to uncouple the variables. To do so, we write $\text{cond}(\chi)=ql$ and detect the conditions $(n_i,q_id)=1$ by M\"obius inversion. We thus estimate the above
\begin{align*}
	&\ll \log_2 x \sum_{\substack{\ell\ll D \\ q\ll Q \\ lq>R\\ (l,q)=1}} \sum_{\substack{r_1,r_2,s\\ r_i\ll Q/q \\ s\ll D/\ell}} \sum_{\substack{q_i \asymp Q, d\asymp D\\ r_iq \mid q_i \\ sl\mid d}} \frac{1}{q_1q_2 d^2} \sum_{\substack{ \chi \text{ primitive }\\\chi \pmod {q\ell}\\ }}\Big|\sum_{\substack{n_i\\ (n_i, a_2)=1}}\beta_{r_1sn_1}\overline{\beta_{r_2sn_2}}\chi(n_1\overline{n_2})\Big|\\
&\ll \frac{(\log x)^3}{D} \sum_{\substack{r_1,r_2,s\\ r_i\ll Q \\ s\ll D}} \frac{1}{r_1 r_2 s} \sum_{\substack{\ell\ll D \\ q\ll Q \\ lq>R}} \frac{1}{q^2l}  \sum_{\substack{ \chi \text{ primitive }\\\chi \pmod {q\ell}\\ }}\Big|\sum_{\substack{n_i\\ (n_i, a_2)=1}}\beta_{r_1sn_1}\overline{\beta_{r_2sn_2}}\chi(n_1\overline{n_2})\Big|.
\end{align*}
By Cauchy-Schwarz and the symmetry between $n_1$ and $n_2$, we arrive at
\begin{align}
\mathcal X\ll  \frac{(\log x)^{O(1)}}{D}\sum_{r\ll x}\frac{\tau(r)}{r}\sum_{\substack{\ell\ll D \\ q\ll Q \\ lq>R\\ (l,q)=1}}\frac{1}{q^2 \ell}\sum_{\substack{\chi \text{ primitive}\\\chi\pmod {q\ell}}}\Big|\sum_{\substack{n\asymp N\\ (n,a_2)=1} }\beta_{rn}\chi (n) \Big|^2.
\end{align}
By the large sieve inequality for Dirichlet characters (\cite[Lemma 3.3]{Dr2}) and partial summation, the contribution from $q>R^{1/3}$ can be bounded by
\begin{align}
& \ll (\log x)^{O(1)}\frac{1}{D}\Big(\frac{(D^2Q^2+N)N}{Q^2}+\int_{R^{1/3}}^{2Q} \frac{(D^2t^2+N)N}{t^3}dt\Big) \ll (\log x)^{O(1)}\Big(ND+\frac{N^2}{DR^{2/3}}\Big).
\end{align}
Similarly the contribution from $\ell\geq R^{2/3}, q\leq R^{1/3}$ can be bounded by 
\begin{align}
&\ll (\log x)^{O(1)}\frac{1}{D}\Big( \frac{(R^{2/3}t^2+N)N}{t}\Big\vert_{R^{2/3}}^D+\int_{R^{2/3}}^{2D} \frac{(R^{2/3}t^2+N)N}{t^2}dt\Big)\\
&\ll (\log x)^{O(1)}\Big(R^{2/3}N+\frac{N^2}{DR^{2/3}}\Big).
\end{align}
Therefore, we have $$\sum_{d\asymp D}\sup_{\substack{w\pmod d^\times\\a_3\mid d}}| X_1-X_3|\ll (\log x)^{O(1)} \frac{N^2}{DR^{2/3}}$$ if $D\ll N/R^{1/3}$. Since $N\gg x^\eta$, we conclude that for $\delta$ small enough (in terms of $\eta$) equation \eqref{dispersion} holds and so \eqref{bilinearsum} follows. This completes the proof of Proposition \ref{TypeII}.

 \section{Proofs of Theorems \ref{thm1} and \ref{thm1'}}\label{prime-section}
 

\subsection{Primes in arithmetic progressions}  Before we start with the analysis of the main term, we give a result on averaged equidistribution of primes in long arithmetic progressions, which will be used to handle the error terms in the sieving process.
\begin{prop}\label{PrimeAp}
	There exists some absolute constant $\varpi>0$ such that the following holds. 
	Let $x\geq 2$, $c_0,c,d, C\in \mathbb N, (c_0, c)=1, a_1,a_2\in \mathbb{Z}\setminus \{0\}$ such that for 
	\begin{align}
	Q\leq x^{1/2+\varpi},\quad |a_1|\leq x^{1+\varpi}, \quad |a_2|\leq x^\varpi, \quad |\lambda_d|\ll \tau(d)^C
	\end{align} 
	we have 
	\begin{align}
	\sum_{d\leq x^{\varpi}}\lambda_d\sup_{\substack{w\pmod d^\times}}\Big|\sum_{\substack{q\leq Q\\ (q,a_1a_2d)=1\\ q\equiv c_0 \pmod c}}\Big(\sum_{\substack{n\leq x\\n\equiv a_1\overline{a_2}\pmod q\\ n\equiv w\pmod d}}\Lambda(n)-\frac{1}{\phi(qd)}\sum_{\substack{n\leq x\\ (n, qd)=1}}\Lambda(n)\Big)\Big|\ll_{C,A} cx(\log x)^{-A}.
	\end{align}
\end{prop}
\textbf{Proof.}
This is the analogue of \cite[Theorem 2.1]{ABL} and can be proved in the exact same way as at the end of \cite[Section 4]{ABL}. We only need to replace  \cite[Proposition 4.1]{ABL} with Proposition \ref{TypeII} with the choice $R=(\log x)^B$ for some large $B$ depending on $C, A$. The contribution from $\tau_3$-type sums is negligible by choosing $\varpi$ small enough since \cite[Lemma 2]{BFI2} saves a small power of $x$. Using $\sum_{d\leq x^\kappa} \lambda_d/\phi(d)\ll (\log x)^{O(1)}$ and the Siegel-Walfisz Theorem we also have that the contribution from $\chi\pmod{qd}$ with $\operatorname{cond}(\chi)\leq R $ can be bounded by $x R (\log x)^{O(1)}\exp(-b\sqrt{\log x})$, which is negligible by the choice of $R$.
	
\medskip	
	
For future reference we remark that  we can replace the von Mangoldt function by the characteristic function on primes.
 \subsection{Preparing for the sieve}\label{52}
 In this section, we use $\kappa$ to denote a sufficiently small (depending on   $\varpi$  in Proposition \ref{PrimeAp} and $A, \eta$ as in Theorem \ref{prime}) positive constant.
 We use $\varepsilon$ to denote an arbitrarily small positive constant. 
 
 Let $\textbf{d} = (d_1, d_2)$ denote a pair of two square-free numbers with $d_1, d_2 \leq n^{\kappa}$. In preparation for a sieve we define
 $$\mathcal{A}_{\textbf{d}}(n) = \{(p, x_1, x_2) \mid p + d_1^2x_1^2 + d_2^2x_2^2 = n\}$$
(where of course $p$ denotes a prime).  As usual, we denote by $c_p(n)$  the Ramanujan sum.  The key input for the sieve is the following:

\begin{prop}\label{sieve}With the above notation, there exists some absolute constant $\kappa>0$ such that uniformly for $d_1, d_2\leq n^{\kappa}$ and $\lambda_{\bd}\ll (\tau(d_1d_2))^C$ we have 
\begin{align}\label{averageAd}
\operatornamewithlimits{\sum\sum}_{d_1,d_2\leq n^\kappa}\lambda_{\bd}\Big|\#\mathcal{A}_{\textbf{d}}(n)  -  \operatorname{Li}(n)\mathfrak S_{\textbf{d}} (n)\Big|\ll_{C,A} n(\log n)^{-A},
\end{align}
where the singular series $\mathfrak S_{\textbf{d}} (n)$ is given by 
\begin{equation}\label{euler}
\frac{\pi}{d_1d_2}\prod_{\substack{p\mid (d_1,d_2)}}\Big(1-\frac{c_p(n)}{p-1}\Big)\prod_{\substack{p\mid d_1d_2\\ p\nmid(d_1, d_2)\\ p\nmid 2n}}\Big(1-\frac{(\frac{n}{p})}{p-1}\Big)\prod_{p\nmid d_1d_2n}\Big(1+\frac{\chi_{-4}(p)}{p(p-1)}\Big)\prod_{\substack{p\mid n\\ p\nmid 2d_1d_2}}\Big(1-\frac{\chi_{-4}(p)}{p}\Big).
\end{equation}
\end{prop}

 \emph{Remark:} One can check that this matches  the main term coming from a formal application of the circle method. In particular, we see that $\mathfrak S_{\textbf{d}} (n) = 0$ if $(d_1, d_2, n) > 1$ as it should. \\

\textbf{Proof.}  Let $\delta=(d_1, d_2), d_1'=d_1/\delta, d_2'=d_2/\delta$ as in Lemma \ref{lem1}. By Lemma \ref{lem1} and Theorem \ref{prime} (with $\kappa$ sufficiently small depending on the constants $A, \eta$) we have 
\begin{equation}\label{asymp}
\#\mathcal{A}_{\textbf{d}}(n) = \sum_{\delta_1\delta_1'=d_1'}\sum_{\delta_2\delta_2'=d_2'}\frac{w_{\delta_1'\delta_2'}}{h_{\delta_1'\delta_2'}}G_{\delta_1'\delta_2'}\sum_{w\in\mathcal G_{\delta_1'\delta_2'}}\sum_{\substack{p\equiv n \pmod {(\delta\delta_1'\delta_2')^2}\\n-p\equiv w\pmod {\delta_1'\delta_2'}}}r\left(\frac{n-p}{\delta^2 \delta_1^2\delta_2^2}\right) + O(n^{1-\kappa}),
\end{equation}
where $\mathcal{G}_{\delta_1'\delta_2'} \subseteq (\Bbb{Z}/\delta_1'\delta_2'\Bbb{Z})^{\ast}$ is as in \eqref{residues} 
and  $r = \textbf{1} \ast \chi_{-4}$ as before. From \eqref{classnumber} we obtain
 \begin{align}
\frac{w_{\delta_1'\delta_2'}}{h_{\delta_1'\delta_2'}}
=\frac{4}{\delta_1'\delta_2'}\prod_{p\mid \delta_1'\delta_2'}\Big(1-\frac{\chi_{-4}(p)}{p}\Big)^{-1}
\end{align}
in all cases, including $\delta_1'\delta_2' = 1$. Note that we have 
\begin{align}\label{whgbound}
\frac{w_{\delta_1'\delta_2'}}{h_{\delta_1'\delta_2'}}G_{\delta_1'\delta_2'}\# \mathcal G_{\delta_1'\delta_2'}\ll (\tau(\delta_1'\delta_2'))^{O(1)}.
\end{align}We apply the Dirichlet's hyperbola method to the $r$-function so that the innermost $p$-sum in \eqref{asymp} can for some parameter $Y > 0$ be written as 
\begin{align}\sum_{Y\leq p\leq n}\sum_{\substack{ab=(n-p)/\delta^2\delta_1^2\delta_2^2\\ n-p\equiv w\pmod{\delta_1'\delta_2'}, }}\chi_{-4}(b)+O(Y^{1+\varepsilon}) =S_1+S_2+O(Y^{1+\varepsilon}),
\end{align}
say, where 
\begin{align}
S_1&=\sum_{\substack{a\leq \sqrt{n}\\ (a, \delta_1'\delta_2')=1}}\sum_{\substack{Y\leq p\leq n\\ p\equiv n\pmod {a\delta^2\delta_1^2\delta_2^2}\\ n-p\equiv w \pmod{\delta_1'\delta_2'}}}\chi_{-4}\Big(\frac{n-p}{a\delta^2\delta_1^2\delta_2^2} \Big),\\
S_2&=\sum_{\substack{ b\leq \frac{n-Y}{\sqrt{n}\delta^2\delta_1^2\delta_2^2}\\ (b, \delta_1'\delta_2')=1}}\sum_{\substack{ Y\leq p\leq n-b\sqrt{n}\delta^2\delta_1^2\delta_2^2\\ p\equiv n \pmod {b\delta^2\delta_1^2\delta_2^2}\\ n-p\equiv w\pmod{ \delta_1'\delta_2'}}}\chi_{-4}(b).
\end{align}
We choose $Y=n^{1-3\kappa}$ so that the error term is acceptable after summing over $d_1,d_2\leq n^\kappa$. We will evaluate $S_1$ and $S_2$ on average over $d_1,d_2$ by Proposition \ref{PrimeAp}.  This requires some preparation.  It turns out that the contribution from $S_1$ will not appear in the main term $\operatorname{Li}(n)\mathfrak S_\bd(n)$, but as there are some complications from the powers of two, we start with a detailed treatment for $S_1$. 
 
We use the notation $$r_2=(\delta\delta_1\delta_2,2^\infty), \quad d=\frac{\delta\delta_1\delta_2}{r_2}, \quad d'=\frac{\delta_1'\delta_2'}{(\delta_1'\delta_2',2)}, \quad a_2=(a, 2^\infty), \quad a_d=(a, d^\infty).$$ We see that $S_1$ equals
\begin{align}
&\sum_{\substack{a_2\mid 2^\infty\\ (a_2, \delta_1'\delta_2')=1}}\sum_{\substack{a\leq \sqrt{n}/a_2\\ (a, 2\delta_1'\delta_2')=1}}\sum_{\substack{Y\leq p\leq n\\ p\equiv n\pmod{aa_2r_2^2d^2}\\n-p\equiv w\pmod{\delta_1'\delta_2'}}}\chi_{-4}(a)\chi_{-4 }\Big(\frac{n-p}{a_2r_2^2}\Big)\\
=&\sum_{\substack{a_2\mid 2^\infty\\ (a_2, \delta_1'\delta_2')=1}}\sum_{a_d\mid d^\infty}\chi_{-4}(a_d)\underset{u \, (\text{mod } 4)}{\left.\sum\right.^{\ast}} \chi_{-4}(u )\sum_{v\pmod 4}\chi_{-4}(v)\sum_{\substack{a\leq \sqrt{n}/a_2a_d\\ (a, 2d_1d_2)=1\\ a\equiv u\pmod 4}}\sum_{\substack{Y\leq p\leq n \\ p\equiv n\pmod{aa_2a_dr_2^2d^2}\\n-p\equiv w\pmod {d'}\\ n-p \equiv a_2r_2^2v \, (\text{mod }4a_2r_2^2)
}}1.
\end{align}
Here we have used the fact that if $2\mid \delta_1'\delta_2'$ then we must have $a_2r_2=1$ and $(w,2)=1$.  
The contribution when one of $a_2,a_d$ is at least $n^{4\kappa}$ can be bounded by
\begin{align}
\sum_{\substack{a\mid (2d)^\infty\\ a\geq n^{4\kappa}}}\sum_{q\leq \sqrt{n}/a}\frac{n}{\phi(a)\phi(q)}\ll n^{1-3\kappa+\varepsilon},
\end{align}
which is negligible after summing over $d_1,d_2\leq n^\kappa$. So from now on we restrict to $ a_2,a_d\leq n^{4\kappa}$ (note that $r_2,  d', d$ are automatically at most $n^{2\kappa}$). We see that the $a, p$-sum equals
$$ \sum_{\substack{a\leq \sqrt{n}/a_2a_d\\(a, ndd')=1\\a\equiv u \pmod 4}}\sum_{\substack{Y\leq p\leq n \\ p\equiv n \pmod {a} \\p\equiv n- w\pmod {d'}\\p\equiv n-va_2r_2^2\pmod{4a_2r_2^2}\\  p\equiv n \pmod{a_dd^2}}}1+O(n^{1/2+\varepsilon})$$
where the error term comes from the artificially added condition $(a, n) = 1$. For the same reason we can restrict to $(n, d) =(n-w,d')=1$. After combining the last three congruence conditions modulo $4a_2r_2^2d'a_dd^2$ to a single condition and noting $d_1d_2\mid 4a_2r_2^2d'a_dd^2$ and $\lambda_\bd\ll \tau(d_1d_2)^C$, we are in a position to apply Proposition \ref{PrimeAp} (since $\kappa$ is sufficiently small in terms of $\varpi$) together with the prime number theorem to the two innermost sums and recast the previous display as
\begin{align}
&\mathds{1}_{(n-w, d')=(n, d)=1} \sum_{\substack{a\leq \sqrt{n}/a_2a_d\\ a\equiv u \pmod 4\\ (a, 2d_1d_2n)=1}}\frac{\operatorname{Li}(n)}{\phi(a)\phi(d')\phi(4a_2a_dr_2^2d^2)}+r_{\textbf{d}}
\end{align}
where the error terms $r_\bd$ satisfy
\begin{align}
\sum_{d_1,d_2\leq n^\kappa}\lambda_{\bd} \sum_{\delta_1\delta_1'=d_1'}\sum_{\delta_1\delta_2'=d_2'}\frac{w_{\delta_1'\delta_2'}}{h_{\delta_1'\delta_2'}}G_{\delta_1'\delta_2'}\sum_{w\in \mathcal G_{\delta_1'\delta_2'}}|r_\bd|\ll x(\log x)^{-A}.
\end{align}
Plugging back, we obtain
\begin{align}
S_1 =  &\mathds{1}_{(n-w,d')=(n, d)=1}\sum_{\substack{a_2\mid 2^\infty \\ a_2 \leq n^{4\kappa}}}\sum_{\substack{a_d\mid d^\infty\\ a_d \leq n^{4\kappa}}}\sum_{\substack{ v\pmod {4}\\ (n-va_2r_2^2, 4)=1}}\chi_{-4}(v)\underset{u\pmod {4}}{\left.\sum\right.^{\ast}}\chi_{-4}(u)\\
 &\quad\quad\quad\quad\quad\times
\sum_{\substack{a\leq \sqrt{n}/a_2a_d\\ a\equiv u \pmod 4\\ (a, 2d_1d_2n)=1}}\frac{\operatorname{Li}(n)}{\phi(a)\phi(d')\phi(4a_2a_dr_2^2d^2)} + r_{\bd}. 
\end{align}
But now the $v$-sum vanishes: with $h = a_2r_2^2 \mid 2^{\infty}$ we have 
\begin{displaymath}
  \sum_{\substack{ v\pmod {4}\\ (n-vh, 4)=1}}\chi_{-4}(v) = \sum_{\ell \mid 4} \mu(\ell) \sum_{\substack{v\, (\text{mod } 4) \\ vh \equiv n \, (\text{mod } \ell)}}\chi_{-4}(v).
\end{displaymath}
Since $\chi_{-4}$ is primitive, by \cite[(3.9)]{IK} the inner sum vanishes, unless $4 \mid \ell/(\ell, h)$ (which can only happen if $h=1$), but then the M\"obius function vanishes. 
We conclude the contribution from $S_1$ in \eqref{averageAd} is acceptable. \\

We now turn to $S_2$ 
 where the calculation is similar, and we obtain 
\begin{align}
S_2
&=\sum_{b_d\mid {(\delta\delta_1\delta_2)^\infty}}\chi_{-4}(b_d)\sum_{\substack{ b\leq\frac{n-Y}{ \sqrt{n}(\delta\delta_1\delta_2)^2b_d}\\ (b, d_1d_2)=1}}\chi_{-4}(b)\sum_{\substack{ Y\leq p\leq n-bb_d\sqrt{n}\delta^2\delta_1^2\delta_2^2\\ p\equiv n \pmod {b_d\delta^2 \delta_1^2\delta_2^2}\\ p\equiv n \pmod b\\ n-p\equiv w \pmod{\delta_1'\delta_2'}}}1\\
&=\underset{u \, (\text{mod } 4)}{\left.\sum\right.^{\ast}}\chi_{-4}(u) \sum_{b_d\mid {(\delta\delta_1\delta_2)^\infty}}\chi_{-4}(b_d)\sum_{\substack{ b\leq\frac{n-Y}{ \sqrt{n}(\delta\delta_1\delta_2)^2b_d}\\ b\equiv u\pmod 4\\(b, d_1d_2)=1}}\sum_{\substack{ Y\leq p\leq n-bb_d\sqrt{n}\delta^2\delta_1^2\delta_2^2\\ p\equiv n \pmod {b_d\delta^2 \delta_1^2\delta_2^2}\\ p\equiv n \pmod b\\p\equiv n-w\pmod {\delta_1'\delta_2'}}}1. 
\end{align}
As before we restrict to  $b_d\leq n^{4\kappa}$ at the cost of an error of $O(x^{1-\kappa})$ (after summing over $d_1,d_2\leq n^\kappa$). The contribution from $(\delta\delta_1\delta_2,n)>1$ is also negligible, so let us from now on assume $(\delta\delta_1\delta_2,n)=1$. We cannot apply Proposition \ref{PrimeAp} to evaluate the $b, p$-sum on average over $d_1,d_2$ directly due to the dependence of the length of the intervals. 
To remedy this, we split the $p$-sum into intervals of the form $((1 - \Delta)P, P]$ with $\Delta = (\log n)^{-K}$ for some sufficiently large $K=K(C)$. 
 The condition $p\leq n-bb_d\sqrt{n}\delta^2\delta_1^2\delta_2^2$ interferes in at most one of such intervals, whose contribution we estimate by the Brun-Titchmarsh inequality at the cost of an admissible error  $O(n \Delta (\log x)^{O(1)})$ which is negligible by choosing $K$ large enough in terms of $C$. For all other intervals we can split the $b$-sum in to intervals of the shape $((1-\Delta)P,P]$. The condition $b\leq (n-Y)/\sqrt{n}(\delta\delta_1\delta_2)^2b_d$ interferes in one of such intervals, whose contribution we can be bound by $O(nR(\log x)^{O(1)}\Delta)$ using the trivial bound for $\mathfrak u_R(n, q,a,d,w)\ll \frac{R\tau(q)\tau(d)}{\phi(q)\phi(d)}$. Here $R$ is chosen to be some large power of $\log x$ in the proof of Proposition \ref{PrimeAp} and so the error term from separating the $b$-variable is negligible when $K$ is sufficiently large. For the rest of the intervals, we can apply Proposition \ref{PrimeAp}. Assembling the main term, we recast $S_2$ as
 \begin{align}
&\underset{u \, (\text{mod } 4)}{\left.\sum\right.^{\ast}}\chi_{-4}(u) \sum_{ \substack{b_d\mid {(\delta\delta_1\delta_2)^\infty}\\ b_d \leq n^{4\kappa}}}\chi_{-4}(b_d)\sum_{\substack{ b\leq\frac{n-Y}{ \sqrt{n}(\delta\delta_1\delta_2)^2b_d}\\ b\equiv u\pmod 4\\(b, d_1d_2n)=1}}
\sum_{\substack{ Y\leq p\leq n-bb_d\sqrt{n}\delta^2\delta_1^2\delta_2^2\\ (p,b_d\delta^2 \delta_1^2\delta_2^2)=1\\ (p, b)=1\\(p,\delta_1'\delta_2')=1}}\frac{\mathds{1}_{(n-w,\delta_1'\delta_2')=(n, \delta\delta_1\delta_2)=1}}{\phi(b)\phi(b_d\delta^2\delta_1^2\delta_2^2)\phi(\delta_1'\delta_2')}+r'_{\bd} \\
&=\operatorname{Li}(n)\sum_{\substack{ b_d\mid {(\delta\delta_1\delta_2)^\infty}\\ b_d \leq n^{4\kappa}}}\frac{\chi_{-4}(b_d)}{\phi(b_d\delta^2\delta_1^2\delta_2^2)}\sum_{\substack{ b\leq\frac{n-Y}{ \sqrt{n}(\delta\delta_1\delta_2)^2b_d}\\ (b, d_1d_2n)=1}}\frac{\chi_{-4}(b)}{\phi(b)}
\frac{\mathds{1}_{(n-w, \delta_1'\delta_2')=(n, \delta\delta_1\delta_2)=1}}{\phi(\delta_1'\delta_2')}+r'_\bd,
\end{align}
where the error terms $r'_\bd$ satisfy 
\begin{align}
\sum_{d_1,d_2\leq n^\kappa}\lambda_{\bd} \sum_{\delta_1\delta_1'=d_1'}\sum_{\delta_1\delta_2'=d_2'}\frac{w_{\delta_1'\delta_2'}}{h_{\delta_1'\delta_2'}}G_{\delta_1'\delta_2'}\sum_{w\in \mathcal G_{\delta_1'\delta_2'}}|r'_\bd|\ll x(\log x)^{-A}.
\end{align}
As a last step, we compute the $b$-sum in the main term in $S_2$ using \cite[Lemma 5.2]{ABL} and complete the sum over $b_d$, and recast the previous display up to an admissible error as
$$\operatorname{Li}(n)L(1, \chi_{-4})\prod_{p}\Big(1+\frac{\chi_{-4}(p)}{p(p-1)}\Big) \mathfrak c(d_1d_2n)\sum_{ b_d\mid {(\delta\delta_1\delta_2)^\infty}}\frac{\chi_{-4}(b_d)}{\phi(b_d\delta^2\delta_1^2\delta_2^2)}\frac{\mathds{1}_{(n-w,d')=(n, \delta\delta_1\delta_2)=1}}{\phi(\delta_1'\delta_2')},$$
where
\begin{equation}\label{cf}
\mathfrak c(f) = \prod_{p \mid f} \Big(1 + \frac{\chi_{-4}(p)}{p(p-1)}\Big) ^{-1}\Big(1 - \frac{\chi_{-4}(p)}{p}\Big).
\end{equation}
This finishes the analysis of $S_2$.

We now return to \eqref{asymp} and summarize that there exists some absolute constant $\kappa>0$ such that 
\begin{align}
\operatornamewithlimits{\sum\sum}_{d_1,d_2\leq n^\kappa}\lambda_{\bd}\Big|\#\mathcal{A}_{\textbf{d}}(n)  -  \operatorname{Li}(n)\mathfrak S_{\textbf{d}} (n)\Big|\ll_{C,A} n(\log n)^{-A}\end{align}
where the singular series $\mathfrak S_{\textbf{d}}(n)$ is given by
\begin{displaymath}
\begin{split}
\mathfrak  S_{\bd} (n)=  &\underset{(\delta\delta_1\delta_2, n) = 1}{ \sum_{\delta_1\delta_1'=d_1'}\sum_{\delta_2\delta_2'=d_2'}}\frac{4}{\delta_1'\delta_2'}\prod_{p\mid \delta_1'\delta_2'}\Big(1-\frac{\chi_{-4}(p)}{p}\Big)^{-1}\frac{G_{\delta_1'\delta_2'}}{\phi(\delta_1'\delta_2')}\sum_{\substack{w\pmod{\delta_1'\delta_2'}\\(n-w,\delta_1'\delta_2')=1\\ w\in \mathcal G_{\delta_1'\delta_2'}}}1
\\&\times 
\sum_{ b_d\mid {(\delta\delta_1\delta_2)^\infty}}\frac{\chi_{-4}(b_d)}{\phi(b_d\delta^2\delta_1^2\delta_2^2)}L(1, \chi_{-4}) \prod_{p}\Big(1+\frac{\chi_{-4}(p)}{p(p-1)}\Big)\mathfrak c(d_1d_2n).\end{split}
\end{displaymath}
The sum over $w\pmod {\delta_1'\delta_2'}$ can be computed using \eqref{residues} as 
$$\prod_{\substack{p\mid \delta_1'\delta_2'\\p\nmid 2n}}\Big(\frac{p-1}{2}-\frac{\chi_{p^*}(n)+1}{2}\Big)\prod_{\substack{p\mid \delta_1'\delta_2'\\ p\mid n}}\Big(\frac{p-1}{2}\Big)\mathds{1}_{(n-1, \delta_1'\delta_2',2)=1}.$$
It is now a straightforward exercise with Euler products, noting that $\chi_{p^{\ast}}(n) = (n/p)$ for odd $p$ and
 $L(1, \chi_{-4}) = \pi/4$, to obtain the expression \eqref{euler}. This completes the proof.  

\subsection{Completion of the proof of Theorem \ref{thm1}}\label{53}

With Proposition \ref{sieve} available, Theorem \ref{thm1} follows now easily from application of a sieve. Let $\mathcal A=(a_k)$ be a finite sequence of non-negative real numbers.  
We introduce the notation
\begin{align}\label{sievedset}
 S(\mathcal A, z)=\sum_{\substack{(k, P_1(z))=1}}a_k, \quad P_q(z)=\prod_{\substack{p\leq z,p\nmid q}}p, \quad 
 A_c=\sum_{\substack{k\equiv 0\pmod c}}a_k,\quad A=A_1.
\end{align}
With these notations, we recall the following lemma from sieve theory (see e.g.\ \cite[Theorem 1]{IwS}). 
\begin{lemma}\label{sievelemma}
	 Let $\gamma, L>0$ be some fixed constants. Let $\Omega(c)$ be a multiplicative function satisfying $0\leq \Omega(p)<p$ and
	\begin{align}\label{sievedimension}
	\prod_{w\leq p\leq w'}\Big(1-\frac{\Omega(p)}{p}\Big)^{-1}\leq \Big(\frac{\log w'}{\log w}\Big)^{\gamma}\Big(1+\frac{L}{\log w}\Big)
	\end{align}	
 for all $2\leq w\leq w'$. 	Then we have 
	\begin{align}
S(\mathcal A, z)\geq A \prod_{p\leq z}\Big(1-\frac{\Omega(p)}{p}\Big)\Big(f_\gamma(s)-\frac{e^{\sqrt{L}}Q(s)}{(\log \mathcal D)^{1/3}}\Big)-\sum_{\substack{c\leq \mathcal D\\ c\mid P_1(z)}}\Big| A_c- \frac{\Omega(c)}{c}A\Big|,
	\end{align}
	where $s= \log \mathcal D/\log z$,
 $Q(s)< \exp(-s\log s+s\log\log3s+O(s))$, and $f_\gamma(s)$ is some continuous function such that $0< f_\gamma(s)<1$ and $f(s)=1+O(e^{-s})$ as $s\rightarrow \infty$. 
\end{lemma}
We are going to sieve the sequence
 $\mathcal{A}=\mathcal A(n)=(a_k)$ where $a_k=\#\{k=x_1x_2: p+x_1^2+x_2^2=n\}$. 
Then from Proposition \ref{sieve}, there exists some $\kappa>0$ such that for $c\leq n^{\kappa}$ with $\mu^2(c)=1$ we have
\begin{align}\label{bc}
A_c=\mu(c)\sum_{\substack{\textbf{d}\\ p\mid d_1d_2\Leftrightarrow p\mid c}}\mu(d_1)\mu(d_2)\frac{\omega(\textbf{d};n)}{d_1d_2}\mathfrak S(n)\operatorname{Li}(n) + r_c,
\end{align}
where 
\begin{align}
\mathfrak S(n)&=\pi\prod_{p\nmid n}\Big(1+\frac{\chi_{-4}(p)}{p(p-1)}\Big)\prod_{\substack{p\mid n}}\Big(1-\frac{\chi_{-4}(p)}{p}\Big),\\
\omega(\textbf{d};n)&=\prod_{\substack{p\mid (d_1,d_2)}}\Big(1-\frac{c_p(n)}{p-1}\Big)\prod_{\substack{p\mid d_1d_2\\ p\nmid(d_1, d_2)\\ p\nmid 2n}}\Big(1-\frac{(\frac{n}{p})}{p-1}\Big)\prod_{\substack{p\nmid n\\ p\mid d_1d_2}}\Big(1+\frac{\chi_{-4}(p)}{p(p-1)}\Big)^{-1}\prod_{\substack{p\mid n\\ p\mid 2d_1d_2}}\Big(1-\frac{\chi_{-4}(p)}{p}\Big)^{-1}\\
&:=\prod_{p^\nu\parallel d_1d_2}\omega_\nu(p;n),
\end{align}
and 
\begin{align}\label{rcbound}
\operatornamewithlimits{\sum}_{c\leq n^\kappa}\Big|A_c- \mu(c)\sum_{\substack{\textbf{d}\\ p\mid d_1d_2\Leftrightarrow p\mid c}}\mu(d_1)\mu(d_2)\frac{\omega(\textbf{d};n)}{d_1d_2}\mathfrak S(n)\operatorname{Li}(n)\Big|\ll_{K} n(\log n)^{-K}\end{align}
for any $K>0$. We next examine the multiplicative function $\omega_\nu(\bullet;n)$.
For odd $p$ we compute explicitly 
\begin{align}
\omega_1(p;n)&=
\begin{cases}
\displaystyle\Big (1-\frac{(n/p)}{p-1}\Big)\Big(1+\frac{\chi_{-4}(p)}{p(p-1)}\Big)^{-1},& p\nmid n,\\
\displaystyle \Big(1-\frac{\chi_{-4}(p)}{p}\Big)^{-1}, & p\mid n,
\end{cases}
\end{align}
and
\begin{align}
\omega_2(p;n)&=
\begin{cases}
\displaystyle\Big(1-\frac{c_p(n)}{p-1}\Big)\Big(1+\frac{\chi_{-4}(p)}{p(p-1)}\Big)^{-1},& p\nmid n,\\
\displaystyle\Big(1-\frac{\chi_{-4}(p)}{p}\Big)^{-1}, & p\mid n,
\end{cases}
\end{align}
and for $p=2$ we have 
\begin{align}
\omega_1(2;n)=1, \quad  \omega_2(2;n)= 1-\frac{c_2(n)}{p-1}.
\end{align}
Let $\Omega(\bullet;n)$ be the multiplicative function defined by 
$$ \Omega(p;n)=2\omega_1(p;n)-\omega_2(p;n)p^{-1}$$
so that 
we have from \eqref{bc} and \eqref{rcbound} that for any $\mathcal D\leq n^{\kappa}$
\begin{align}
\sum_{c\leq \mathcal D} \mu^2(c) \Big|  A_c-\frac{\Omega(c;n)}{c} A\Big| \ll   n(\log n)^{-K}
\end{align}
where $A=\mathfrak S(n)\operatorname {Li}(n)$. Note that we have $$0\leq \Omega(p;n)< 3, \quad \Omega(p;n)=2+O(1/p).$$
Thus there exists some absolute constant $L$ such that \eqref{sievedimension} holds for $\gamma=2$. We choose $\mathcal D=n^{\kappa_0}$ for some absolute small constant $\kappa_0$ (e.g.\ $\kappa_0=\kappa/3$) and choose $C=C(\kappa_0)$ sufficiently large such that $f_2(C\kappa_0)>0$ holds. Since $\Omega(2;n)=2$ if and only if $n\equiv 0 \pmod 2$ and $\Omega(3;n)=3$ if and only if $n\equiv 2\pmod 3$, it follows that $0\leq \Omega(p)<p$ if $n\equiv 1,3 \pmod 6$.  Then, on taking \(s=C\kappa_0,\) it follows from Lemma \ref{sievelemma} that for $n\equiv 1,3\pmod 6$ we have 
\begin{align}
S(\mathcal A, n^{1/C})&\geq \mathfrak S(n) \operatorname{Li} (n) \prod_{p\leq n^{1/C} }\Big(1-\frac{\Omega(p;n)}{p}\Big) \Big( f_2(C\kappa_0)- \frac{ e^{\sqrt{L}}Q(C\kappa_0)}{(\log \mathcal D)^{1/3}} \Big)-O( n (\log n)^{-A})\\
&\gg \mathfrak S(n)\operatorname{Li}(n)(\log n)^{-2}.
\end{align}
 This completes the proof of Theorem \ref{thm1}. 

\subsection{Proof of Theorem \ref{thm1'}} 
Denote \begin{align}
& \mathcal B=\mathcal B(x)=(b_k), \quad b_k=\#\{k=x_1x_2: p=x_1^2+x_2^2+1, p\leq x\}
\end{align}
Then we have 
\begin{align}
\#\{(p,x_1,x_2): p=x_1^2+x_2^2+1, P^-(x_1x_2)\geq p^{1/C},p\leq x\}\geq S(\mathcal B, x^{1/C}),
\end{align}
and thus it is enough to prove
\begin{align}\label{SB}
S(\mathcal B, x^{1/C})\gg x(\log x)^{-3}
\end{align} for some sufficiently large $C$. We can follow the proof of Theorem \ref{thm1} to obtain \eqref{SB}. To be more precise, 
let $B_c=B_c(x)=\sum_{k\equiv 0 \pmod c}b_k$. Then for $\mu^2(c)=1$ we have 
\begin{align}
B_c=\mu(c)\sum_{\substack{\bd=(d_1,d_2)\\ p\mid d_1d_2\Leftarrow p\mid c}}\mu(d_1)\mu(d_2)\#\mathcal B_\bd(x),
\end{align}
where \begin{align}
&\mathcal B_\bd(x)=\{ (p, x_1,x_2):p=d_1^2x_1^2+d_2^2x_2^2+1, p\leq x\}.
\end{align}
After an application of Lemma \ref{lem21} and Theorem \ref{prime}, we see that  uniformly for $d_1,d_2\leq x^\kappa$ with $\kappa$ sufficiently small we have
\begin{align}
\#\mathcal B_\bd(x)=\sum_{\delta_1\delta_1'=d_1'}\sum_{\delta_2\delta_2'=d_2'}\frac{w_{\delta_1'\delta_2'}}{h_{\delta_1'\delta_2'}}G_{\delta_1'\delta_2'}\sum_{w\in\mathcal G_{\delta_1'\delta_2'}}\sum_{\substack{x^{1-\kappa}\leq p\leq x\\ p\equiv 1 \pmod {(\delta\delta_1'\delta_2')^2}\\p-1\equiv w\pmod {\delta_1'\delta_2'}}}r\left(\frac{p-1}{\delta^2 \delta_1^2\delta_2^2}\right) + O(x^{1-\kappa+\epsilon}).
\end{align}
Then we can follow the proof of Proposition \ref{sieve} to evaluate the main term and obtain that there exists some absolute constant $\kappa>0$ such that uniformly for $d_1, d_2\leq n^\kappa$ and $\lambda_\bd\ll \tau(d_1d_2)^{C'}$ we have
\begin{align}\label{Bd1d2}
\operatornamewithlimits{\sum\sum}_{d_1, d_2\leq x^\kappa}\lambda_\bd \big| \#\mathcal B_\bd(x)-\operatorname{Li}(x)\mathfrak S_\bd(1)\big|\ll_{C', A} x(\log x)^{-A},
\end{align}
where $\mathfrak S_\bd(n)$ is as in \eqref{euler}. With \eqref{Bd1d2}, we can follow the proof and notation in Subsection \ref{53} with $n$ replaced by $1$ in the singular series calculations to see that 
\begin{align}\label{Bc1}
\sum_{c\leq n^\kappa}\Big|B_c-\frac{\Omega(c;1)}{c}\mathfrak S(1)\operatorname{Li}(x)\Big|\ll_A x(\log x)^{-A}
\end{align} for any $A>0$. Note that have $\Omega(p;1)<p$ for all $p$ and $\mathfrak S(1)>0$. Lemma \ref{sievelemma} together with \eqref{Bc1} then yields \eqref{SB} for $C$ sufficiently large. For the second claim in Theorem \ref{thm1'}, we
write $$R_C(n)=\#\{(x_1,x_2):n=x_1^2+x_2^2+1, P^-(x_1x_2)\geq n^{1/C}\}.$$ Since
$
R_C(n)\leq r(n-1)\ll n^\epsilon,
$
we see that 
\begin{align}
\sum_{p\leq x}\mathds{1}_{R_C(p)>0}&\gg x^{-\epsilon}\#\{(p,x_1,x_2): p=x_1^2+x_2^2+1, P^-(x_1x_2)\geq p^{1/C}, p\leq x\}\\&\gg x^{-\epsilon}S(\mathcal B, x^{1/C})\gg x^{1-\epsilon'},
\end{align}
which shows that there are infinitely many primes that can be written as one plus the sum of two squares of integers having no more than $C/2$ prime factors.

 \section{Proof of Theorem \ref{thm2}}\label{smooth-section}
 
 The rest of the paper is devoted to the proof of Theorem \ref{thm2}. The basic idea is similar to the proof of Theorem \ref{thm1}, but it is technically and structurally substantially more involved. Again we start with the distribution of smooth numbers in arithmetic progressions to large moduli.

\subsection{Smooth numbers in arithmetic progressions}
Denote 
$$\Psi(x, y)=\sum_{\substack{m\in S(x,y)}}1,\quad  \Psi_q(x, y)=\sum_{\substack{ m\in S(x,y)\\ (m, q)=1}}1.$$
We also use the usual notation $$u=\frac{\log x}{\log y}, \quad H(u)=\exp \Big(\frac{u}{\log^2 (u+1)}\Big), \quad u\geq 1.$$
For an integer $n$ we set $n_y=\prod_{\substack{p^\nu\|n\\ p\leq y}}p^{\nu}$ and $n_o=\frac{n}{(n,2^\infty)}$.  
We have the following analogue of Proposition \ref{PrimeAp} for smooth numbers. 
\begin{prop}\label{SmoothAp}
	Suppose $|\lambda_d|\ll \tau(d)^C$ for some fixed $C>0$. Then there exist some constants $\varpi, \delta, D >0$ such that the following holds. 
	Let $x\geq 2$, $c_0,c,d\in \mathbb N, (c_0, c)=1, a_1,\ell\in \mathbb{Z}\setminus \{0\}$ such that for 
	\begin{align}
	Q\leq x^{1/2+\varpi},\quad |a_1|\leq x^{1+\varpi}, \quad |a_2|\leq x^\varpi, \quad (\log x)^D\leq y\leq x, \quad a_2\mid \ell, \quad \omega(\ell_y)\ll \log x
	\end{align} 
	we have 
	\begin{align}
	&\sum_{d\leq x^{\varpi}}\lambda_d\sup_{\substack{w\pmod d^\times\\a_3\mid d}}\Big|\sum_{\substack{q\leq Q\\ (q,a_1a_2d)=1\\ q\equiv c_0 \pmod c}}\Big(\sum_{\substack{m\in S(x,y)\\\mu^2(m_o)=(m, \ell)=1\\ m\equiv a_1\overline{a_2a_3} \pmod q\\ m\equiv w\pmod d}}1-\frac{1}{\phi(qd)}\sum_{\substack{m\in S(x,y)\\\mu^2(m_o)=(m,qd\ell)=1}}1\Big)\Big|\ll_{C, A}c \Psi_{\ell}(x,y) (\log x)^{-A}
	\end{align}
	for any $A>0$.
\end{prop}
\textbf{Remark.} One of the important features of Proposition \ref{SmoothAp} is that we make the dependence on $\ell$ explicit. This is necessary for our application since $\Psi_\ell(x,y)\asymp \prod_{p\mid \ell_y}(1-p^{-\alpha})\Psi(x,y)$ with $\alpha = \alpha(x, y)$ as in \eqref{alphadef}, which could be much smaller than $\Psi(x,y)(\log x)^{-A}$ when $y$ is small. Furthermore, the condition $a_2\mid \ell$ ensures that $(n,a_2)=1$ in an application of Proposition \ref{TypeII}. \\

As a preparation for the proof we recall the following result that is a variation of the result in \cite[Section 3.3]{Ha} or \cite[Lemme 5]{Dr1}. 
\begin{lemma}\label{SmoothCharSum}
	Suppose $|\lambda_q|\ll\tau(q)^C$ for some fixed $C>0$ and $\ell \in \mathbb{Z}\setminus\{0\}$. There exist some constants $D, \eta, \delta>0$ such that the following is true. If $(\log x)^D\leq y\leq x$, $q\leq Q\leq x$ and $\omega(\ell_y)\ll \log x$, we have 
	\begin{align}
	\sum_{q\leq Q}\frac{\lambda_q}{\phi(q)}\sum_{\substack{\chi \pmod q\\ 1<\operatorname{cond}(\chi)\leq x^\eta}}\Big|\sum_{\substack{m\in S(x,y)\\ (m, \ell )=1}}\chi(m) \Big|\ll_{C, A} \Psi_\ell(x,y)(H(u)^{-\delta}(\log x)^{-A}+y^{-\delta}).
	\end{align}
	for any $A>0$.
\end{lemma}
\textbf{Proof.} When $\lambda_q=1$, $\ell=1$, this can be found in \cite[Section 3.3]{Ha}. For general coefficients $\lambda_q$ we can follow the same argument using that 
$$\sum_{q\leq Q}\frac{|\lambda_q|}{\phi(q)} \ll (\log x)^{O(1)}\quad \text{and} \quad  \lambda_q\ll q^\varepsilon\, \text{ for any } \varepsilon>0.$$  

We now deal with general values for $\ell$. 
We recall \cite[Theorem 2.4, Lemme 3.1, eq (2.8)]{lBT}: we have uniformly in $\mathcal L(x) \leq y \leq x$, $1\leq d \leq x/y$ $P^+(\ell)\leq y$ and $\omega(\ell)\leq \sqrt{y}$ the estimate
\begin{equation}\label{psix/d}
\Psi_\ell(x/d, y) \ll \prod_{\substack{p\mid \ell}}(1-p^{-\alpha}) \frac{\Psi(x, y)}{d^{\alpha}}
\end{equation}
where $\alpha$ is as in \eqref{alphadef} satisfying
\begin{align}\label{alphaasym}
\alpha= \alpha(x, y) = 1-\frac{\xi(u) + O(1)}{\log y}, \quad u=\frac{\log x}{\log y}
\end{align} 
and $\xi(t)$ is defined implicitly by $$\quad e^{\xi(t)}  = 1 + t \xi(t).$$ It is easy to see that $\xi(u) = \log (u \log u) + O(1)$.
When $\exp((\log x)^{2/5})< y\leq x$, we have 
\begin{align}
\prod_{p\mid \ell_y}\Big(1-\frac{1}{p^{\alpha}}\Big)\gg\exp\Big(-\sum_{p\leq y}p^{-\alpha}\Big)\gg \exp\Big(-\frac{u\log u}{\log y}\Big)\gg 1, 
\end{align}
which implies that $\Psi_\ell(x,y)\asymp \Psi(x,y)$. Thus, assuming without loss of generality that \(p\mid \ell \Rightarrow p\leq y\), and after using M\"obius inversion to detect the condition $(m, \ell)=1$, it is enough to prove that
\begin{align}\label{coprimel}
\sum_{q\leq Q}\frac{\lambda_q}{\phi(q)}\sum_{\substack{ \chi \pmod q\\ 1<\operatorname{cond}(\chi)\leq x^{\eta}}}\sum_{d\mid \ell}\Big|  \sum_{m\leq x/d}\chi(m)\Big|\ll \Psi(x,y)\big(H(u)^{-\delta}(\log x)^{-A}+y^\delta\big).
\end{align} 
The contribution from $d\geq x^{c_0}$ can be bounded trivially by $x^{1-c_0+\eta+\epsilon}$. The contribution from $d\leq x^{c_0}$ can be bounded individually using the result for $\ell=1$ (and adjusting the constant $\eta$) so that together with \eqref{psix/d} we obtain the left hand side of \eqref{coprimel}
 can be bounded by
\begin{align}
&\sum_{\substack{ d\mid \ell}}\sum_{q\leq Q}\frac{\lambda_q}{\phi(q)} \sum_{\substack{ \chi \pmod q\\ 1< \operatorname{cond}(\chi)\leq x^{\eta}}} \Big|\sum_{m\leq x/d}\chi(m)\Big|\\
& \ll \sum_{\substack{ d\mid \ell \\d\leq x^{c_0}}}\Psi\Big( \frac{x}{d},y\Big)\Big(H\Big(u-\frac{\log d}{\log y}\Big)^{-\delta}(\log x)^{-A}+y^{-\delta}\Big)+O(x^{1-c_0+\eta+\epsilon})\\
&\ll\sum_{d\mid \ell}\frac{1}{d^{\alpha}}d^{\frac{2\delta}{\log y(\log u)^2}}\Psi(x,y)\big(H(u)^{-\delta}(\log x)^{-A}+y^{-\delta}\big)+O(x^{1-c_0+\eta+\epsilon})\\
&\ll \Psi(x,y)\big(H(u)^{-\delta}(\log x)^{-A}+y^{-\delta}\big),
\end{align} 
where we used the the assumption of $\ell$ and the range of $y$ in the last step.

It remains 
to incorporate the condition $(m, \ell)=1$ for $(\log x)^D \leq y\leq \exp((\log x)^{2/5})$. To this end, we first generalize \cite[Proposition 1, Theorem 3]{Ha} with $\Psi(x,y)$ replaced by $\Psi_\ell(x,y)$. This can be done by following the proof of \cite[Proposition 1, Theorem 3]{Ha} and replacing $L(s, \chi;y)=\prod_{p\leq y}(1-\chi(p)p^{-s})^{-1}$  by $$L_\ell(s,\chi; y):=\prod_{\substack{p\nmid \ell, p\leq y}}(1-\chi(p)p^{-s})^{-1}.$$ Using $\omega(\ell_y)\ll \log x$ we see that the contribution from $\sum_{(n, \ell)>1}\Lambda(n)\chi(n)n^{-\sigma}$ can be bounded by $$\sum_{n\ll \log x\log_2x}\Lambda(n)n^{-\sigma}\ll \frac{(\log x\log_2x)^{1-\sigma}}{1-\sigma},$$ which is $o(\log x)$ when $\sigma\in[\alpha-1/300, \alpha]$ using \eqref{alphaasym}. Thus we still have the bound 
$$|\log L_\ell(\sigma+it, \chi;y)-\log L_\ell(\alpha+it, \chi;y)| \leq \frac{(\alpha-\sigma)\log x}{2},$$ which is enough to prove the analogue of \cite[Theorem 3]{Ha} with $\Psi(x,y)$ replaced by $\Psi_\ell(x,y)$ when the summand is restricted to $(m, \ell)=1$. We can also replace \cite[Smooth Numbers Result 1, 3]{Ha} by \cite[Theorem 2.1, 2.4]{lBT}. 
 Then we can follow the proof in \cite[Section 3.3]{Ha} in the case $\ell=1$ and replace $\Psi(x,y)$ by $\Psi_\ell(x,y)$ so that condition $(m,\ell)=1$ is preserved. On the way we need the estimate 
 \begin{align}\label{lowerubound}
\sum_{\substack{p\leq y\\ p\nmid q\ell}}\frac{1-\Re(\chi(p)p^{-it})}{p^{\alpha}}\gg\frac{u}{\log^2(u+1)}
 \end{align}
 which was proved on \cite[pp. 16-17]{Ha} when $\ell=1$. For general $\ell$, we simply note that the contribution from $p\mid \ell$ can be trivially bounded by 
 \begin{align}
\sum_{p\mid \ell} \frac{1}{p^\alpha}\ll \frac{(\sqrt{y}\log \sqrt{y})^{1-\alpha}}{(1-\alpha)\log y}\ll \sqrt{u}
 \end{align} 
 so that \eqref{lowerubound} still holds since $u\gg (\log_2x)^2$ by our current assumption on $y$.
\\

\textbf{Proof of Proposition \ref{SmoothAp}:}
The proof follows by combining the proofs of \cite[Lemme 2.1]{FT90} and \cite[Theorem 2.1]{ABL}. We use $\varepsilon$ as an arbitrary positive constant, not necessarily the same at each occurrence. Let $R=x^\eta$ for some sufficiently small $\eta$ to be determined later. 
We use the notation $\mathfrak u_R(n;q, a,d, w)$ as defined in \eqref{urqd} to write 
\begin{align}
\sum_{d\leq x^\varpi}\lambda_d\sup_{\substack{w\pmod d^\times\\a_3\mid d}}\Big|\sum_{\substack{q\leq Q\\ (q, a_1a_2d)=1\\ q\equiv c_0 \pmod c}}\Big(\sum_{\substack{m\in S(x,y)\\\mu^2(m_o)=(m, \ell)=1\\ m\equiv a_1\overline{a_2a_3} \pmod q\\ m\equiv w\pmod d}}1-\frac{1}{\phi(qd)}\sum_{\substack{m\in S(x,y)\\ \mu^2(m_o)=(m, qd\ell)=1}}1\Big)\Big|\leq \mathcal S_1+\mathcal S_2, \end{align}
where 
\begin{align}
&\mathcal S_1=\sum_{d\leq x^\varpi}|\lambda_d|\sup_{\substack{w\pmod d^\times\\a_3\mid d}}\Big|\sum_{\substack{q\leq Q\\ (q, a_1a_2d)=1\\ q\equiv c_0\pmod c}}\sum_{\substack{m\in S(x,y)\\ \mu^2(m_o)=(m, \ell)=1}}\frac{1}{\phi(q)\phi(d)}\sum_{\substack{\chi \pmod q\\\psi \pmod d\\ 1<\operatorname{cond}(\chi\psi)\leq R}}\chi(m\overline{a_1}a_2a_3)\psi(m\overline{w})\Big|,\\
&\mathcal S_2=\sum_{d\leq x^\varpi}|\lambda_d|\sup_{\substack{w\pmod d^\times\\a_3\mid d}}\Big|\sum_{\substack{q\leq Q\\ (q, a_1a_2d)=1\\ q\equiv c_0\pmod c}}\frac{1}{\phi(d)}\sum_{\substack{m\in S(x,y)\\ \mu^2(m_o)=(m,\ell)=1}}u_R(m;q,a_1\overline{a_2a_3},d,q)\Big|.
\end{align}
For $\mathcal S_1$, we use Lemma \ref{SmoothCharSum} and $(q,d)=1$ to conclude that for $\eta$ small enough, 
 \begin{align}
\mathcal S_1\ll & \sum_{d\leq x^\varpi}|\lambda_d|\sup_{\substack{w\pmod d^\times\\a_3\mid d}}\Big|\sum_{\substack{ q\leq Q\\ (q, a_2a_2d)=1\\ q\equiv c_0 \pmod c}}\frac{1}{\phi(qd)}\sum_{\substack{ \chi \pmod q\\ \psi \pmod d\\1<\operatorname{cond}(\chi\psi) \leq R}}\sum_{\substack{m\in S(x, y)\\ \mu^2(m_o)=(m,\ell)=1}}\chi(m\overline{a_1}a_2a_3)\psi(m\overline{w})\Big|\\
 &\ll \sum_{k\leq x^\eta}\sum_{q\leq Qx^\varpi}\frac{\tilde{\lambda}_q}{\phi(q)}\sum_{\substack{\chi \pmod q\\ 1< \operatorname{cond}(\chi)\leq R}}\Big|\sum_{\substack{m\in S(x/k^2,y)\\ (m, \ell)=1}}\chi(m)\Big|+O(x^{1-\eta+\epsilon})\\
 &\ll \sum_{k\leq x^\eta}\Psi_\ell(x/k^2,y)(H(u-\frac{2\log k}{\log y})^{-\delta}(\log x)^{-A}+y^{-\delta})\ll \Psi_{\ell}(x,y) (\log x)^{-A'}
 \end{align}
 for some large constant $A'$ depending on $y, \delta$. Here we used \eqref{psix/d}, $y\geq (\log x)^D$ and $\alpha\geq 1-1/D+o(1)$ in the last step.
Using ideas as in the   proof of Theorem \ref{smooth}, we now complete the argument by showing that 
 there exists some  $\varpi, \delta_0>0$ such that
 \begin{align}\label{largeq}
\mathcal S_2\ll cx (\log x)^{O(1)}R^{-\delta_0}.
 \end{align} 
 To prove \eqref{largeq}, we consider separately the cases  $y\geq x^\eta$ and $y< x^\eta$ for some sufficiently small $\eta>0$.  

 When $y\geq x^\eta$, we use Buchstab's identity \eqref{SmoothBuchstab} to write 
\begin{align}
\mathcal S_2\leq  S_0+ \sum_{i=1}^3 S_i+\overline{S},
\end{align}
where 
\begin{align}
S_0&=\sum_{d\leq x^\varpi}|\lambda_d|\sup_{\substack{w\pmod d^\times\\a_3\mid d}}\Big|\sum_{\substack{q\leq Q\\ (q, a_1a_2d)=1\\q\equiv c_0 \pmod c}}\Big(\sum_{\substack{ m\equiv \overline{a_1}a_2a_3\\\mu^2(m_o)=(m, \ell)=1\\ m\leq x}}1-\frac{1}{\phi(qd)}\sum_{\substack{\chi \pmod q\\\psi \pmod d\\ \operatorname{cond}(\chi\psi)\leq R}}\sum_{\substack{\mu^2(m_o)=1\\(m,\ell)=1\\m\leq x}}\chi(m\overline{a_1}a_2a_3)\psi(m\overline{w})\Big)\Big|,\\
S_j&=\sum_{d\leq x^\varpi}|\lambda_d|\sup_{\substack{w\pmod d^\times\\a_3\mid d}}\Big|\sum_{\substack{ q\leq Q\\ (q, a_1a_2d)=1\\ q\equiv c_0 \pmod c}}\sum_{\substack{y<p_1< \dots< p_j\\\mu^2(m_op_1\cdots p_j)=1\\( mp_1\cdots p_j,\ell)=1\\ mp_1\cdots p_j\leq x}}u_R(mp_1\cdots p_j;q,a_1\overline{a_2a_3},d,w)\Big|,\\
\overline{S}&=\sum_{d\leq x^\varpi}|\lambda_d|\sup_{\substack{w\pmod d^\times\\a_3\mid d}}\Big|\sum_{\substack{ q\leq Q\\ (q, a_1a_2d)=1\\ q\equiv c_0 \pmod c}}\sum_{\substack{y<p_1< \dots< p_4\\  \mu^2(m_op_1\cdots p_4)=1\\(mp_1\cdots p_4,\ell)=1\\mp_1\cdots p_4\leq x}}u_R(mp_1\cdots p_4;q,a_1\overline{a_2a_3},d,w)g_{p_4}(m)\Big|.
\end{align}
Consider $\overline{S}$ first. Note that we can assume $y< p_1\ll x^{1/4}$ since otherwise the sum is empty. The condition $\mu^2(m_op_1\cdots p_4)=1$ is equivalent to $\mu^2(m_o)=(m, p_1\cdots p_4)=1$. We use M\"obius inversion to detect $(m, p_1\cdots p_4)=1$. Since we have that the primes $p_i>y\geq x^\eta$, the contribution from $k\mid p_1\cdots p_4$ with $k<x^{\eta/2}$ comes only from the term $k=1$. Similar to the proof of Theorem \ref{smooth}, we localize $m$ and $p_i$ into intervals $m\in (M, M(1+\mathcal{Z})]$ and $p_i\in (P_j, P_{j+1}]$ with $P_j=y(1+\mathcal Z)^j$ and  with $\mathcal Z=R^{-\eta_0}$ for some $\eta_0>0$ so that 
\begin{align}
\overline{S}\ll \mathcal Z^{-4}(\log x)^{O(1)}\sup_{ M, P_{j_i}}\tilde{\mathcal S} + x(\log x)^{O(1)}\mathcal Z+O(x^{1-\eta/2+\epsilon})
\end{align}
where 
\begin{align}
\tilde{\mathcal S}=\sum_{d\leq x^\varpi}|\lambda_d|\sup_{\substack{w\pmod d^\times\\a_3\mid d}}\Big|\sum_{\substack{q\leq Q\\ (q, a_1a_2d)=1\\ q\equiv c_0\pmod c}}\sum_{\substack{p_i \in (P_{j_i}, (1+\mathcal Z)P_{j_i}]\\ m\in (M, M(1+\mathcal Z)]\\ (mp_1\cdots p_4, \ell)=1\\ 1\leq i\leq 4}}g_{p_4}(m)\mu^2(m_o)u_R(mp_1\cdots p_4;q,a_1\overline{a_2},d,w)\Big|. 
\end{align}
If $M \prod_{i=1}^4P_{j_i}\leq x^{1-\eta/2}$, we can bound $\tilde{\mathcal S}\ll x^{1-\eta/2+\varepsilon}$ trivially. Otherwise, we can apply the bilinear estimate Proposition \ref{TypeII} with $\bm{m}=p_1$ to $\tilde{\mathcal S}$ since $x^{\eta}\leq y\leq P_{j_1}\ll x^{1/4}\ll (M\prod_{i=1}^4P_{j_i})^{1/4+\eta}$. We  obtain 
 \begin{align}\label{sbar}
 \overline{S}\ll \mathcal Z^{-4}cx (\log x)^{O(1)}R^{-1/3}+x (\log x)^{O(1)}\mathcal Z+O(x^{1-\eta/2+\epsilon})
 \end{align}
 as long as 
 \begin{align}\label{QRconditions}
Q\leq x^{1/2+\varpi_0}, \quad c, R, |a_2|\leq x^{\varpi_0},\quad a_1\leq x^{1+\varpi_0}
 \end{align} for some constant $\varpi_0>0$. This suffices for the bound in \eqref{largeq} by choosing a suitable $\eta_0$.
 
 For $S_k$, $1\leq k\leq 3$, we can use M\"obius to detect $(m, p_1\cdots p_j)=1$ as before, replace the characteristic function of the primes by the von Mangoldt function, apply Heath-Brown's identity \eqref{HeathBrownId} to each prime variable $p_i$ and then localize the new variables so that we have 
 \begin{align}
 S_k\ll \mathcal Z^{-24}(\log x)^{O(1)}\sup_{\substack{ M,  M_i,  N_j\\ 1\leq r\leq 12\\ 1\leq s\leq 12}}\dot{\mathcal S}+ x(\log x)^{(1)}\mathcal Z+O(x^{1-\eta/2+\epsilon}),
 \end{align}
 where 
  \begin{align}
 \dot{\mathcal S}=&\sum_{d\leq x^\varpi}\sup_{\substack{w\pmod d^\times\\a_3\mid d}}\Big|\sum_{\substack{ q\leq Q\\ (q, a_1a_2d)=1\\ q\equiv c_0 \pmod c}}\operatornamewithlimits{\sum\sum\sum}_{\substack{m\in (M,(1+\mathcal Z) M]\\m_i\in (M_i, (1+\mathcal Z)M_i]\\n_j\in(N_j, (1+\mathcal Z)N_j]\\(mm_in_j,\ell)=1\forall i,j}} \mu^2(m_o)\prod_{i=1}^r\mu(m_i) u_R\Big(m\prod_{i=1}^rm_i\prod_{j=1}^sn_j;q,a_1\overline{a_2a_3},d,w\Big)\Big|
 \end{align}
 with  $$ M\prod_{i=1}^rM_i\prod_{j=1}^sN_j\ll x, \quad M_i\leq x^{1/4}. $$
 As before, we can assume $M\prod_{i=1}^rM_i \prod_{j=1}^s N_j\geq x^{1-\eta/2}$ since otherwise we can bound $\dot {\mathcal S}\ll x^{1-\eta/2+\varepsilon}$ trivially. Now we argue similarly as in the proof of Theorem \ref{smooth}. If there is a subset $\mathcal I \subseteq \{M, M_i, N_j\}$ such that $K:=\prod_{I\in \mathcal I}I$ such that $x^{\eta}\leq K \leq x^{1/4+\eta}$, then we can apply Proposition \ref{TypeII} 
  to obtain 
 $$ \dot {\mathcal S}\ll cx(\log x)^{O(1)}R^{-1/3}$$
 as long as \eqref{QRconditions} is satisfied. Otherwise, we can use M\"obius inversion to detect $\mu^2(m_o)=1$ if necessary, and it is enough to show there exists some constant $\varpi>0$ such that
 \begin{align}\label{smoothd3}
 \sup_{\substack{w\pmod d^\times\\a_3\mid d}}\Big|\sum_{\substack{q\leq Q\\ (q, a_1a_2d)=1\\q\equiv c_0 \pmod c}} \operatornamewithlimits{\sum\sum\sum}_{\substack{m\asymp M, k\asymp K, l\asymp L\\ (mkl,\ell)=1}}\alpha_ru_R(rmkl;q, a_1\overline{a_2a_3},d,w) \Big|\ll x^{1-\varpi}
 \end{align}
for $r\leq x^{\eta}, x^{1-\eta/2}\ll rMKL\leq x$ and $\alpha_r\ll r^\varepsilon$. After an application of M\"obius inversion to detect the condition $(mkl,\ell)=1$,  we are left with
\begin{align}\label{smoothd32}
	\sup_{\substack{w\pmod d^\times\\a_3\mid d}}\Big|\sum_{\substack{\lambda\mid \ell}}\mu(d)\sum_{\substack{q\leq Q\\ (q, a_1a_2d)=1\\q\equiv c_0 \pmod c}} \operatornamewithlimits{\sum\sum\sum}_{\substack{m\asymp M, k\asymp K, l\asymp L\\}}\alpha_ru_R(r\lambda mkl;q, a_1\overline{a_2a_3},d,w)\Big| \ll x^{1-\varpi}
\end{align}
for $r\leq x^{\eta}, x^{1-\eta/2}\ll r\lambda MKL\leq x$ and $\alpha_r\ll r^\varepsilon$. The contribution from $\lambda>x^{\eta'}$ can be bounded trivially by $x^{1-\eta'/2}$. The contribution for any fixed $\lambda< x^{\eta'}$ can be rewritten and then  estimated in the same way as in \cite[eq.\ (4.5)]{ABL}, where ultimately Deligne's estimates for exponential sums over algebraic varieties over finite fields are used. By choosing $\mathcal Z, \eta'$ suitably, we see that there exists some $\varpi, \delta_0>0$ such that $S_k\ll cx(\log x)^{O(1)}R^{-\delta_0}$ for $1\leq k\leq 3$.

We can also treat $S_0$ by \eqref{smoothd3}, which completes the proof of \eqref{largeq} when $y\geq x^\eta$. 

\medskip
 
 When $y\leq x^\eta$, we can apply \cite[Lemme 3.1]{FT96} in the form of \eqref{smoothfactor} so that it is enough to show that there exists some absolute constant $\varpi, \delta_0>0$ such that 
 \begin{equation}\label{smally}
 \sum_{d\leq x^\varpi}\lambda_d\sup_{\substack{w\pmod d^\times\\a_3\mid d}}\Big|\sum_{\substack{q\leq Q\\ (q, a_1a_2d)\\ q\equiv c_0\pmod c}}\sum_{\substack{ lm\leq x,(lm, \ell)=1\\ \mu^2((lm)_o)=1\\ P^+(l)\leq P^-(m)\\ M\leq m\leq yM}}g_y(l)g_y(m)u_R(lm;q,a_1\overline{a_2a_3},d,w)\Big|\ll cx(\log x)^{O(1)}R^{-\delta_0}. 
 \end{equation}
 Since $P^+(\ell)\leq P^-(m)$ we see that  \begin{align}\label{smally2}
 &\sum_{d\leq x^\varpi}\lambda_d\sup_{\substack{w\pmod d^\times\\a_3\mid d}}\Big|\sum_{\substack{q\leq Q\\ (q, a_1a_2d)\\ q\equiv c_0\pmod c}}\sum_{\substack{ lm\leq x,(lm, \ell)=1\\ \mu^2((lm)_o)=1\\ P^+(l)\leq P^-(m)\\ M\leq m\leq yM}}g_y(l)g_y(m)u_R(lm;q,a_1\overline{a_2a_3},d,w)\Big|\\
 &\ll 
 \sum_{d\leq x^\varpi}|\lambda_d|\sup_{\substack{w\pmod d^\times\\a_3\mid d}}\Big|\sum_{\substack{q\leq Q\\ (q, a_1a_2d)\\ q\equiv c_0\pmod c}}\sum_{\substack{ lm\leq x,(lm, \ell)=1\\ P^+(l)< P^-(m)\\ M\leq m\leq yM}}\mu^2(l_o)g_y(l)\mu^2(m_o)g_y(m)u_R(lm;q,a_1\overline{a_2a_3},d,w)\Big|\\
 & \quad+\sum_{d\leq x^\varpi}|\lambda_d|\sup_{\substack{w\pmod d^\times\\a_3\mid d}}\Big|\sum_{\substack{q\leq Q\\ (q, a_1a_2d)\\ q\equiv c_0\pmod c}}\sum_{p\leq y}\sum_{\substack{ lm\leq x/p^2,(lm, \ell)=1\\ P^+(l)< p<P^-(m)\\ M\leq mp\leq yM}}\mu^2(l_o)g_y(l)\mu^2(m_o)g_y(m)u_R(lm;q,a_1\overline{a_2a_3p^2},d,w)\Big|.
 \end{align}
 After separating variables $l$ and $m$ as in \cite[eq. (3.37)]{Dr1} and localizing $l, m$ into short intervals $l\asymp L$ and $m\asymp M$, we can apply Proposition \ref{TypeII} for the first sum above to obtain \eqref{smally}, as long as $x^{\eta_1}\leq M$ and $My\leq (LM)^{1/4+\eta_1}$, which can be satisfied by choosing $M=x^{\eta_2}$ for some suitable $\eta_2>0$ when $LM\geq x^{1-\eta_3}$. If $LM\leq x^{1-\eta_3}$, we can obtain \eqref{smally} trivially. For the second sum, we move the summation of $p\leq y$ outside so that contribution from $p\leq x^{\eta_4}$ can be estimated using Proposition \ref{TypeII} similarly as above and the contribution from $y\geq x^{\eta_4}$ can be bounded trivially. This complete the proof of \eqref{largeq} when $y\leq x^\eta$.

\subsection{Implementing sieve weights}
The most straightforward approach to Theorem \ref{thm2} would be to  use Proposition \ref{SmoothAp} to evaluate $$\sum_{m\in S(n,y)}r_{d_1,d_2}(n-m)$$ 
for square-free numbers $d_1, d_2$ and $y = g(n)$, and then apply a sieve. Unfortunately  the main term in the asymptotic formula is not a multiplicative function of $d_1d_2$. As mentioned in the introduction, this makes it problematic to implement a sieve, so instead we will work with the sieve weights directly. We recall \cite[Lemma 3]{IwS2}.
\begin{lemma}\label{lowerboundsieve}
	Let $\gamma, L>0$ be some fixed constants. Then
	there exists a sequence $\lambda_c^-$ supported on square-free integers less than $\mathcal D$ such that 
	$$\lambda_1^-=1, \quad  |\lambda_c^-|\leq 1,  \quad \mu*1 \geq \lambda^-* 1$$
	and that for all multiplicative functions $w$ satisfying $0\leq w(p)<1$ and 
	\begin{align}
	\prod_{W\leq p\leq W'}(1-w(p))^{-1}\leq \Big(\frac{\log W'}{\log W}\Big)^{\gamma}\Big(1+\frac{L}{\log w}\Big)
	\end{align}	
	 for all $2\leq W\leq W'$ we have 
	\begin{align}
	&\sum_{\substack{c\mid P(z)}}\lambda_c^-w(c)\geq\prod_{p\leq z}(1-w(p))\Big(f_\gamma(s)+O\Big(\frac{e^{\sqrt{L}-s}}{(\log \mathcal D)^{1/3}}\Big)\Big),  \quad  s=\frac{\log \mathcal D}{\log z}\geq 2, 
	\end{align}
	where $f_\gamma(s)$ is some continuous function such that $0< f_\gamma(s)<1$ and $f_\gamma(s)=1+O(e^{-s})$ as $s\rightarrow\infty$. 
\end{lemma}
 The idea is to estimate the sieve weights $\lambda_c^-$ with multiplicative factors $w(c)$ first to get a lower bound using Lemma \ref{lowerboundsieve}, and then sum over the smooth numbers $m$.  
 In order to avoid major technical difficulties, we need to prepare the set-up very carefully.
\begin{itemize}
\item 
Many estimates are sensitive to the prime factorization of $m$, and it simplifies our life to restrict $m$ to be square-free. On the other hand, if $n \equiv 2$ (mod 4), we will need $4\mid m$ in order to represent $n$ by $m$ and two odd integer squares. Thus we only restrict $m$ to have an odd square-free part. This restriction will not change the order of the number of solutions, but the fact that square-free numbers are not equidistributed among {\em all} residue classes will leads to different constants for different $n$. 
\item It turns out that a factor $\chi_{-4}(h)$ with $h\mid (n,m)$ appears in the main term. It is convenient to force  $\chi_{-4}(h)>0$ (so that a lower bound suffices), and thus we also restrict $m$ to be  coprime to all primes $p \equiv 3$ (mod 4) dividing $n$. Since $n-m$ is a sum of two squares, this is a harmless maneuver and affects only square factors,
 but it makes an important difference for the computation if we  implement this. Note that we could simplify the computations by restricting $m$ to be coprime to $n$, but this could lose some order of magnitude for certain cases of $n$, e.g.\ when $n$ is the product of the first few primes that are $1 \pmod 4$.
		\item Finally, the sieve weights behave a bit erratically for small primes, so we sieve them out directly by M\"obius inversion.
		\end{itemize}
With these general remarks in mind, we fix the following notation. We write
$$\mathfrak N=\prod_{\substack{p\mid n\\ p\equiv 3\pmod 4}}p, \quad m_{o}=\prod_{\substack{p^\nu\parallel m\\ p>2}}p^{\nu} = \frac{m}{(m, 2^{\infty})}.$$ Recall the notation \eqref{sievedset}.  Let $$\mathcal{A}= \mathcal{A}(n)=(a_k),\quad a_k=\#\Big\{k=x_1x_2: \begin{array}{ll}
m+x_1^2+x_2^2=n,\\ P^+(m)\leq y,
\end{array}
\begin{array}{ll}
\mu^2(m_{o})=1,\\(m,\N)=1
\end{array}\Big\}.$$
Then for $c$ square-free, we have
\begin{align}
\mathcal{A}_c=\mu(c)\sum_{\substack{\bd\\ p\mid d_1d_2\Leftrightarrow p\mid c}}\mu(d_1)\mu(d_2)\sum_{\substack{m\leq n\\ P^+(m)\leq y\\ \mu^2(m_{o})=(m,\N)=1}}r_{d_1,d_2}(n-m). 
\end{align}
 Let $\mathcal Q$ be some fixed absolute constant (we will later choose $\mathcal{Q} = 30$). Let $\lambda_c^-$ be some lower bound sieve weights as in Lemma \ref{lowerboundsieve} of dimension $\gamma$, level $\mathcal D=n^{\kappa_0}$ for some sufficiently small $\kappa_0$ and $z=n^{1/C}$ for some sufficiently large $C$ so that $f_\gamma(C\kappa_0)>0$.
Then we have the following lower bound for $ S(\mathcal A,z)$:
\begin{align}\label{SAz}
S(\mathcal A, z)&=\sum_{\substack{m\leq n\\ P^+(m)\leq y \\ \mu^2(m_o)=(m,\N)=1}}\sum_{\substack{x_1^2+x_2^2=n-m\\ (x_1x_2,\mathcal Q)=1}}\sum_{c\mid (x_1x_2,P_{\mathcal{Q}}(z))}\mu(c)\\
&\geq \sum_{\substack{m\leq n\\ P^+(m)\leq y \\ \mu^2(m_o)=(m,\N)=1}}\sum_{\substack{x_1^2+x_2^2=n-m}}\sum_{\substack{u_1,u_2\mid \mathcal Q\\ u_i\mid x_i}}\mu(u_1)\mu(u_2)\sum_{c\mid (x_1x_2,P_\mathcal Q (z))}\lambda_c^-\\
&=\sum_{u\mid \mathcal Q}\sum_{c\mid P_\mathcal Q(z)}\lambda_c^-\mu(c)\sum_{\substack{d_1,d_2\\ p\mid d_1d_2\Leftrightarrow p\mid uc}}\mu(d_1)\mu(d_2)\sum_{\substack{m\leq n\\ P^+(m)\leq y\\ \mu^2(m_o)=(m,\N)=1}}r_{d_1,d_2}(n-m).
\end{align}
Our next goal is to evaluate the  innermost sum using Lemma \ref{lem1}. 

\subsection{The cuspidal contribution}

We first treat the contribution  to the $m$-sum in \eqref{SAz} that corresponds to   \eqref{sg3}. To this end we need to show that there exists some absolute constant $\eta>0$ (which eventually will depend on the constants in Theorem \ref{smooth})
such that 
\begin{align}\label{smoothcuspidal}
\sum_{\substack{m\leq n\\ m\equiv n\pmod \Delta\\m \equiv w\pmod {q} \\\mu^2(m_o) = (m, \N)=1\\ P^+(m)\leq y}}\lambda\Big(\frac{n-m}{\Delta}\Big) \ll n^{1-\eta}
\end{align}
uniformly for $\Delta, q\leq n^\eta$ with $(\Delta w, q) = 1$ (specifically, in the notation of \eqref{sg3} we have $\Delta = (\delta\delta_1\delta_2)^2$,  $q = \delta_1'\delta_2'$) and Hecke eigenvalues of a cusp form $\phi$ whose conductor is less than $n^\eta$.  By M\"obius inversion, the previous display equals
\begin{align}
\sum_{\substack{d\leq \sqrt{n}\\(d, 2\N)=1\\ P^+(d)\leq y}}\mu(d)\sum_{\substack{m\leq n/d^2\\ d^2m\equiv n \pmod \Delta\\ d^2m\equiv w\pmod q\\(m,\N)=1\\ P^+(m)\leq y}}\lambda\Big(\frac{n-d^2m}{\Delta}\Big). 
\end{align}
The contribution from $d\geq n^{4\eta}$ can be bounded by 
\begin{align}
\sum_{d\geq n^{4\eta}}\Big(\sum_{\substack{m\leq n/d^2}} 1\Big)^{3/4} \Big(\sum_{\substack{m\leq n/d^2}}\Big|\lambda\Big(\frac{n-d^2m}{\Delta}\Big)\Big|^4\Big)^{1/4}& \ll  n^{\varepsilon}  \sum_{d\geq n^{4\eta}} \frac{n}{d^{3/2}}  \ll n^{1+\varepsilon-2\eta},\end{align}
using the  bound  \eqref{L4bound} for the fourth moment of Hecke eigenvalues and recalling that presently the conductor $C_\phi \leq n^{\eta}$. 

For $d\leq n^{4\eta}$, we detect the condition $(m, \N)=1$ using M\"obius inversion to write 
\begin{align}
\sum_{\substack{ m\leq n/d^2\\ d^2m \equiv n \pmod \Delta\\ d^2m \equiv w \pmod q\\ (m, \N)=1\\P^+(m)\leq y}}\lambda\Big( \frac{n-d^2m}{\Delta}\Big)=\sum_{\substack{\tau\mid \N\\ P^+(\tau)\leq y}}\mu(\tau)\sum_{\substack{ m\leq n/d^2\tau\\ d^2\tau m\equiv n \pmod \Delta \\ d^2\tau m \equiv w\pmod q\\ P^+(m)\leq y}}\lambda\Big( \frac{n-d^2 \tau m}{\Delta}\Big).
\end{align}
The contribution from $\tau \geq n^{4\eta}$ can be bounded by 
\begin{align}
n^{1+\varepsilon} \sum_{\substack{ \tau\mid n \\ \tau \geq n^{4\eta}}}\frac{1}{\tau}\ll \frac{n^{1+\varepsilon} \tau(n) }{n^{4\eta}} \ll n^{1+2\varepsilon-4\eta}  . 
\end{align}
Both error terms are acceptable. Thus, for $d, \tau \leq n^{4\eta}$ and automatically $(d\tau, q) = 1$ and $d_0 := (d^2 \tau, \Delta) \mid n$,  we are left with bounding
$$\sum_{\substack{ m\leq n/d^2\tau\\  m\equiv \frac{n}{d_0}\overline{\frac{d^2\tau}{d_0}} \pmod {\frac{\Delta}{d_0}} \\  m \equiv w \overline{d^2\tau}\pmod q\\ P^+(m)\leq y}}\lambda\Big( \frac{n-d^2 \tau m}{\Delta}\Big),
$$
and the desired power saving as in \eqref{smoothcuspidal} follows (after obvious smoothing) from Theorem \ref{smooth}.

\subsection{Computing the main term I -- summing over arithmetic progressions}\label{part1} We now return to  \eqref{SAz} and insert the contribution corresponding to  \eqref{nl2}. In view of the bounds in the previous subsection, we have for $\mathcal D=n^{\kappa_0}$ for sufficiently small $\kappa_0$ (depending on the constants $\varpi$ in Proposition \ref{SmoothAp} and $\eta$ is as in \eqref{smoothcuspidal})
\begin{align}\label{S}
S(\mathcal{A}, z) \geq \sum_{u\mid \mathcal Q}\sum_{c\mid P_\mathcal Q(z)}\lambda_c^-\mu(c)\sum_{\substack{d_1,d_2\\ p\mid d_1d_2\Leftrightarrow p\mid uc}}\mu(d_1)\mu(d_2) S_{d_1,d_2} + O(n^{1-\kappa_0})
\end{align}
where 
\begin{align}
S_{d_1,d_2}&=\operatornamewithlimits{\sum\sum}_{\substack{\delta=(d_1,d_2)\\ d_i'=d_i/\delta\\ \delta_1\delta_1'=d_1'\\\delta_2\delta_2'=d_2'}}\frac{4}{\delta_1'\delta_2'}\prod_{p\mid \delta_1'\delta_2'}(1-\frac{\chi_{-4}(p)}{p})^{-1}G_{\delta_1'\delta_2'}\sum_{w\in \mathcal G_{\delta_1'\delta_2'}}\sum_{\substack{m\equiv n \pmod{\delta^2 \delta_1^2\delta_2^2}\\ n-m\equiv w\pmod{\delta_1'\delta_2'}\\ P^+(m)\leq y\\ \mu^2(m_o)=(m, \N)=1}}r\Big(\frac{n-m}{\delta^2\delta_1^2\delta_2^2}\Big). 
\end{align}
As in the proof of Proposition \ref{sieve}, we apply Dirichlet's hyperbola method to the $r$-function. Let $\eta_0$ be some small positive constant that will be chosen appropriately later. Let $$x=n(1-1/Z),  \quad  Z=n^{\eta_0}$$ for some $\eta_0>0$ to be chosen later. In particular, in terms of orders of magnitude, $n$ and $x$ can (and will) be used interchangeably. Note that the contribution from $m\geq x$ or $m\leq n-x$ can be bounded by $Z^{-1}n^{1+\varepsilon}\ll n^{1-\eta_0/2}$. Thus we can write $$S_{d_1,d_2}=:S_1+S_2+O(n^{1-\eta_0/2}),$$
where 
\begin{align}
 S_1&=\operatornamewithlimits{\sum\sum}_{\substack{\delta=(d_1,d_2)\\ d_i'=d_i/\delta\\ \delta_1\delta_1'=d_1'\\\delta_2\delta_2'=d_2'}}\frac{4}{\delta_1'\delta_2'}\prod_{p\mid \delta_1'\delta_2'}\Big(1-\frac{\chi_{-4}(p)}{p}\Big)^{-1}G_{\delta_1'\delta_2'}\sum_{w\in \mathcal G_{\delta_1'\delta_2'}}\sum_{\substack{ a\leq \sqrt{n}\\ (a, \delta_1'\delta_2')=1}}\sum_{\substack{n-x\leq m\leq x\\ P^+ (m)\leq y\\ \mu^2(m_o)=(m,\N)=1\\ m\equiv n \pmod {a\delta^2\delta_1^2 \delta_2^2}\\ n-m\equiv w\pmod {\delta_1'\delta_2'}}}\chi_{-4}\Big(\frac{n-m}{a\delta^2\delta_1^2\delta_2^2}\Big),\\
 S_2&=\operatornamewithlimits{\sum\sum}_{\substack{\delta=(d_1,d_2)\\ d_i'=d_i/\delta\\ \delta_1\delta_1'=d_1'\\\delta_2\delta_2'=d_2'}}\frac{4}{\delta_1'\delta_2'}\prod_{p\mid \delta_1'\delta_2'}\Big(1-\frac{\chi_{-4}(p)}{p}\Big)^{-1}G_{\delta_1'\delta_2'}\sum_{w\in \mathcal G_{\delta_1'\delta_2'}}\sum_{\substack{ b\leq \frac{x}{(\delta\delta_1\delta_2)^2\sqrt{n}}\\ (b, \delta_1'\delta_2')=1}}\sum_{\substack{ n-x\leq m\leq x\\ m\leq n-b(\delta\delta_1\delta_2)^2\sqrt{n}\\ P^+ (m)\leq y\\ \mu^2(m_o)=(m,\N)=1\\m\equiv n \pmod {b \delta^2 \delta_1^2 \delta_2^2}\\ n-m\equiv w\pmod{\delta_1'\delta_2'}}}\chi_{-4}(b).
\end{align}
Again we expect that $S_1$ gives a negligible contribution, as in the proof of Proposition \ref{sieve}, due to a particular behaviour at the prime 2 which makes the main term disappear, but this is not so easy to see in the present situation. We will carry out the computation for both $S_1$ and $S_2$ and combine them to a joint main term for which we obtain the desired lower bound. 
As before, the computation of $S_1$ is slightly harder, so we focus on this part in detail, while the computation of $S_2$ is similar, but technically slightly easier. Our first step is to establish the asymptotic evaluations \eqref{S1} and \eqref{S2} below. 

 We would like to choose $\eta_0$ sufficiently small so that we can evaluate the sums by Proposition \ref{SmoothAp} on average of $d_1,d_2$, but this requires some preparation.  We consider $S_1$ first. As before we write $$r_2=(\delta \delta_1\delta_2, 2^\infty), \quad d=\delta\delta_1\delta_2/r_2, \quad  a_2=(a, 2^\infty), \quad a_d=(a, d^\infty), \quad  d'=\delta_1'\delta_2'/(\delta_1'\delta_2',2),$$ so that $a, m$-sum in $S_1$ becomes
\begin{align}\label{sum0}
&\sum_{\substack{a_2\mid 2^\infty\\ (a_2, \delta_1'\delta_2')=1}}\sum_{a_d\mid d^\infty}\chi_{-4}(a_d)\underset{u \, (\text{mod } 4)}{\left.\sum\right.^{\ast}} \chi_{-4}(u )\sum_{\substack{v\pmod{4}}}\chi_{-4}(v)
\sum_{\substack{a\leq \sqrt{n}/a_2a_d\\ (a, 2d_1d_2)=1\\ a\equiv u\pmod 4}}\sum_{\substack{n-x\leq m\leq x ,P^+(m)\leq y\\ \mu^2(m_o)=(m, \mathfrak N)=1\\ m\equiv n\pmod{aa_2a_dr_2^2d^2}\\n-m\equiv w\pmod {d'}\\ 
		\frac{n-m}{a_2r_2^2}\equiv v\pmod {4}
	}}1.
\end{align}
Here we have used the fact that if $2\mid \delta_1'\delta_2'$, then we must have $r_2a_2=1$ and $(w, 2)=1$ so that the condition $n-m\equiv w\pmod {2}$ is incorporated in the $v$-sum.
As before, we can restrict $a_2, a_d\leq n^{\eta_0}$ at the cost of an error of size $O(n^{1+\varepsilon-\eta_0/2})$. The $a,m$-sum in \eqref{sum0} can be written as
\begin{align}\label{sum1}
&\sum_{\substack{ a\leq \sqrt{n}/a_2a_d\\ (a, 2d_1d_2)=1\\ a\equiv u \pmod 4}} \sum_{\substack{ n-x\leq m\leq x, P^+(m)\leq y\\\mu^2(m_o)=(m, \mathfrak N)=1\\m \equiv n \pmod a\\m \equiv n\pmod {a_dd^2}\\ m\equiv n-va_2r_2^2 \pmod{4a_2r_2^2}\\ m\equiv n-w \pmod {d'}}}1.
\end{align}
By the Chinese remainder theorem, the congruence conditions modulo $a_dd^2, 4a_2r_2^2, d'$ can be written as a single congruence condition $\frac{m}{g_1}\equiv c_1\pmod {F_1}$ with $(c_1, F_1)=1$, where 
\begin{equation}\label{f1g1}
F_1=\frac{4a_2r_2^2 a_dd^2d'}{g_1}, \quad g_1=(n-va_2r_2^2, 4a_2r_2^2)(n-w, d')(n, a_dd^2).
\end{equation}
 In order for the $m$-sum to be nonempty, we must have 
 \begin{align}\label{g1conditions}
 P^+(g_1)\leq y,\quad \mu((g_1)_o)^2=(g_1, \mathfrak N)=1,
  \end{align}
 as well as 
 $$\quad P^+((a,n))\leq y, \quad \mu((a,n))^2=((a,n),\mathfrak N)=1,$$ 
   in which case we can re-write the $a, m$-sum in \eqref{sum1} as
\begin{align}\label{Sumam}
\sum_{\substack{ a\leq \sqrt{n}/a_2a_d\\ (a, 2d_1d_2)=1\\ a\equiv u \pmod 4}}\sum_{\substack{(n-x)/g_1\leq m\leq x/g_1\\  P^+(m)\leq y\\ \mu^2(m_o)=(m,  (g_1)_o \mathfrak N)=1 \\ m\equiv n \bar g_1\pmod a\\ m \equiv c_1\pmod {F_1}}}1=\sum_{\substack{h\mid n \\ (h, 2d_1d_2)=1\\ \mu^2(h)=(h, \mathfrak N)=1\\ P^+(h) \leq y}}\sum_{\substack{a\leq \sqrt{n}/a_2a_dh\\ (a, 2d_1d_2)=1\\(a, n/h)=1\\ a\equiv u \bar h \pmod 4}}\sum_{\substack{n-x\leq mg_1h\leq x \\ P^+(m)\leq y\\\mu^2(m_o)=(m_o, g_1h \mathfrak N)=1 \\ m\equiv \frac{n}{h} \bar g_1 \pmod a\\ m \equiv c_1\bar h \pmod {F_1}}}1.
\end{align}
Note that from the construction of the sieve weights we have $d_1, d_2 \leq \mathcal Q\mathcal D\ll n^{\kappa_0}$ and thus we automatically have $g_1 \ll n^{2\eta_0 + 4\kappa_0}$, a small power of $n$. Also note that $g_1 \mid (2d_1d_2)^{\infty}$, so that automatically $(g_1, h) = 1$. 

We can again restrict  $h\leq n^{\eta_0}$ with a total error of size $O(n^{1+\varepsilon-\eta_0/2})$. Finally we detect the conditions $\mu^2(m_o) = 1$ by M\"obius inversion. We summarize that the $a, m$-sum in \eqref{sum0} can, up to an admissible error of size $O(n^{1+\varepsilon-\eta_0/2})$, be replaced with
$$\sum_{\substack{h\mid n \\ (h, 2d_1d_2)=1\\ \mu^2(h)=(h, \mathfrak N)=1\\  P^+(h) \leq y\\ h \leq n^{\eta_0}}}\sum_{\substack{a\leq \sqrt{n}/a_2a_dh\\ (a, 2d_1d_2)=1\\ (a, n/h)=1\\ a\equiv u \bar h\pmod 4}}\sum_{\substack{ n-x  \leq mg_1h\leq x\\\mu^2(m_o)=1\\(m_o, g_1h\mathfrak N)=1\\ P^+(m)\leq y\\   m\equiv \frac{n}{h} \overline{  g_1} \pmod a\\   m\equiv c_1  \overline{ h}\pmod {F_1}}}1,$$
where $a_2, a_d \ll n^{\eta_0}, g_1\ll n^{2\eta_0+4\kappa_0}$ are sufficiently small powers of $n$. 
We evaluate this as  
\begin{align}\label{hklambdaam}
\sum_{\substack{h\mid n \\ (h, 2d_1d_2)=1\\ \mu^2(h)=(h, \mathfrak N)=1\\ P^+(h) \leq y\\ h \leq n^{\eta_0}}}\Big(\sum_{\substack{ a\leq \sqrt{n}/a_2a_dh\\ (a, 2d_1d_2)=1\\ (a,n/h)=1\\ a\equiv u \bar{ h}\pmod 4}}\sum_{\substack{ n-x \leq m  g_1h \leq x\\ P^+(m)\leq y\\ (m, aF_1)=1\\\mu^2(m_o)=1\\(m_o, g_1h\mathfrak N)=1}}\frac{1}{\phi(a)\phi(F_1)}+\tilde{r}_1\Big),
\end{align}
where 
\begin{align}
\tilde{r}_1=\sum_{\substack{a\leq \sqrt{n}/a_2a_dh\\ (a, 2d_1d_2)=1\\ (a, n/h)=1\\ a\equiv u \bar h\pmod 4}}\sum_{\substack{ n-x  \leq mg_1h\leq x\\\mu^2(m_o)=1\\(m_o, g_1h\mathfrak N)=1\\ P^+(m)\leq y\\   m\equiv \frac{n}{h} \overline{  g_1} \pmod a\\   m\equiv c_1  \overline{ h}\pmod {F_1}}}1-\sum_{\substack{ a\leq \sqrt{n}/a_2a_dh\\ (a, 2d_1d_2)=1\\ (a,n/h)=1\\ a\equiv u \bar{ h}\pmod 4}}\sum_{\substack{ n-x \leq m  g_1h \leq x\\ P^+(m)\leq y\\ (m, aF_1)=1\\\mu^2(m_o)=1\\(m_o, g_1h\mathfrak N)=1}}\frac{1}{\phi(a)\phi(F_1)}.
\end{align}
We want to show that $\tilde{r}_1$ that is small on average using Proposition \ref{SmoothAp} with $\mathbf{q}=a, \mathbf{a}_1=\frac{n}{h}, \bd=F_1$. It is not enough to apply Proposition \ref{SmoothAp} with $\mathbf{a}_2=g_1$ and sum over $g_1$ directly due to sparseness of smooth numbers when $y$ is small and thus we need to be more precise about $g_1$.
Note that $a_2r_2^2a_dd^2d'\mid g_1F_1$ and so $(h,d_1d_2)=1$ and $(a, 2d_1d_2)$ is equivalent to $(a, g_1F_1)=1$. The coefficients of $h$ can be bounded by $1$. The coefficients of $g_1, F_1$ is bounded by $\tau(dd')^C\ll \tau(g_1F_1)^C$. From \eqref{f1g1}, we see that if $p\mid g_1, p\nmid F_1$, then we must have either $p\mid (n-w,d')$ or $p\mid 2n$ and $p^2\mid g_1$. If $p^2\mid g_1$, we also have $p\mid 2n$. Therefore we can write $g_1=g_0g_1'g_3$ where $g_1'=(n-w, d')$ is square-free, $g_0\mid (2n)^\infty$ is square-full and $g_3\mid F_1$. We can then change the order of summation to move $h, g_0,g_1'$ outside the summation over $u, c$ or equivalently $d_1,d_2$.  
After dividing the $m$-sum into short intervals to separate the variables $m$ from $g_1$ similarly as in the proof of Proposition \ref{sieve}, we can apply Proposition \ref{SmoothAp} for the $c,a_2,a_d,a,m$-sum with $\mathbf{d}=F_1, \mathbf{a}_2=g_0g_1', \mathbf{a}_3=g_3$ for fixed $h, g_0g_1'$ (with $h\leq n^{\eta_0}$ and $g_0g_1'\ll n^{2\eta_0+4\kappa_0}$ for $\eta_0, \kappa_0$ sufficiently small)
so that the total contribution of the error $\tilde{r}_{1}=\tilde{r}_1(g_3)$ in \eqref{S} can be bounded by 
\begin{align}
\mathcal E_1
 & \ll \sum_{\substack{h\mid n, h\nmid \mathfrak N\\ h\leq n^{\eta_0} \\ P^+(h)\leq y}}\sum_{\substack{g_1\leq n^{2\eta_0+4\kappa_0}\\P^+(g_1)\leq y\\g_1=g_0g_1'g_3\\g_0\mid 2^\infty n\\ p\mid g_0\Rightarrow p^2\mid g_0 }}\sum_{\substack{F_1\ll n^{4\kappa_0}\\ g_3\mid F_1}}\tau(g_1F_1)^C\sum_{\substack{d'\mid F_1g_1'\\ \mu^2(d')=1}}\frac{\tau(d')}{d'}\sum_{\substack{w\in \mathcal G_{d'}\\ (n-w, d')=g_1'}}|\tilde{r}_1(g_3)|\\
 &\ll \sum_{\substack{h\mid n, h\nmid \mathfrak N\\ h\leq n^{\eta_0} \\ P^+(h)\leq y}}\sum_{\substack{g_0g_1'\leq n^{2\eta_0+4\kappa_0}\\P^+(g_0g_1')\leq y\\g_0\mid 2^\infty n\\ p\mid g_0\Rightarrow p^2\mid g_0 }}\tau(g_0g_1')^C\sum_{\substack{F_1\ll n^{4\kappa_0}\\ }}\tau(F_1)^{2C}\max_{g_3\mid F_1}\sum_{\substack{d'\mid F_1g_1'\\ \mu^2(d')=1}}\frac{\tau(d')}{d'}\sum_{\substack{w\in \mathcal G_{d'}\\ (n-w, d')=g_1'}}|\tilde{r}_{1}(g_3)|.
\end{align}
Note that 
\begin{align}
\sum_{\substack{d'\mid F_1g_1'\\\mu^2(d')=1}}\frac{\tau(d')}{d'}\sum_{\substack{w\in \mathcal G_{d'}\\ (n-w, d')=g_1'}}1\ll\frac{\tau(F_1g_1')^2}{g_1'}  
\end{align}
and thus for $y\geq\mathcal L(x) = (\log x)^{D}$ for some sufficiently large $D$, Proposition \ref{SmoothAp} gives
\begin{align}
 \mathcal E_1
 & \ll_C\sum_{\substack{h\mid n, h\nmid \mathfrak N\\ h\leq n^{\eta_0}\\ P^+(h)\leq y}}\sum_{\substack{g_0g_1'\leq n^{2\eta_0+4\kappa_0}\\ P^+(g_0g_1')\leq y\\ g_0\mid 2^\infty n\\p\mid g_0\Rightarrow p^2\mid g_0}}\frac{\tau(g_0g_1')^C}{g_1'}\Psi_{hg_0g_1'\mathfrak N}\Big(\frac{x}{g_0g_1'h},y\Big)(\log x)^{-A}.
\end{align}

We recall \eqref{psix/d}, \eqref{alphaasym} and note that for $y\geq \mathcal L(x)$ we have $\alpha \geq 1-\frac{1}{D}+o(1)$.  
Thus when $y\geq \mathcal L(x)$ we have
\begin{align}
\mathcal E_1& \ll_C\sum_{\substack{h\mid n, h\nmid \mathfrak N\\ h\leq n^{\eta_0}\\ P^+(h)\leq y}}\sum_{\substack{g_0g_1'\leq n^{2\eta_0+4\kappa_0}\\ P^+(g_0g_1')\leq y\\ g_0\mid 2^\infty n\\p\mid g_0\Rightarrow p^2\mid g_0}}\frac{\tau(g_0g_1')^C}{g_1'}\Psi_{hg_0g_1'\mathfrak N}\Big(\frac{x}{g_0g_1'h},y\Big)(\log x)^{-A}\\
&\ll \sum_{\substack{h\mid n , h\nmid \mathfrak N\\P^+(h)\leq y}}\frac{1}{h^\alpha}\prod_{\substack{p\mid \mathfrak N\\ p\leq y}}(1-p^{-\alpha})\sum_{\substack{g_0\mid 2^\infty n\\ p\mid g_0\Rightarrow p^2\mid g_0}}\frac{\tau(g_0)^C}{g_0^\alpha}\sum_{g_1'\leq n^{2\eta_0+4\kappa_0}}\frac{\tau(g_1')^C}{(g_1')^{1+\alpha}}\Psi(x,y)(\log x)^{-A}\\
& \ll \prod_{\substack{p\mid n\\ p\leq y}}\Big(1+\frac{\chi_{-4}(p)}{p^\alpha}\Big)\Psi(x, y)(\log x)^{-A} 
\end{align}
by choosing $A$ large enough in terms of $C$. Since $\prod_{p\mid n}(1+\frac{1}{p})\ll \log_2 n$, we have $$\mathcal E_1\ll \mathcal F(n, y)\Psi(x,y)(\log x)^{-A}$$ for any $A>0$ with $\mathcal F(n, y)$ defined as \eqref{Fdef}.

For the main term, we can  substitute the main term in \eqref{hklambdaam} back into \eqref{sum0} and evaluating the $u$-sum, we obtain
  (writing the condition \eqref{g1conditions} on $g_1$ in terms of the remaining variables and noting that $\mu^2(d')=1, a_d\mid d^\infty$) that we can replace $S_1$ by $\tilde{S}_1$  in \eqref{S} up to a negligible error  $\mathcal E_1$ where 
  \begin{equation}\label{S1}
  \tilde{S}_1= \sum_{\delta_1\delta_1'=d_1'}\sum_{\delta_2\delta_2'=d_2'}\frac{4}{\delta_1'\delta_2'}\prod_{p\mid \delta_1'\delta_2'}\Big(1-\frac{\chi_{-4}(p)}{p}\Big)^{-1}G_{\delta_1'\delta_2'}\mathds{1}_{(d, \mathfrak N)=1}\mathcal M_1,
  \end{equation}
  with 
  \begin{align}
  \mathcal M_1=&\sum_{\substack{a_2\mid 2^\infty\\ (a_2, \delta_1'\delta_2')=1\\ a_2\leq n^{\eta_0} }} \sum_{\substack{a_d\mid d^\infty,a_d\leq n^{\eta_0}\\P^+((n, a_dd^2))\leq y\\ \mu^2((n, a_dd^2))=1}} \chi_{-4}(a_d) \sum_{\substack{v\pmod{4}}}\chi_{-4}(v) 
  \\ &\times \sum_{\substack{w\pmod {\delta_1'\delta_2'}\\ (w,\delta_1'\delta_2')=1\\w\in \mathcal G_{\delta_1'\delta_2'}\\ P^+((n-w,d'))\leq y\\  (n-w,d',\N)=1}}\sum_{\substack{h\mid n\\ (h, 2d_1d_2\mathfrak N)=1\\P^+(h)\leq y\\ \mu^2(h)=1\\h\leq n^{\eta_0}}}\chi_{-4}(h)\sum_{\substack{a\leq \sqrt{n}/a_2a_dh\\ (a, 2d_1d_2)=1\\ (a, n/h)=1}}\sum_{\substack{n-x \leq mg_1h \leq x\\ P^+(m)\leq y\\ (m,aF_1(g_1)_oh\N )=1\\\mu^2(m_o)=1 }}\frac{\chi_{-4}(a)}{\phi(a)\phi(F_1)}.
  \end{align}
Next we change the summation order of $m$ and $a$, evaluate the $a$-sum by \cite[Lemma 5.2]{ABL} and then complete the $a_2, a_d,h$-sums with an admissible error to finally obtain (with e.g.\ $\eta_0=10\kappa_0$), uniformly for $d_1,d_2\ll n^{\kappa_0}$ (with $\kappa_0$ sufficiently small), that
  \begin{align}
  \mathcal M_1=&\frac{\pi}{4}\prod_{p}\Big(1+\frac{\chi_{-4}(p)}{p(p-1)}\Big)\sum_{\substack{a_2\mid 2^\infty\\ (a_2, \delta_1'\delta_2')=1 }} \sum_{\substack{a_d\mid d^\infty\\P^+((n, a_dd^2))\leq y\\ \mu^2((n, a_dd^2))=1}} \chi_{-4}(a_d) \sum_{\substack{v\pmod{4}}}\chi_{-4}(v) 
  \\ &\times  \sum_{\substack{w\pmod {\delta_1'\delta_2'}\\  w\in \mathcal G_{\delta_1'\delta_2'}\\ P^+((n-w,d'))\leq y\\  (n-w,d',\N)=1}}\sum_{\substack{h\mid n\\ (h, 2d_1d_2\N)=1\\P^+(h)\leq y\\ \mu^2(h)=1}}\chi_{-4}(h)\sum_{\substack{  m\leq x/g_1h\\ P^+(m)\leq y\\ (m,F_1(g_1)_oh\N )=1\\\mu^2(m_o)=1\\ }}\frac{\mathfrak{c}(2d_1d_2mn/h)}{ \phi(F_1)}  + O(n^{1 - 3\kappa_0})
  \end{align}
with $\mathfrak{c}(f)$ as in \eqref{cf}.\\

   The same strategy can be applied to $S_2$. We give a sketch of the manipulations highlighting the differences from those for $S_1$. To begin with, we do not need to separate $r_2$ from $\delta\delta_1\delta_2$, so that the $b,m$-sum in $S_2$ becomes
   \begin{align}\label{bm-sum}
   	&\sum_{b_d\mid (\delta\delta_1\delta_2)^\infty}\chi_{-4}(b_d)\underset{u \, (\text{mod } 4)}{\left.\sum\right.^{\ast}} \chi_{-4}(u )
   	\sum_{\substack{b\leq \frac{x}{b_d(\delta\delta_1\delta_2)^2\sqrt{n}}\\ (b, d_1d_2)=1\\ b\equiv u\pmod 4}}\sum_{\substack{n-x\leq m\leq x, P^+(m)\leq y\\ m\leq n-bb_d(\delta\delta_1\delta_2)^2 \sqrt{n}\\  \mu^2(m_o)=(m, \mathfrak N)=1\\ m\equiv n\pmod{bb_d(\delta \delta_1\delta_2)^2}\\n-m\equiv w\pmod{\delta_1'\delta_2'}}}1.
   \end{align}
   We can restrict $b_d\leq n^{\eta_0}$ with an error of size $O(n^{1+\varepsilon-\eta_0/2})$. After splitting $m$ into residue classes modulo $d_1d_2$, using the Chinese remainder theorem to rewrite the congruence condition modulo $b_d(\delta\delta_1\delta_2)^2\delta_1'\delta_2'$ as $m\equiv c_2\pmod{F_2}$ with 
   \begin{align}\label{f2g2}
   F_2=\frac{b_d(\delta\delta_1\delta_2)^2\delta_1'\delta_2'}{g_2},\  g_2=(n-w,\delta_1'\delta_2')(n, b_d\delta^2\delta_1^2\delta_2^2)
   \end{align} 
   and then extracting $(b,n)$, we can write the $b,m$-sum in  \eqref{bm-sum} as 
   \begin{align}
  \mathds{1}_{ \substack{P^+(g_2)\leq y\\(g_2,\mathfrak N)=1\\\mu^2((g_2)_o)=1}}\sum_{\substack{h\mid n \\ (h, d_1d_2\mathfrak N)=1\\ P^+(h)\leq y \\ \mu^2(h)=1}}\chi_{-4}(h)
   \sum_{\substack{ b\leq \frac{x}{b_d(\delta\delta_1\delta_2)^2\sqrt{n}h}\\ (b, d_1d_2)=1\\(b, n/h)=1\\ b\equiv u\pmod 4}}\sum_{\substack{n-x\leq mg_2h\leq x, P^+(m)\leq y\\ mg_2h\leq n-bb_d(\delta\delta_1\delta_2)^2 \sqrt{n}\\  \mu^2(m_o)=(m_o, g_2h\mathfrak N)=1\\ m\equiv \frac{n}{h}\overline{g}_2\pmod{b}\\ 
   		m\equiv c_2\overline{h}\pmod {F_2}}}1.
   \end{align}
         We can again restrict $h\leq n^{\eta_0}$ with a total error of size $O(n^{1+\varepsilon-\eta_0/2})$. Additionally, we 
   need to separate $b,m$ in the summation condition as well as their dependency on $g_2,F_2$ in order to apply Proposition \ref{SmoothAp}. As usual, this can be achieved by splitting $m$ into intervals of the form $[M, (1+\Delta)M]$ with $\Delta=R^{-\eta_0}$ so that at most one interval interferes with the condition $mg_2h\leq n-bb_d\sqrt{n}(\delta\delta_1\delta_2)^2$, in which case the  contribution can be bounded by $O(n^{1+\varepsilon}\Delta)$ trivially. Here $R=\exp(\sqrt{\log x})$. Similarly we can split the $b$-sum into intervals of the form $[M, (1+\Delta)M]$ so that the condition $b\leq \frac{x}{b_d(\delta\delta_1\delta_2)^2\sqrt{n}h}$ only interferes with one interval in which case the contribution is also negligible. After M\"obius inversion to detect $\mu^2(m_o)=1$, we can finally apply Proposition \ref{SmoothAp} so that the $b,m$-sum in the last display can be evaluated as
   \begin{align}\label{bmsum}
   \sum_{\substack{ b\leq \frac{x}{b_d(\delta\delta_1\delta_2)^2\sqrt{n}h}\\ (b, d_1d_2)=1\\(b,n/h)=1\\ b\equiv u\pmod 4}}\sum_{\substack{n-x\leq mg_2h\leq x, P^+(m)\leq y\\ mg_2h\leq n-bb_d(\delta\delta_1\delta_2)^2 \sqrt{n}\\  \mu^2(m_o)=(m, \mathfrak N)=1\\ (m, bF_2(g_2)_oh\mathfrak N)=1}}\frac{1}{\phi(b)\phi(F_2)}+ \tilde{r}'_2.
   \end{align}  
   Here the total contribution from the error terms $\tilde{r}'_2$ in \eqref{S} can be bounded using Proposition \ref{SmoothAp} the same way as in $\mathcal E_1$ so that \begin{align}
   \mathcal E_2:=&\sum_{u\mid \mathcal Q}\sum_{c\mid P_{\mathcal Q}(z)}\lambda_c^-|\mu(c)|\sum_{\substack{d_1, d_2\\ p\mid d_1d_2\Leftrightarrow p\mid uc}}|\mu(d_1)\mu(d_2)|\operatornamewithlimits{\sum\sum}_{\substack{\delta=(d_1,d_2)\\ d_i'=d_i/\delta\\ \delta_1\delta_1'=d_1'\\\delta_2\delta_2'=d_2'}}\frac{4}{\delta_1'\delta_2'}\prod_{p\mid \delta_1'\delta_2'}\Big(1-\frac{\chi_{-4}(p)}{p}\Big)^{-1}G_{\delta_1'\delta_2'}\sum_{w\in \mathcal G_{\delta_1'\delta_2'}}
   |\tilde{r}_2'|\\
   &\ll_C \prod_{\substack{p\mid n\\ p\leq y}}\Big(1+\frac{\chi_{-4}(p)}{p^\alpha}\Big)\Psi(x,y)(\log x)^{-A}\ll_C\mathcal F(n, y)\Psi(x,y)(\log x)^{-A+1}
   \end{align}
   for any $A>0$.  
 After changing the order of summation, evaluating the $b$-sum in \eqref{bmsum} by \cite[Lemma 5.2]{ABL} (note that the upper limit for $b_d$ is always $\gg n^{1/2-\eta_0}$ unless $|n-mg_2h|\leq n^{1-\eta_0}$ whose contribution can be bounded by $O(n^{1-\eta_0/2})$) and completing the $b_d,h$-sums, we obtain, uniformly for $d_1,d_2\leq n^{\kappa_0}$, that we can replace $S_2$ by $\tilde{S}_2$ in \eqref{S} up to a negligible error $\mathcal E_2$ where  
   \begin{align}\label{S2}
   \tilde{S}_2=\sum_{\delta_1\delta_1'=d_1'}\sum_{\delta_2\delta_2'=d_2'}\frac{4}{\delta_1'\delta_2'}\prod_{p\mid \delta_1'\delta_2'}(1-\frac{\chi_{-4}(p)}{p})^{-1}G_{\delta_1'\delta_2'}\mathds{1}_{(\delta\delta_1\delta_2, \mathfrak N)=1}\mathcal M_2
   \end{align}
  and
    \begin{align}
    &\mathcal M_2=\frac{\pi}{4}\prod_{p}\Big(1+\frac{\chi_{-4}(p)}{p(p-1)}\Big) \sum_{\substack{b_d\mid(\delta\delta_1\delta_2)^\infty\\P^+((n, b_d(\delta\delta_1\delta_2)^2))\leq y\\ \mu^2((n, b_d(\delta\delta_1\delta_2)^2)_o)=1}} \chi_{-4}(b_d) 
    \\ &\times \sum_{\substack{w\pmod {\delta_1'\delta_2'}\\ w\in \mathcal{G}_{\delta_1'\delta_2'}\\ P^+((n-w,\delta_1'\delta_2'))\leq y\\  (n-w,\delta_1'\delta_2',\N)=1}}\sum_{\substack{h\mid n\\ (h, d_1d_2)=1\\P^+(h)\leq y\\ \mu^2(h)=1\\(h, \N)=1}}\chi_{-4}(h)\sum_{\substack{m\leq x/g_2h\\ P^+(m)\leq y\\ (m,F_2(g_2)_oh\N )=1\\\mu^2(m_o)=1\\ }}\frac{\mathfrak{c}(d_1d_2mn/h)}{ \phi(F_2)} + O(n^{1-3\kappa_0}).
    \end{align}

    \subsection{Computing the main term II -- Euler products} 
   We substitute \eqref{S1} and \eqref{S2} back into \eqref{S}.   By choosing $\mathcal D=n^{\kappa_0}$ sufficiently small, it follows from \eqref{SAz} together with \eqref{smoothcuspidal} that
    \begin{equation}\label{saz}
    S(\mathcal A, z) \geq  \pi \prod_{p}\Big(1+\frac{\chi_{-4}(p)}{p(p-1)}\Big) (S_1^-+S_2^-) +O(\mathcal E_1+\mathcal E_2)+O(n^{1-\kappa_0}),
    \end{equation}
  where 
  \begin{align}
  S_1^-&=\sum_{u\mid \mathcal Q}\sum_{\substack{c\mid P_\mathcal Q(z)}}\lambda_{c}^-\mu(c)\sum_{\substack{d_1,d_2\\ p\mid d_1d_2\Leftrightarrow p\mid uc}}\mu(d_1)\mu(d_2)\operatornamewithlimits{\sum\sum}_{\substack{\delta=(d_1,d_2)\\\delta_i\delta_i'=d_i'=d_i/\delta\\ r_2d=\delta\delta_1\delta_2, r_2\mid 2, (d,2\mathfrak N)=1\\ d'=\delta_1'\delta_2'/(\delta_1'\delta_2',2)}}\frac{1}{\delta_1'\delta_2'}\prod_{p\mid \delta_1'\delta_2'}\Big(1-\frac{\chi_{-4}(p)}{p}\Big)^{-1}\\
  &\times \sum_{\substack{a_2\mid 2^\infty\\ (a_2, \delta_1'\delta_2')=1  }}\sum_{\substack{a_d\mid d^\infty\\P^+((n, a_dd^2))\leq y\\ \mu^2((n, a_dd^2))=1}} \chi_{-4}(a_d)  \sum_{\substack{v\pmod{4}}}\chi_{-4}(v) \\
  &
 \times G_{d'}\sum_{\substack{w\pmod {d'}\\ w\in \mathcal G_{d'}\\ P^+((n-w,d'))\leq y\\  (n-w,d',\N)=1}}\sum_{\substack{h\mid n\\ (h, 2d_1d_2\N)=1\\P^+(h)\leq y\\\mu^2(h)=1}}\chi_{-4}(h) \sum_{\substack{m\leq x/g_1h\\ P^+(m)\leq y\\ (m, F_1(g_1)_oh\N )=1\\ \mu^2(m_o)=1}}\frac{ \mathfrak{c}(2d_1d_2mn/h)}{ \phi(F_1)}
  \end{align}
  with $F_1, g_1$ as in \eqref{f1g1}, and 
  \begin{align}
  S_2^-&=\sum_{u\mid \mathcal Q}\sum_{c\mid P_\mathcal Q(z)}\lambda_{c}^-\mu(c)\sum_{\substack{d_1,d_2\\ p\mid d_1d_2\Leftrightarrow p\mid uc}}\mu(d_1)\mu(d_2)\operatornamewithlimits{\sum\sum}_{\substack{ \delta=(d_1,d_2)\\\delta_i\delta_i'=d_i'=d_i/\delta\\ d=\delta\delta_1\delta_2, (d, \mathfrak N)=1\\ d'=\delta_1'\delta_2'}}\frac{1}{\delta_1'\delta_2'}\prod_{p\mid \delta_1'\delta_2'}\Big(1-\frac{\chi_{-4}(p)}{p}\Big)^{-1}\\
  &\times G_{d'}\sum_{\substack{b_d\mid d^\infty\\P^+((n, b_dd^2))\leq y\\ \mu^2((n, b_dd^2)_o)=1}} \chi_{-4}(b_d)  
  \sum_{\substack{w\pmod {d'}\\  w\in \mathcal G_{d'}\\ P^+((n-w,d'))\leq y\\ (n-w,d',\N)=1}}\sum_{\substack{h\mid n\\ (h, d_1d_2\N)=1\\P^+(h)\leq y\\\mu^2(h)=1}}\chi_{-4}(h)\sum_{\substack{ m \leq x/g_2h\\ P^+(m)\leq y\\(m, F_2(g_2)_oh\N)=1\\ \mu^2(m_o)=1}}\frac{\mathfrak c(d_1d_2mn/h)}{\phi(F_2)}
  \end{align}
  with $F_2, g_2$ as in \eqref{f2g2}.
Now we rename $mg_i$ as $m$, change the order of summation and evaluate $S_1^-, S_2^-$ in roughly the form $$\sum_{\substack{h\mid n\\ (\cdots)}}\sum_{\substack{m\leq x/h\\(\cdots)}}\sum_{\substack{u\mid \mathcal Q\\ (u, h)=1}}\sum_{\substack{c\mid P(z)\\ (c, \mathcal Qh)=1}}\lambda_c^-w_{m,n,h}(c)$$
for some multiplicative function $w_{m,n,h}(c)$. In particular, we postpone the summation over $m$ to the last moment. Precisely, we write
\begin{equation}\label{s1-}
\begin{split}
S_1^-&=\sum_{\substack{h\mid n\\ P^+(h)\leq y\\\mu^2(h)=1\\(h,\N)=1}}\chi_{-4}(h)\sum_{\substack{m\leq x/h\\(m,\N)=1\\ P^+(m)\leq y\\ \mu^2(m_oh)=1}}\sum_{\substack{u\mid \mathcal Q\\ (u,h)=1}}\sum_{\substack{c\mid P(z)\\ (c,\mathcal Qh)=1}}\lambda_{c}^-\mu(c)\sum_{\substack{d_1,d_2\\p\mid d_1d_2\Leftrightarrow p\mid uc}}\mu(d_1)\mu(d_2)\mathfrak c\Big(2d_1d_2\frac{n}{h}m\Big)
\\&\operatornamewithlimits{\sum\sum}_{\substack{\delta=(d_1,d_2),\delta_i\delta_i'=d_i'=d_i/\delta\\ r_2d=\delta\delta_1\delta_2, r_2\mid 2, (d,2)=1\\ d'=\delta_1'\delta_2'/(\delta_1'\delta_2',2)}}\frac{G_{d'}}{\delta_1'\delta_2'}\prod_{p\mid \delta_1'\delta_2'}\Big(1-\frac{\chi_{-4}(p)}{p}\Big)^{-1}
\operatornamewithlimits{\sum\sum\sum\sum}_{\substack{ a_2,a_d, v,w \\ (a_2,a_d, v,w)\in \mathcal R_1}}\frac{\chi_{-4}(a_d)\chi_{-4}(v)}{\phi(F_1)}
\end{split}
\end{equation} 
  where $\mathfrak c( f)$ is as in \eqref{cf}, $F_1$ is as in \eqref{f1g1}, and $$\mathcal R_1=\left\{ 
  \begin{array}{cc} a_2\mid 2^\infty ,a_d\mid d^\infty\\  v\pmod{4},w\pmod{d'} \end{array}: 
  \begin{array}{cc}
  (a_2,\delta_1'\delta_2')=1\\ (m, 4a_2r_2^2)=(n-va_2r_2^2,4a_2r_2^2)\\ (m, a_dd^2)=(n, a_dd^2)\end{array} \begin{array}{cc}
  w\in \mathcal G_{d'}\\ (w, d')=1\\(m,d')=(n-w, d')\end{array}\right\}.$$ 
 To see this, we recall the conditions \eqref{g1conditions} and observe that after replacing $mg_1$ with $m$, the new conditions $P^+(m)\leq y$, $\mu^2(m_o)=1$, $(m, \mathfrak N)=1$ take care of the old conditions $P^+((n,a_dd^2))\leq y$,  $P^+((n-w, d'))\leq y$, $\mu^2((n,a_dd^2))=(n, a_dd^2, \mathfrak N)=(n-w, d', \mathfrak N)=1$. Moreover, the old condition $(m, F_1)=1, mg_1\leq x/h$ is equivalent to the new condition $(m, 4a_2r_2^2a_dd^2d')=g_1, m\leq x/h$ after replacing $mg_1$ with $m$.
 Note also that $\mathcal G_{\delta_1'\delta_2'}\simeq\mathcal G_{d'}$. 
 Similarly, we write \begin{equation}\label{s2-def}
  \begin{split}
  S_2^-&=\sum_{\substack{h\mid n\\ P^+(h)\leq y\\\mu^2(h)=1\\(h,\N)=1}}\chi_{-4}(h)\sum_{\substack{m\leq x/h\\(m,\N)=1\\ P^+(m)\leq y\\ \mu^2(m_oh)=1}}\sum_{\substack{u\mid \mathcal Q\\ (u,h)=1}}\sum_{\substack{c\mid P(z)\\ (c,\mathcal Qh)=1}}\lambda_{c}^-\mu(c)\sum_{\substack{d_1,d_2\\p\mid d_1d_2\Leftrightarrow p\mid uc}}\mu(d_1)\mu(d_2)\mathfrak c\Big(d_1d_2\frac{n}{h}m\Big)
  \\&\operatornamewithlimits{\sum\sum}_{\substack{\delta=(d_1,d_2)\\\delta_i\delta_i'=d_i'=d_i/\delta\\ d=\delta\delta_1\delta_2\\ d'=\delta_1'\delta_2'}}\frac{G_{d'}}{\delta_1'\delta_2'}\prod_{p\mid \delta_1'\delta_2'}\Big(1-\frac{\chi_{-4}(p)}{p}\Big)^{-1}\operatornamewithlimits{\sum\sum}_{\substack{ b_d, w \\ (b_d,w)\in \mathcal R_2}}\frac{\chi_{-4}(b_d)}{\phi(F_2)}
  \end{split}
  \end{equation} 
  where $F_2$ is as in \eqref{f2g2} and 
$$\mathcal R_2=\left\{ 
\begin{array}{cc} b_d\mid d^\infty\\ w\pmod{d'} \end{array}:
\begin{array}{cc}
(m, b_dd^2)=(n, a_dd^2)\\
w\in \mathcal G_{d'}
\end{array} \begin{array}{cc}\\
(w, d')=1\\(m,d')=(n-w, d')\end{array}\right\}.$$ 

Our next aim is to write the   second line  in \eqref{s1-} as an Euler product which will be given in \eqref{secondline} below. The  $v$-sum modulo $4$ in \eqref{s1-} contributes
\begin{align}\label{vmod4}
\sum_{\substack{ v\pmod 4\\ (n-va_2r_2^2, 4a_2r_2^2)=(m, 4a_2r_2^2)}}\frac{\chi_{-4}(v)}{\phi(\frac{4a_2r_2^2}{(n-va_2r_2^2, 4a_2r_2^2)})}=\mathds{1}_{2a_2r_2^2\mid m}\mathds{1}_{a_2r_2^2\parallel n}\chi_{-4}\Big(\frac{n}{a_2r_2^2}\Big)\times\begin{cases}
-1 & 4a_2r_2^2\nmid m,\\
1 & 4a_2r_2^2\mid m.
\end{cases}
\end{align}
In particular, we see that for each $r_2$, only one $a_2$ with $(a_2, \delta_1'\delta_2')=1$ gives a non-zero contribution.
The $w$-sum together with $G_{d'}$ can be evaluated as (treat each prime $p 
\mid d'$ at a time)
\begin{equation}\label{wsum}
G_{d'}\sum_{\substack{w\pmod {d'}\\ w\in \mathcal G_{d'}\\ (n-w, d')=(m, d')}}\frac{1}{\phi\big(d'/(n-w, d')\big)}=\mathds{1}_{(d',m,n)=1}
\prod_{\substack{p\mid d'\\  p\mid mn}}(1+\chi_{p^*}(n))\prod_{\substack{p\mid d'\\ p\nmid mn}}\Big(1-\frac{1+\chi_{p^*}(n)}{p-1}\Big).
\end{equation}
Finally, the sum over $a_d\mid d^\infty$ gives
\begin{align}\label{adsum}
&\sum_{\substack{ a_d\mid d^\infty \\ (n, a_dd^2)=(m, a_dd^2)}}\frac{\chi_{-4}(a_d)}{\phi(a_dd^2/(n, a_dd^2))} = \prod_{p \mid d}  \mathfrak d(p;m,n)
\end{align}
where
\begin{align}\label{dpmn}
\mathfrak d(p;m,n)=\begin{cases}
\frac{1}{(p-1)(p-\chi_{-4}(p))}, & p\nmid mn,\\
\frac{p}{(p-1)(p-\chi_{-4}(p))}, & \nu_p(m)=\nu_p(n)=1,\\
 0, & \text{else}.
\end{cases}
\end{align}
Here we use crucially that (the odd part of) $m$ is square-free. For more general $m$, the formula would look much more complicated, in particular the contribution could be much larger which would create considerable technical problems in applying Lemma \ref{lowerboundsieve}. 

We summarize: the second line in  \eqref{s1-}
can be written as an Euler product
\begin{align}\label{secondline}
&\sum_{\substack{\delta_1' \delta_1 =  d_1'\\ \delta_2'\delta_2 = d_2'}}\frac{G_{\delta_1'\delta_2'}}{\delta_1'\delta_2'}\prod_{p\mid \delta_1'\delta_2'}\Big(1-\frac{\chi_{-4}(p)}{p}\Big)^{-1}\sum_{a_2,a_d,v,w}(\dots)=:\prod_{p^{\nu}\parallel d_1d_2}\mathfrak f_1(p,\nu;m,n)
\end{align}
where the Euler factor is composed of \eqref{vmod4}, \eqref{wsum} and \eqref{adsum}. 

We  compute $\mathfrak f_1(2,\nu;m,n)$ explicitly using \eqref{vmod4}.
If $2\nmid d_1d_2$, then $r_2=1$ and only $a_2=2^{\nu_2(n)}$ contributes. If $2\parallel d_1d_2$, then we have two cases: when $2\parallel \delta_1'\delta_2'$, we have $r_2=1$ and only $a_2=1$ contributes (since we require $(a_2, \delta_1'\delta_2')=1$); and when $2\parallel \delta_1\delta_2$, we have $r_2=2$ and only $a_2=2^{\nu_2(n)-2}$ contributes. If $2\mid (d_1,d_2)$, then $r_2=2$ and only $a_2=2^{\nu_2(n)-2}$ contributes. Precisely, we have 
\begin{equation}\label{f1}
\mathfrak f_1(2,\nu;m,n)=\begin{cases}
\mathds{1}_{2^{\nu_2(n)+1}\mid m}\chi_{-4}(n_o)(-1)^{\mathds{1}_{2^{\nu_2(n)+2}\nmid m}},& \nu=0,\\
\frac{1}{2} \mathds{1}_{2\mid m}\chi_{-4}(n)(-1)^{\mathds{1}_{4\nmid m}}+ \mathds{1}_{2^{\nu_2(n)+1}\mid m}\mathds{1}_{4\mid n}\chi_{-4}(n_o)(-1)^{\mathds{1}_{2^{\nu_2(n)+2}\nmid m}},& \nu=1 ,\\ 
\mathds{1}_{2^{\nu_2(n)+1}\mid m}\mathds 1_{4\mid n}\chi_{-4}(n_o)(-1)^{\mathds{1}_{2^{\nu_2(n)+2}\nmid m}},& \nu=2.
\end{cases}
\end{equation}
We can also compute $\mathfrak f_1(p,\nu;m,n)$ for $2\nmid p$ using \eqref{wsum} and \eqref{adsum} as 
\begin{equation}\label{r1}
\mathfrak f_1(p, \nu;m,n)=\begin{cases}
1, & \nu=0,\\
\frac{1}{(p-1)(p-\chi_{-4}(p))}+\frac{1}{p}(1-\frac{\chi_{-4}(p)}{p})^{-1}(1-\frac{1+\chi_{p^*}(n)}{p-1}), & p\nmid mn,\nu=1,\\
\frac{1}{p}(1-\frac{\chi_{-4}(p)}{p})^{-1}\big(1+\chi_{p^*}(n)\big), & p\mid mn, p\nmid (m,n),\nu=1,\\
\frac{p}{(p-1)(p-\chi_{-4}(p))},
& p\mid (m,n),  \nu= 1,\\
\mathfrak d(p;m,n), & \nu=2.
\end{cases}
\end{equation}
\\

For $S_2^-$, we have analogous  results with $\mathfrak f_1$ replaced by $\mathfrak f_2$  where
\begin{equation}\label{f2}
\mathfrak f_2(2,\nu;m,n)=\begin{cases}
1,  & \nu=0, \\
\frac{1}{2}\mathds{1}_{2\mid mn, 2\nmid (m,n)}+\mathds{1}_{(m,4)=(n,4)}\frac{1}{\phi(4/(n,4))}, &  \nu=1, \\ 
\mathds{1}_{(m,4)=(n,4)}\frac{1}{\phi(4/(n,4))}, & \nu=2, 
\end{cases}
\end{equation}
(noting that if $2\mid \delta_1'\delta_2'$ then the $w$-sum  is zero unless $2| mn, 2\nmid (m,n)$; and if $2\mid \delta\delta_1\delta_2$ then the $2$-part in the $b_d$-sum becomes $\frac{1}{\phi(4/(n,4))}\mathds{1}_{(m,4)=(n,4)}$),
and for $2\nmid p$ we have
\begin{equation}\label{r2}
\mathfrak f_2(p,\nu;m,n)=\mathfrak f_1(p,\nu;m,n).
\end{equation}
   
\subsection{Computation of the main term III  -- the sieve weights}   

We are now prepared to further evaluate $S_1^{-}$ and $S_2^{-}$. We will bound $S_1^{-}+S_2^{-}$ from below by \eqref{s2-} using properties of the lower bound sieve weights in Lemma \ref{lowerboundsieve}.

To begin with, we substitute \eqref{secondline} into \eqref{s1-}.  We still keep $h$ and $m$ fixed and consider the $u, c, d_1, d_2$-sum in \eqref{s1-} which now equals
$$\sum_{\substack{u\mid \mathcal Q\\ (u, h)=1}}\sum_{\substack{c\mid P(z)\\ (c,\mathcal Qh)=1}}\lambda_c^-\mu(c)\sum_{\substack{d_1,d_2\\ p\mid d_1d_2\Leftrightarrow p\mid uc}}\mu(d_1)\mu(d_2)\mathfrak c\Big(2d_1d_2\frac{n}{h}m\Big)\prod_{p^\nu\parallel d_1d_2}\mathfrak f_1(p,\nu;m,n).$$
We recast this as (noting $\mathfrak c(2)=1$)
$$\mathfrak c\Big(\frac{n}{h}m\Big)\sum_{\substack{u\mid \mathcal Q\\ (u, h)=1}}\sum_{\substack{c\mid P(z)\\ (c,\mathcal Qh)=1}}\lambda_c^- \mu(c)\sum_{\substack{d_1,d_2\\ p\mid d_1d_2\Leftrightarrow p\mid uc}}\mu(d_1)\mu(d_2)\prod_{\substack{p^\nu\parallel d_1d_2\\p\nmid 2}}\mathfrak w(p,\nu;m,n,h)\prod_{2^\nu\parallel d_1d_2}\mathfrak f_1(2,\nu;m,n)
$$
(of course each of the products over primes has only one factor, since $\nu$ is unique) where 
\begin{align}\label{w}
\mathfrak w(p, \nu;m,n,h)=\begin{cases}
\mathfrak c(p)\mathfrak f_1(p,\nu;m,n), & p\nmid 2\frac{n}{h}m,\\
\mathfrak f_1(p,\nu;m,n), &2\nmid p\mid \frac{n}{h}m.
\end{cases}
\end{align}
We observe that
  the sum over $d_1,d_2$ gives 
\begin{align}
\tilde{\mathfrak{f}}_1(2;m, n ; u)\prod_{\substack{p\mid uc\\ p\nmid 2}}(-2\mathfrak w(p,1;m,n,h)+\mathfrak w(p,2;m,n,h)) 
\end{align}
where
\begin{equation}\label{tildef}
\tilde{\mathfrak{f}}_1(2;m, n ; u) =  \begin{cases} 
 -2\mathfrak f_1(2,1;m,n)+\mathfrak f_1(2,2;m,n), & 2\mid u,\\
 \mathfrak f_1(2,0;m,n), &2 \nmid u.\end{cases}
\end{equation}

Using the formulas \eqref{f2} and \eqref{r2} instead of \eqref{f1} and \eqref{r1},  the same computation evaluates the corresponding $u, c, d_1, d_2$-sums of $S_2^-$  in \eqref{s2-def} as
\begin{align}\label{S2-}
&\sum_{\substack{u\mid \mathcal Q\\ (u, h)=1}}\sum_{\substack{c\mid P(z)\\(c, \mathcal Qh)=1}}\lambda_c^-\mu(c)\sum_{\substack{d_1,d_2\\ p\mid d_1d_2\Leftrightarrow p\mid uc}}\mu(d_1)\mu(d_2)\mathfrak c\Big(d_1d_2\frac{n}{h}m\Big)\prod_{p^\nu\parallel d_1d_2}\mathfrak f_2(p,\nu;m,n)\nonumber\\
&=\mathfrak c\Big(\frac{n}{h}m\Big)\sum_{\substack{u\mid \mathcal Q\\ (u,h)=1}}\sum_{\substack{c\mid P(z)\\ (c,\mathcal Qh)=1}}\lambda_c^-\mu(c) \tilde{\mathfrak{f}}_2(2;m, n ; u)
\prod_{ 2\nmid p\mid uc }(-2\mathfrak w(p,1;m,n,h)+\mathfrak w(p,2;m,n,h))\end{align}
where $\tilde{\mathfrak{f}}_2(2;m, n ; u)$ is defined in the same way as \eqref{tildef}, replacing the index $1$ with the index 2. 
Combining 
the above calculations, we have shown
\begin{displaymath}
\begin{split}
S_1^-+S_2^- = &\sum_{j=1}^2\sum_{\substack{h\mid n\\ P^+(h)\leq y\\\mu^2(h)=(h,\N)=1}} \chi_{-4}(h)\sum_{\substack{ hm\leq x\\ P^+(m)\leq y\\ \mu^2(mh)=1}}\mathfrak c\Big(\frac{n}{h}m\Big)\sum_{\substack{u\mid \mathcal Q\\ (u,h)=1}}\sum_{\substack{c\mid P(z)\\ (c,\mathcal Qh)=1}}\lambda_c^-\mu(c)\\
&\times  \tilde{\mathfrak{f}}_j(2;m, n ; u)
\prod_{ 2 \nmid p\mid cu}(-2\mathfrak w(p,1;m,n,h)+\mathfrak w(p,2;m,n,h)). 
\end{split}
\end{displaymath}

We are now finally in a position to apply Lemma \ref{lowerboundsieve} to bound from below the  sum over the sieve weights. We make the following important observations. From \eqref{w},  \eqref{r1}, \eqref{dpmn}, \eqref{cf} we have 
\begin{equation}\label{w12}
\begin{split}
0&\leq 2\mathfrak w(p,1;m,n,h)-\mathfrak w(p,2;m,n,h)\\
&=\left.\begin{cases}
\frac{1}{p}(2-\frac{2\chi_{p^*}(n)+1}{p-1})(1+\frac{\chi_{-4}(p)}{p(p-1)})^{-1}, & p\nmid nm\\
\frac{2}{p}(1-\frac{\chi_{-4}(p)}{p})^{-1}(1+\chi_{p^*}(n)), & p\mid \frac{n}{h}m,p\nmid (m,n)\\
 \frac{1}{p-1}(1-\frac{\chi_{-4}(p)}{p})^{-1}\mathds{1}_{\nu_p(n)=1},& p\mid (m,n)
\end{cases}\right\}
\leq \frac{5}{p},
\end{split}
\end{equation}
for $p\nmid 2h$. So we can choose $\mathcal{Q} = 30$ and use a sieve of dimension $\gamma = 5$. Secondly, since $h\mid n$, $(h, 2\N)=1$, we have $
\chi_{-4}(h)=1$, so that a lower bound for the sum of sieve weights with multiplicative coefficients suffices. This is an important device in the argument. Now applying Lemma \ref{lowerboundsieve} and \eqref{saz}, we obtain
\begin{align}\label{s2-}
 S(\mathcal A, n^{1/C})\geq \mathcal{S}(n)\Big(f_5(C\kappa_0)+O\Big(\frac{e^{-C\kappa_0}}{(\log \mathcal{D})^{1/3}}\Big)\Big) + O_A(\mathcal F(n,y)\Psi(x,y)(\log x)^{-A}),
\end{align}
where  $f_5$ has the same meaning as in Lemma \ref{lowerboundsieve}, $\mathcal{D}=n^{\kappa_0}$ for some  small $\kappa_0>0$ and $C$ a large constant such that $f_5(C\kappa_0)>0$ and 
\begin{equation}\label{Sn}
\mathcal{S}(n)=\sum_{\substack{h\mid n\\ P^+(h)\leq y\\\mu^2(h)=1\\(h,2\N)=1}}\sum_{\substack{ hm\leq x\\ P^+(m)\leq y\\ (m, \mathfrak N)=1\\\mu^2(m_oh)=1}}\mathfrak c\Big(\frac{n}{h}m\Big) \mathfrak F_{\mathcal{Q}}(m,n,h)
\prod_{\substack{p\leq z\\ p\nmid \mathcal Qh}}(1-2\mathfrak w(p,1;m,n,h)+\mathfrak w(p,2;m,n,h)), 
\end{equation}
where $\mathfrak F_{\mathcal{Q}}(m,n,h)=\prod_{\substack{p\mid \mathcal Q\\ p\nmid h}}\mathfrak F_p(m,n,h)$ with 
\begin{align}\label{F2}
\mathfrak F_2(m,n,h)&=\sum_{j=1}^2(\mathfrak f_j(2,0;m,n)-2\mathfrak f_j(2,1;m,n)+\mathfrak f_j(2,2;m,n))
\end{align} (the right hand side is independent of $h$) and
\begin{align}\label{Fp}
\mathfrak F_p(m,n,h)=\begin{cases}
(1-2\mathfrak w(p,1;m,n,h)+\mathfrak w(p,2;m,n,h)), & p \nmid 2h,\\
1, & p\mid h.
\end{cases} 
\end{align}
 Recalling \eqref{w12}, we define multiplicative functions $\mathfrak g_i$  supported on square-free integers with $\mathfrak g_i(p)=1$ except for $ p\mid P_{\mathcal Q}(z)$ where
\begin{equation}\label{defgi}
\begin{split}
&\mathfrak g_0(p)=1-\Big(\frac{2}{p}-\frac{2\chi_{p^*}(n)+1}{p(p-1)}\Big)\Big(1+\frac{\chi_{-4}(p)}{p(p-1)}\Big)^{-1},\\
&\mathfrak g_1(p)=1-\frac{2}{p}\Big(1+\chi_{p^*}(n)\Big)\Big(1-\frac{\chi_{-4}(p)}{p}\Big)^{-1},\\
&\mathfrak g_2(p)=1-\frac{1}{p-1}\Big(1-\frac{\chi_{-4}(p)}{p}\Big)^{-1}\mathds{1}_{\nu_p(n)=1},
\end{split}
\end{equation}
so that 
\begin{displaymath}
\begin{split}
\mathcal{S}(n) = &\sum_{\substack{h\mid n \\ P^+(h)\leq y\\ \mu^2(h)=(h, 2\mathfrak N)=1}}\sum_{\substack{hm\leq x\\ P^+(m)\leq y\\ (m, \mathfrak N)=1\\\mu^2(m_oh)=1}}\mathfrak c\Big(\frac{n}{h}m\Big)\mathfrak{F}_{\mathcal Q}(m, n,h)\prod_{\substack{ p\nmid nm\\p\nmid h}}\mathfrak g_0(p)\prod_{\substack{ p\mid \frac{n}{h}m\\ p\nmid h(m,n)}}\mathfrak g_1(p)\prod_{\substack{p\mid (m,n)\\p\nmid h}}\mathfrak g_2(p). 
\end{split}
\end{displaymath}
We  rename $mh$ by $m$, $h$ by $hg$ where $h\mid (m,n), h\nmid \mathcal Q$, and $2\nmid g\mid (m,n, \mathcal Q)$. 
 In this way we can re-write $\mathcal{S}(n)$ as 
\begin{displaymath}
\begin{split}
&\prod_{p\mid P_{\mathcal Q}(z)}\mathfrak g_0(p) \sum_{\substack{m\leq x\\ P^+(m)\leq y\\(m, \mathfrak N)=1\\ \mu^2(m_o)=1}}\sum_{\substack{g\mid (n, m, \mathcal Q)\\(g, 2\mathfrak N)=1}}\mathfrak{F}_{\mathcal Q}\Big(\frac{m}{g}, n,g\Big)\prod_{\substack{p\mid \mathcal Q\\ p\mid \frac{n}{g}\frac{m}{g}}}\mathfrak c(p)\sum_{\substack{h\mid (m,n)\\ (h, \mathcal Q\mathfrak N)=1}}\frac{1}{\mathfrak g_0(h)}\prod_{\substack{ p\mid(m,n)\\ p\nmid h}}\frac{\mathfrak g_2(p)}{\mathfrak g_0(p)}\prod_{\substack{p \mid mn\\p\nmid (m,n)}}\frac{\mathfrak g_1(p)}{\mathfrak g_0(p)}\prod_{\substack{p\nmid \mathcal Q\\ p\mid \frac{n}{h}\frac{m}{h}}}\mathfrak c(p)\\
&=\prod_{p\mid P_{\mathcal Q}(z)}\mathfrak g_0(p) \sum_{\substack{r\mid 2^\infty\mathcal Q\\(r, \mathfrak N)=1}}\Big(\sum_{\substack{g\mid (n,r, \mathcal Q) \\ \mu^2(g)=1\\(g, 2\mathfrak N)=1}}\mathfrak{F}_\mathcal {Q}\Big(\frac{r}{g}, n,g\Big)\prod_{\substack{p\mid \mathcal Q\\p\mid \frac{n}{g}\frac{r}{g}}}\mathfrak c(p)\Big)\\
&\quad\quad\quad\quad\quad\times \sum_{\substack{m\leq x/r\\ P^+(m)\leq y\\ (m,\mathcal{Q}\mathfrak N)=1\\\mu^2(m)=1}}\Big(\sum_{\substack{h\mid (m,n)}}\frac{1}{\mathfrak g_0(h)}\prod_{\substack{ p\mid(m,n)\\ p\nmid h}}\frac{\mathfrak g_2(p)}{\mathfrak g_0(p)}\prod_{\substack{ p\mid mn\\p\nmid (m,n)}}\frac{\mathfrak  g_1(p)}{\mathfrak g_0(p)}\prod_{\substack{p\nmid \mathcal Q\\ p\mid \frac{n}{h}\frac{m}{h}}}\mathfrak c(p)\Big). 
\end{split}
\end{displaymath}
We first see that the sum over $g$ can be written as (note that $\mathfrak F_p(\bullet,\bullet,p)=1$)
\begin{align}
\prod_{\substack{p\mid \mathcal Q }}\Big(\mathfrak F_{p}(r,n,1)\prod_{p\mid rn}\mathfrak c(p)+\mathds{1}_{\substack{p\mid (n,r)\\p\nmid 2\mathfrak N}}\cdot\prod_{\substack{p\mid \frac{n}{p}\frac{r}{p}}}\mathfrak c(p)\Big)=:\mathfrak F_\mathcal Q(r,n)=\prod_{p\mid \mathcal Q}\mathfrak F_p(r,n).
\end{align}
With $\mathcal Q=30$, $\mu^2(r_o)=1$, this becomes 
\begin{align}
\mathfrak F_2(r,n,1)\mathfrak F_3(r, n,1)\prod_{3\mid nr}\mathfrak c(3)\Big(\mathfrak F_5(r,n,1)\prod_{5\mid nr}\mathfrak c(5)+\mathds{1}_{5\mid (n,r)}\cdot\prod_{p\mid \frac{n}{5} }\mathfrak c(5)\Big).
\end{align}
Using \eqref{F2}, \eqref{f1}, \eqref{f2} and \eqref{Fp}, \eqref{w12}, we compute explicitly 
\begin{align}\label{FF23}
\mathfrak F_2(r,n)=\begin{cases}
\frac{1}{2}\mathds{1}_{2\nmid r}, &  n\equiv 1\pmod 2,\\
2\mathds{1}_{8\mid r},&n \equiv 2\pmod 8,\\
2\mathds{1}_{4\parallel r},&n \equiv 6\pmod 8,\\
\mathds{1}_{2\parallel r},& n\equiv 0 \pmod 4.
\end{cases},\quad
\mathfrak F_3(r,n)=\begin{cases}
\frac{4}{5}\mathds{1}_{3\nmid r},&n \equiv 1\pmod 3,\\
\frac{8}{5}\mathds{1}_{3\mid r},&n \equiv 2\pmod 3,\\
\frac{4}{5}\mathds{1}_{3\nmid r},& n\equiv 0 \pmod 3,
\end{cases}
\end{align}
(note that if $3\mid n$ then $(r, 3) = 1$ since $(r, \mathfrak{N}) = 1$) and (with $\mathfrak c(5)=\frac{16}{21}$)
\begin{align}\label{FF5}
\mathfrak F_5(r,n)=\begin{cases}
\mathfrak c(5)\mathds{1}_{5\nmid r},&n \equiv 1,4\pmod 5,\\
\mathfrak c(5)\mathds{1}_{5\mid r}+\frac{4}{7}\mathds{1}_{5\nmid r},&n \equiv 2,3\pmod 5,\\
2\mathfrak c(5)\mathds{1}_{\substack{5\mid n\\ 5\mid r}}+\frac{1}{2}\mathfrak c(5)\mathds{1}_{5\nmid r},& n\equiv 0 \pmod 5.
\end{cases}
\end{align}
In particular, we can always find some residue class $r_0\pmod{ 120}$ such that $\mathfrak F_{\mathcal Q}(r,n)>0$ for $r\equiv r_0\pmod {120}$.
 
Noting that $(m, \mathcal Q)= \mu^2(m)=1$ and $\mathfrak g_i(p)=1$ if $p\mid \mathcal Q$, we next evaluate the inner sum over $h$ as
\begin{equation}\label{defGG}
\begin{split}
&\sum_{\substack{h\mid (m,n)}}\frac{1}{\mathfrak g_0(h)}\prod_{\substack{ p\mid(m,n)\\ p\nmid h}}\frac{\mathfrak g_2(p)}{\mathfrak g_0(p)}\prod_{\substack{ p\mid mn\\p\nmid (m,n)}}\frac{\mathfrak g_1(p)}{\mathfrak g_0(p)}\prod_{\substack{p\nmid \mathcal Q\\ p\mid \frac{n}{h}\frac{m}{h}}}\mathfrak c(p)\\
& =\sum_{h\mid (m,n)}\frac{1}{\mathfrak g_0(h)}\prod_{p\mid h}\frac{\mathfrak g_0(p)}{\mathfrak g_2(p)}\prod_{\substack{p\mid h\\ \nu_p(n)=1}}\frac{1}{\mathfrak c(p)}\prod_{p\mid (m,n)}\frac{\mathfrak g_2(p)}{\mathfrak g_0(p)}\prod_{\substack{p\mid m\\p\nmid n}}\frac{\mathfrak g_1(p)}{\mathfrak g_0(p)}\mathfrak c(p)\prod_{\substack{p\mid n\\ p\nmid m}}\frac{\mathfrak g_1(p)}{\mathfrak g_0(p)}\prod_{\substack{p\mid n \\ p\nmid m \mathcal Q}}\mathfrak c(p)\prod_{\substack{ p\mid (m,n)}}\mathfrak c(p)\\
& =
\sum_{\substack{h\mid (m,n)}}\Big(\frac{1}{\mathfrak g_0(h)}\prod_{p\mid h}\frac{\mathfrak g_0(p)}{\mathfrak g_2(p)}\prod_{\substack{p\mid h\\ \nu_p(n)=1}}\frac{1}{\mathfrak c(p)}\Big)\prod_{\substack{ p\mid (m,n)}}\frac{\mathfrak g_2(p)}{\mathfrak g_1(p)}\prod_{\substack{p\mid m \\ p\nmid n}}\frac{\mathfrak g_1(p)}{\mathfrak g_0(p)}\mathfrak c(p)\prod_{\substack{p\nmid \mathcal Q\\ p\mid n}}\frac{\mathfrak g_1(p)}{\mathfrak g_0(p)}\mathfrak c(p)\\
& =\prod_{\substack{p\mid (m,n)\\ \nu_p(n)=1}}\Big(1+ \frac{1}{\mathfrak c(p) \mathfrak g_2(p)}\Big)\prod_{\substack{p\mid (m,n)\\ \nu_p(n)>1}}\Big(1+ \frac{1}{\mathfrak g_2(p)}\Big)\prod_{\substack{ p\mid (m,n)}}\frac{\mathfrak g_2(p)}{\mathfrak g_1(p)}\prod_{\substack{ p\mid m \\ p\nmid n}}\frac{\mathfrak g_1(p)}{\mathfrak g_0(p)}\mathfrak c(p)\prod_{\substack{p\nmid \mathcal Q\\ p\mid n}}\frac{\mathfrak g_1(p)}{\mathfrak g_0(p)}\mathfrak c(p)\\
& =:\mathfrak G(m;n) \prod_{\substack{p\nmid \mathcal Q\\ p\mid n}}\frac{\mathfrak g_1(p)}{\mathfrak g_0(p)}\mathfrak c(p),
\end{split}
\end{equation}
where $\mathfrak G(m;n)$ is a multiplicative function of $m$ with $\mathfrak G(p^\nu;n)=0$ for $\nu\geq 2$ and $\mathfrak G(p;n)=0$ for $p\mid \mathcal Q\mathfrak N$. We summarize
\begin{equation}\label{summarize}
\mathcal{S}(n) = \prod_{p\mid P_{\mathcal Q}(z)}\mathfrak g_0(p) \sum_{\substack{r\mid 2^\infty\mathcal Q\\ (r, \mathfrak N)=1}}\mathfrak{F}_{\mathcal Q}(r, n)\prod_{\substack{p\nmid \mathcal Q\\ p\mid n}}\frac{\mathfrak g_1(p)}{\mathfrak g_0(p)}\mathfrak c(p)\sum_{\substack{m\leq x/r\\ P^+(m)\leq y}}\mathfrak G(m;n).
\end{equation}

\subsection{Computation of the main term IV -- summing over smooth square-free numbers}\label{part4}   

Our last goal is a lower bound for $\mathcal{S}(n)$ based on the formula \eqref{summarize}. 

In order to carry out the sum over $m$ in \eqref{summarize}, it is useful to 
 define $\mathfrak h(\bullet; n)=\mathfrak G(\bullet; n)*\mu$. Recall that the definition of $\mathfrak{G}(m;n)$ is implicit in \eqref{defGG}, so that  
  $\mathfrak h(p^\nu;n) = 0$ if $\nu> 2$,  and for $\nu \leq 2$ we have  
  \begin{align}
\mathfrak h(p;n)=\begin{cases}
(1+\frac{1}{\mathfrak c(p)\mathfrak g_2(p)})\frac{\mathfrak g_2(p)}{\mathfrak g_1(p)}-1, & p\mid n,p\nmid \mathcal Q\mathfrak N, \nu_p(n)=1,\\
(1+\frac{1}{\mathfrak g_2(p)}) \frac{\mathfrak g_2(p)}{\mathfrak g_1(p)}-1, & p\mid n,p\nmid \mathcal Q\mathfrak N, \nu_p(n)>1,\\
-1, &  p\mid \mathcal Q\mathfrak N,\\
\frac{\mathfrak g_1(p)}{\mathfrak g_0(p)}\mathfrak c(p)-1, & p\nmid n \mathcal Q
\end{cases}
\end{align}
and
\begin{align}
\mathfrak h(p^2;n)= - \mathfrak{h}(p, n). 
\end{align}
In particular, from \eqref{defgi} and \eqref{cf} we have
\begin{align}\label{hpn}
\mathfrak h(p;n)= \begin{cases}
1+\frac{4}{p-3}, &p\mid n, p\nmid \mathcal Q\mathfrak N,\\
-1, & p\mid \mathcal Q\mathfrak N,\\
O(\frac{1}{p}), &p\nmid n\mathcal Q,
\end{cases}\quad
\mathfrak h(p^2;n) = O(1),
\end{align}
where we can take $5$ as implied constants. Using \eqref{psix/d} and $\alpha\geq 1-\frac{1}{D}+o(1)$ with $D$ sufficiently large, the contribution from $r\geq K=n^{1/100}$ in \eqref{summarize} can be bounded by 
\begin{align}\label{largetk}
\Psi(x, y) \sum_{t_1t_2^2 \geq K} \frac{(t_1, n)}{t_1 (t_1t^2_2)^{\alpha}} \ll \Psi(x,y) \frac{1}{K^{2\alpha/5}}\sum_{t_1, t_2} \frac{(t_1, n)}{t_1 (t_1t^2_2)^{3\alpha/5}}\ll K^{-1/3} n^{\varepsilon} \Psi(x, y). 
\end{align}
To evaluate the contribution from $r\leq K$, we consider two cases depending on the size of $y$.
 \subsubsection{Case I: $y\geq \exp((\log_2x)^2)$}
We write 
 \begin{align}\label{G-sum}
 \sum_{\substack{m\leq x/r\\ P^+(m)\leq y}}\mathfrak G(m;n)=\sum_{\substack{t\leq x/r\\ P^+(t)\leq y}}\mathfrak h(t;n) \Psi\left(\frac{x}{tr},y\right).
 \end{align}
 We can truncate the $t$-sum using \eqref{psix/d} as in \eqref{largetk}
so that the contribution from $t\geq K$ \eqref{G-sum} can be bounded by $K^{-1/3}n^\varepsilon\Psi(x,y)$.
 For $t,r \leq K$ we evaluate  the sum by \cite[Theorem 2 with $q=1$]{FT91} 
 getting 
 \begin{equation}\label{tksum}
 \begin{split}
 &\sum_{\substack{t\leq K\\ P^+(t)\leq y }} \mathfrak{h}(t;n)  \frac{x}{tr}\int_{0^-}^\infty\rho\Big( \frac{\log x}{\log y}-\frac{\log (tr)}{\log y}-v\Big)dR(y^v)\Big(1+O\big(\exp(-(\log y)^{0.55})\big)\Big).
 \end{split}
 \end{equation}
 where
 $R(x) = [x]/x$. Strictly speaking, the application of 
 \cite[Theorem 2]{FT91} requires $y \leq x/K^2 \leq x/(tr)$, which we can assume without loss of generality since $\Psi(x, y) \asymp \Psi(x, x)$ for $y > x^{1/2}$, say. 
 
 To estimate the error term, we 
 use \eqref{psix/d} again, and obtain the bound
 \begin{equation}\label{smalltk}
 \begin{split}
 \sum_{t\leq K} |\mathfrak{h}(t; n)| \Psi\Big( \frac{x}{tr}, y\Big)\exp(-(\log y)^{0.55})& \ll  \prod_{ p\mid n }\Big(1+\frac{O(1)}{p^{\alpha}}\Big)  \Psi(x, y) \exp(-(\log y)^{0.55})\\
 & \ll \Psi(x, y)  \frac{(\log n)^{O(1)}}{ \exp((\log y)^{0.55})} \ll \frac{\Psi(x, y)}{(\log x)^{100}}
 \end{split}
 \end{equation}
 for $y \geq \mathcal \exp((\log_2x)^2)$. 
 
 We now focus on the main term in \eqref{tksum} which we substitute back as the $m$-sum in \eqref{summarize}. By partial summation, the main term equals 
 \begin{align}\label{smoothmain1}
 \prod_{p\mid P_{\mathcal Q}(z)}\mathfrak g_0(p) \cdot x \int_{0^-}^\infty \rho(u-v) dW(y^v)
 \end{align}
 where
 \begin{align}
 W(U)&=\frac{1}{U}\prod_{\substack{p\nmid \mathcal Q\\  p\mid n}}\frac{\mathfrak g_1(p)}{\mathfrak g_0(p)}\mathfrak c(p)\sum_{\substack{r\leq K \\r\mid 2^\infty\mathcal Q\\ (r, \mathfrak N)=1}}\mathfrak F_{\mathcal Q}(r,n)
 \sum_{\substack{t\leq K \\ P^+(t)\leq y}}\mathfrak{h}(t;n) \sum_{\substack{ m\leq U /tr}}1. 
 \end{align}
 It is clear that $W(U) = 0$ for $U < 1$ and $W(U) = c + O(1/U)$ for some constant $c$ (depending on $n$ and  $z$). 
 In preparation for an application of \cite[Lemme 4.1]{FT90} we define
 $$ M(s)=\int e^{-sv}dW(e^v).$$
 We compute $M(s)$ as $$M(s)= 
 \mathfrak S_K(n,s) \frac{s}{s+1}\zeta(s+1)$$
 where
 \begin{align}
 \mathfrak S_K(n,s)= \prod_{\substack{p\nmid \mathcal Q\\ p\mid n}}\frac{\mathfrak g_1(p)}{\mathfrak g_0(p)} \mathfrak c(p) \sum_{\substack{r\leq K\\ r\mid 2^\infty\mathcal Q\\ (r, \mathfrak N)=1}}\frac{\mathfrak F_{\mathcal Q}(r,n)}{r^{s+1}}\sum_{\substack{t\leq K \\ P^+(t)\leq y}}\frac{\mathfrak h(t;n)}{t^{s+1}}.
 \end{align}
 Recalling \eqref{defgi} and \eqref{cf} we have 
 \begin{equation}\label{arith}
 \frac{60}{67} \leq \frac{\mathfrak g_1(p)}{\mathfrak g_0(p)}\mathfrak c(p)=1- \frac{\chi_{-4}(p)}{p}+O(p^{-2}), \quad p\mid n, \,\, p\nmid \mathcal Q = 30. 
 \end{equation}
 We can extend $r, t$ to infinity getting
 \begin{align}
 &\mathfrak S(n,s)= \prod_{\substack{p\nmid \mathcal Q\\ p\mid n}}\frac{\mathfrak g_1(p)}{\mathfrak g_0(p)} \mathfrak c(p)\sum_{k=0}^\infty \frac{\mathfrak F_2(2^k,n)}{2^{k(s+1)}}\prod_{\substack{p\mid \mathcal Q\\ (p, 2\mathfrak N)=1}}\Big(\mathfrak F_p(1,n)+ \frac{\mathfrak F_p(p,n)}{p^{s+1}}\Big)\prod_{\substack{ p\leq y}}\Big(1+\frac{ \mathfrak h(p;n)}{p^{s+1}}+\frac{\mathfrak h(p^2;n)}{p^{s+2}}\Big)
 \end{align}
 where uniformly in $|s| \leq 1/2$ we have  (once again by Rankin's trick)
 \begin{equation}\label{errorK}
 \mathfrak S(n,s) - \mathfrak S_K(n,s) \ll \frac{1}{K^{1/2}}\prod_{p\mid n}\Big(1+ \frac{1}{p^{\Re s+1/2}}\Big)\Big) \ll \frac{1}{K^{1/4}} = n^{-1/400},
 \end{equation}
 and by Cauchy's formula this bound remains true for all fixed derivatives with respect to $s$. Note that by \eqref{hpn} we have 
 \begin{align}\label{hpbound}
 \frac{6}{7}\leq 1 + \frac{\mathfrak{h}(p; n)}{p} + \frac{\mathfrak{h}(p^2; n)}{p^2} = 1 + \frac{\chi_{-4}(p)}{p} + O(p^{-2}),\quad p\mid n.
 \end{align}
 Combing \eqref{FF23}, \eqref{FF5}, \eqref{hpn}, \eqref{arith} and \eqref{hpbound},  we obtain for all $j\geq 0$ that 
 \begin{equation}\label{boundsS}
 \frac{d^j}{ds^j} \mathfrak S(n,s)|_{s = 0} \asymp \prod_{\substack{p\mid n\\ p\leq \min(y,z)}}\Big(1+O\Big(\frac{(\log p)^j}{p^{2}}\Big)\Big)\prod_{\substack{p\mid n\\ p\geq\min( y,z)}}\Big(1+O\Big(\frac{(\log p)^j}{p}\Big)\Big) \ll 1
 \end{equation}
 for $y \geq \mathcal \exp((\log_2x)^2)$ and $\log z\asymp \log n$,  and we also see 
 $\mathfrak S(n, 0) \gg 1.$
 
 We are now prepared to apply \cite[Lemme 4.1]{FT90} to evaluate \eqref{smoothmain1}, which we complement with \cite[Lemme 4.1]{FT91} for bounds of $\rho^{(j)}(u)$. Using again 
 $
 y\geq \mathcal L(n)
 $ and $\log z\asymp \log n$ and recalling the error terms \eqref{smalltk}, \eqref{largetk} and \eqref{errorK}, we conclude 
 \begin{equation}\label{lower}
 \begin{split}
 \mathcal{S}(n)&\asymp \Big(\prod_{p\mid P_{\mathcal Q}(z)}\mathfrak g_0(p) \Big)
 \mathfrak S(n,0) x\rho(u)\Big(1+O\Big(\frac{\log \log n}{\log y}\Big)\Big)  + O\Big(\frac{\Psi(x, y)}{(\log x)^{100}} \Big)\\
 &\gg \frac{
 	x \rho(u)}{(\log z)^2} +   O\Big(\frac{\Psi(x, y)}{(\log x)^{100}} \Big) \asymp \frac{ \Psi(x, y)}{(\log z)^{2}}.
 \end{split}
 \end{equation}

 We plug this lower bound back into \eqref{s2-}, noting that $\mathcal F(n,y)\asymp 1$ when $y\geq \exp((\log_2x)^2)$, and choose $C=C(\kappa_0)$ a sufficiently large constant so that $f_5(C\kappa_0) > 0$ getting
 $$S(\mathcal{A}, n^{1/C}) \gg 
 (\log z)^{-2} \Psi(x, y)$$
 provided that $y \geq \mathcal \exp((\log_2x)^2)$ and $\mathcal{D}=n^{\kappa_0}$ for some sufficiently small $\kappa_0>0$. 
 This concludes the proof of Theorem \ref{thm2} when $ y\geq \exp((\log_2x)^2)$. \\

 \textbf{Remark.} We close this computation by a short comment on the ``singular series'' $\mathfrak{S}(n, 0)$ encountered in the final lower bound, which might look a bit unexpected. This is due to our restrictions on $m$, in particular the fact that we impose a square-free condition on $m$.  As  square-free numbers are not equidistributed among all residue classes,  the singular series $\mathfrak S(n,0)$ will genuinely depend on the factorization of $n$ in our set-up. For example, we can see from \eqref{FF23} that there are fewer solutions to $n=m+x^2+y^2$ with $(xy,3)=1$, $\mu^2(m_o)=1$ when $n\equiv 2\pmod 3$ than in the other cases, whilst the number of solutions to $n=m+x^2+y^2$ with $(xy,3)=1$ and no extra conditions on $m$ would be approximately equal for all residue classes of $n\pmod 3$.

\subsubsection{Case II: $(\log x)^D\leq y\leq \exp((\log_2x)^2)$}
When $y$ is small, the factor $\prod_{p\mid n}(1+\frac{O(1)}{p^\alpha})$ in \eqref{smalltk} can be much bigger than a large power of $\log x$ and so we cannot bound $\mathfrak h(t;n)$ by absolute values. Instead, we apply Perron's formula to evaluate the $m$-sum in \eqref{summarize} directly as follows. Write
\begin{equation}\label{integral}
\mathcal S_0(n,x/r):=\sum_{\substack{m\leq x/r\\ P^+(m)\leq y}}\mathfrak G(m;n)=\frac{1}{2\pi i}\int_{(\alpha)}G(s;n)\Big(\frac{x}{r}\Big)^s\frac{ds}{s}
\end{equation}
where $G(s;n)=\sum_{m=1}^\infty\mathfrak G(m;n)n^{-s}$ and as before $\alpha$ is defined in \eqref{alphadef}. (Recall that for $D$ large enough we have $9/10 \leq \alpha \leq 1$.) 
Denote $$\zeta_m(s,y)=\prod_{\substack{p\leq y\\ p\nmid m}}(1-p^{-s})^{-1}, \quad \zeta(s,y)=\zeta_1(s,y).$$
Using \eqref{hpbound} we can write
 \begin{align}\label{Gsdef}
G(s;n)=\sum_{m=1}^\infty\frac{\mathfrak h(m;n)}{m^s}\zeta(s,y)=\prod_{\substack{p\mid n_y}}\Big(1-\frac{1}{p^s}\Big)^{-\chi_{-4}(p)}\zeta(s,y)\tilde{G}(s;n)
\end{align}
 for some function $\tilde{G}(s;n)$ which is holomorphic and uniformly bounded in $\Re(s)\geq 2/3$, say.
 We use the saddle point method to evaluate \eqref{integral}. The main term comes from a small neighbourhood of the point $s=\alpha$. 
  Precisely, we truncate the integral in \eqref{integral} at height $1/\log y$ to write
  \begin{align}\label{S0n}
  \mathcal S_0(n,x/r)=\frac{1}{2\pi i}\int_{\alpha-i/\log y}^{\alpha+i/\log y}G(s;n)\Big(\frac{x}{r}\Big)^s\frac{ds}{s}+E(x,y),
  \end{align}
   say. 
   To handle the error term $E(x,y)$ we need some estimate of $G(s;n)$ when $s=\alpha+i\tau$ for $|\tau| \geq 1/\log y$. 
  This can be done following the proof of the bound for $\zeta_m(s,y)/\zeta_m(\alpha;y)$ in \cite[(4.42)]{lBT}. In fact, the bound for $\zeta_m(s,y)/\zeta_m(\alpha;y)$ follows from the case of $m=1$ together with  bounds for the contribution from primes satisfying $p\mid m$. In our case, we can bound the contribution from $p\mid n_y$ in the same way (since $\omega(n_y)\ll \sqrt{y}$ for $y\geq \mathcal L(x)$) and thus we can obtain for $s=\alpha+i\tau$ 
  that 
\begin{align}\label{Gbound}
\frac{G(s;n)}{G(\alpha;n)}\ll \begin{cases}
\exp\Big(-\frac{c_0y}{\log y}\log \Big(1+ \frac{\tau^2\sigma_2}{y/\log y}\Big)\Big), & |\tau| \leq 1/\log y,\\
\exp \Big(\frac{-c_0\tau^2 u}{(1-\alpha)^2+\tau^2}\Big), & 1/\log y< |\tau|\leq  Y_\varepsilon
\end{cases}
\end{align}
where $Y_\varepsilon=\exp\big((\log y)^{3/2-\varepsilon}\big)$ 
and 
\begin{align}\label{sigkdef}
\sigma_k=\frac{d^k}{ds^k}\log \zeta(s,y)\Big\vert_{s=\alpha}.
\end{align} 
Note that we have (see e.g. \cite[eq.\ (3.9)]{lBT})
\begin{align}\label{sigkbound}
\quad \sigma_k\asymp u(\log y)^k\asymp \log x(\log y)^{k-1}
\end{align}
for $k \geq 1$. Similarly to the proof of \cite[eq.\ (4.49)]{lBT}, we can show that for $1\leq z\leq Y_\varepsilon$ we have
\begin{align}\label{s0short}
S_0(n, x+x/z)-S_0(n,x)\ll x^{\alpha}G(\alpha;n)(1/z+e^{-cu})
\end{align}
for some small $c>0$. Set 
\begin{align}\label{Rdef}
R:=e^{-c'u/(\log 2u)^2}+1/Y_\varepsilon
\end{align} for some $c'>0$. With \eqref{Gbound} and  \eqref{s0short}, we can follow \cite[Lemma 10]{HT} to bound the contribution from  $1/\log y\leq |\tau|\leq T$ and $|\tau|\geq T$ in \eqref{integral} respectively with $T:=1/R^2,$ so that \eqref{S0n} holds with
\begin{align}\label{Exy}
E(x,y)\ll (x/r)^{\alpha}G(\alpha;n)R,
\end{align}
which can be bounded using \eqref{Rdef} and \eqref{sigkbound} by
\begin{align}\label{Ebound}
E(x,y)\ll \frac{x^\alpha G(\alpha;n)}{r^\alpha\alpha \sqrt{\sigma_2}u}.
\end{align}

It remains to evaluate the integral in \eqref{S0n}.
Let $T_0=(u^{1/3}\log y)^{-1}$.
Using \eqref{Gbound}, the contribution from $|\tau|\leq T_0$ can be bounded (see also \cite[p.\ 281]{HT})
\begin{align}\label{>T0}
\frac{x^\alpha}{r^\alpha}\frac{G(\alpha;n)}{\alpha}\int_{T_0}^{1/\log y}\Big(1+\frac{\tau^2\sigma_2}{y/\log y}\Big)^{-c_0y/\log y}d\tau&\ll \frac{x^\alpha}{r^\alpha}\frac{G(\alpha;n)}{\alpha \sqrt{\sigma_2}}\int_{T_0 \sqrt{\sigma_2}}e^{-c_0\tau^2/2}d\tau \\
&\ll \frac{x^\alpha}{r^\alpha}\frac{G(\alpha;n)e^{-c_1u^{1/3}}}{\alpha \sqrt{\sigma_2}}.
\end{align}
The contribution from $|\tau|\leq T_0$ in the integral \eqref{S0n} gives the main term. We can show that 
\begin{align}\label{<T0}
\frac{1}{2\pi i}\int_{\alpha-T_0}^{\alpha+T_0 }G(\alpha;n)\frac{(x/r)^s}{s}ds=\frac{x^\alpha G(\alpha; n)}{r^\alpha\alpha \sqrt{ 2\pi \sigma_2}}\Big(1+O\Big(\frac{1}{u^{1/3}}\Big)\Big).
\end{align}
To prove \eqref{<T0}, we compare it with the estimation for $\Psi(x,y)$. Using \cite[last equation on p. 280 in the proof of Lemma 11]{HT} which states that uniformly for $x\geq y\geq 2$ we have 
\begin{align}
\frac{1}{2\pi i}\int_{\alpha-iT_0}^{\alpha+iT_0}\zeta(s,y)\frac{x^s}{s}ds=\frac{x^\alpha \zeta(\alpha, y)}{\alpha \sqrt{ 2\pi \sigma_2}}\Big(1+O\Big(\frac{1}{u}\Big)\Big),
\end{align} we can write  
\begin{align}\label{S02}
&\frac{1}{2\pi i}\int_{\alpha-iT_0}^{\alpha+iT_0}G(\alpha;n)\frac{(x/r)^s}{s}ds\\
=&\frac{x^\alpha G(\alpha;n)}{r^\alpha\alpha \sqrt{ 2\pi \sigma_2}}\Big(1+O\Big(\frac{1}{u}\Big)\Big)+\frac{1}{2\pi i}\int_{\alpha-iT_0}^{\alpha+iT_0}\Big(\frac{G(s;n)}{\zeta(s,y)}-\frac{G(\alpha;n)}{\zeta(\alpha,y)}\Big)\zeta(s,y)\Big(\frac{x}{r}\Big)^s\frac{ds}{s}.
\end{align}
Thus it remains to show that second term in \eqref{S02} is of smaller order.
Set $$\gamma(s):=\log \frac{G(s;n)}{\zeta(s,y)}.$$ Then we have uniformly for $1/2<\sigma<1$ and $\Re(s)=\sigma$ that
\begin{align}
&|\gamma'(s)|\leq \gamma'(\sigma)\ll \sum_{p\mid n_y}\frac{\log p}{p^{\sigma}-1}\ll (\log n)^{1-\sigma},\label{gamma'}\\
& |\gamma''(s)|\leq \gamma'(\sigma)\ll \sum_{p\mid n_y}\frac{(\log p)^2p^\sigma}{(p^\sigma-1)^2}\ll (\log n)^{1-\sigma}\log_2n\label{gamma''}.
\end{align}
We have the Taylor expansion around $\tau=0$ for $|\tau|\leq T_0\ll (\sigma_3)^{-1/3}$ (see e.g.\ \cite[eq.\ (4.54)]{lBT})
\begin{align}
\frac{\zeta(s,y)x^s}{s}=\frac{\zeta(\alpha, y)x^\alpha}{s}e^{-\tau^2\sigma_2/2}(1+O(\tau^3\sigma_3+\tau/\alpha))
\end{align}
where $\sigma_i$ is as in \eqref{sigkdef},  
and for $|\tau|\leq T_0$ we have 
\begin{align}
\frac{G(s;n)}{\zeta(s,y)}-\frac{G(\alpha;n)}{\zeta(\alpha,y)}&=\frac{G(\alpha;n)}{\zeta(\alpha;n)}\exp\Big(i\tau \gamma'(\alpha)+O(\tau^2\gamma''(\alpha))-1\Big)\\
&=\frac{G(\alpha;n)}{\zeta(\alpha;n)}\Big(i\tau\gamma'(\alpha)+O\big(\tau^2(|\gamma''(\alpha)|+\gamma'(\alpha)^2)\big)\Big).
\end{align}
Using this together with the formulas
\begin{align}
\int_{-T_0}^{T_0}\tau e^{-\sigma_2\tau^2/2}d\tau=0,\quad \int_{-T_0}^{T_0}|\tau|^k e^{-\sigma_2\tau^2/2}d\tau \ll_k \sigma_2^{-(k+1)/2} \quad (k\geq 1),
\end{align}
we see that
\begin{align}\label{T0err}
&\frac{1}{2\pi i}\int_{\alpha-iT_0}^{\alpha+iT_0}\Big(\frac{G(s;n)}{\zeta(s,y)}-\frac{G(\alpha;n)}{\zeta(\alpha,y)}\Big)\zeta(s,y)\Big(\frac{x}{r}\Big)^s\frac{ds}{s} \ll \frac{x^\alpha G(\alpha;n)}{r^\alpha\alpha \sqrt{\sigma_2}}\left(|\gamma'(\alpha)|\Big( \frac{\sigma_3}{\sigma_2^2}+\frac{1}{\sigma_2\alpha}\Big)+\frac{|\gamma''(\alpha)|+\gamma'(\alpha)^2}{\sigma_2}\right).
\end{align} 
Using \eqref{sigkbound}, \eqref{gamma'} and \eqref{gamma''} this can be bounded by (for $\alpha\geq 9/10$)
\begin{align}
\frac{x^\alpha G(\alpha;n)}{r^\alpha\alpha \sqrt{\sigma_2}}\Big( \frac{|\gamma'(\alpha)|}{u\log y}+\frac{|\gamma''(\alpha)|+(\gamma'(\alpha))^2}{u(\log y)^2}\Big)\ll \frac{x^{\alpha}}{r^\alpha}\frac{G(\alpha;n)}{\alpha \sqrt{\sigma_2}}u^{-1/3}.
\end{align}
This completes the proof of \eqref{<T0}.

Combining \eqref{S0n}, \eqref{Ebound}, \eqref{>T0} and \eqref{<T0}, we see that 
\begin{align}
\mathcal S_0(n, x/r)=\frac{x^\alpha G(\alpha;n)}{r^\alpha \alpha \sqrt{2\pi \sigma_2}}\Big(1+ O\big(\frac{1}{u^{1/3}}\big)\Big).
\end{align}
Recalling \eqref{arith}, \eqref{Gsdef} and $\Psi(x,y)\asymp \frac{x^\alpha\zeta(\alpha,y)}{\alpha \sqrt{2\pi \sigma_2}}$ (see e.g.\  \cite[Theorem 1]{HT}), we conclude that when $y\geq \mathcal L(x)$ the formula \eqref{summarize} becomes
\begin{align}
\mathcal S(n)&=\prod_{p\mid P_{\mathcal Q}(z)}\mathfrak g_0(p) \sum_{\substack{r\mid 2^\infty\mathcal Q\\ (r, \mathfrak N)=1}}\mathfrak{F}_{\mathcal Q}(r, n)\prod_{\substack{p\nmid \mathcal Q\\ p\mid n}}\frac{\mathfrak g_1(p)}{\mathfrak g_0(p)}\mathfrak c(p)\sum_{\substack{m\leq x/r\\ P^+(m)\leq y}}\mathfrak G(m;n)\\
&\asymp \prod_{p\mid n}\Big(1-\frac{\chi_{-4}(p)}{p}\Big)\prod_{\substack{p\mid n \\ p\leq y}}\Big(1+\frac{\chi_{-4}(p)}{p^\alpha}\Big) \Psi(x,y)\Big(1+ O\Big(\frac{1}{u^{1/3}}\Big)\Big).
\end{align}
 We plug this lower bound back into \eqref{s2-} and choose $C=C(\kappa_0)$ a sufficiently large constant so that $f_5(C\kappa_0) > 0$ getting
 $$S(\mathcal{A}, n^{1/C}) \gg 
 (\log z)^{-2} \mathcal F(n,y)\Psi(x, y)$$
 provided that $y \geq \mathcal L(x)$ and $\mathcal{D}=n^{\kappa_0}$ for some sufficiently small $\kappa_0>0$. 
 This concludes the proof of Theorem \ref{thm2} when $\mathcal L(x)\leq y\leq \exp((\log_2x)^2)$. \\


\begin{thebibliography}{KMV}

\bibitem[AC]{AC} J. Arthur, L.  Clozel, \emph{Simple Algebras, Base Change, and the Advanced Theory of the Trace Formula},  Annals of Math. Studies \textbf{120},  Princeton University Press 1990

\bibitem[ABL]{ABL} E. Assing, V. Blomer, J. Li, \emph{Uniform Titchmarsh divisor problems}, Adv. Math. \textbf{393} (2021), 108076.

\bibitem[Ba]{Ba} A. Balog, \emph{On additive representation of integers}, 
Acta Math. Hungar. \textbf{54} (1989), 297-301.

\bibitem[BBD]{BBD} V. Blomer, J. Br\"udern, R. Dietmann, \emph{Sums of smooth squares},  Compos. Math. \textbf{145} (2009),  1401-1441.


\bibitem[BFI]{BFI2} E. Bombieri, J. B. Friedlander and H. Iwaniec, \emph{Primes in arithmetic progressions to large moduli. II}, Math. Ann. \textbf{277} (1987) 361-393.

\bibitem[BS]{BS} Z. I. Borevich, I. R. Shafarevich, \emph{Number Theory}, Academic Press 1966

\bibitem[BF]{BF} J. Br\"udern, \'E. Fouvry, \emph{
Lagrange's four squares theorem with almost prime variables}, 
J. Reine Angew. Math. \textbf{454} (1994), 59--96.
 
	\bibitem[lBT]{lBT} R. de la Bret\`eche and G. Tenenbaum, \emph{Propri\'et\'es statistiques des entiers friables}, {Ramanujan J.} \textbf{9} (2005), 139--202.


\bibitem[Co]{Co} D. A. Cox, \emph{Primes of the form $x^2+ny^2$},  John Wiley, 1989.


\bibitem[Dr1]{Dr1}S. Drappeau, \emph{Th\'eor\`emes de type Fouvry-Iwaniec pour les entiers friables}, {Compositio Math.} \textbf{151} (2015), 828--862.

\bibitem[Dr2]{Dr2} S. Drappeau, \emph{Sums of Kloosterman sums in arithmetic progressions, and the error term in the dispersion method}, Proc. Lond. Math. Soc.  \textbf{114} (2017), 684-732.

\bibitem[DT]{DT} S. Drappeau, B. Topacogullari, \emph{Combinatorial identities and Titchmarsh's divisor problem for multiplicative functions},   
Algebra Number Theory \textbf{13} (2019),   2383--2425. 

\bibitem[DFI]{DFI} W. Duke, J. Friedlander,  H. Iwaniec. \emph{Bounds for automorphic L-functions},
Invent. Math. \textbf{112} (1993), 1--8.


\bibitem[FT1]{FT90}\'{E}. Fouvry, G. Tenenbaum, \emph{Diviseurs de Titchmarsh des entiers sans grand facteur premier},  in: Analytic number theory (Tokyo, 1988), Lecture Notes in Mathematics \textbf{1434} (Springer, Berlin, 1990), 86--102.

\bibitem[FT2]{FT91}\'{E}. Fouvry, G. Tenenbaum, \emph{Entiers sans grand facteur premier en progressions arithm\'etiques},  Proc. London Math. Soc. (3) \textbf{63} (1991), 449--494.
	
\bibitem[FT3]{FT96}  E. Fouvry, G. Tenenbaum, \emph{R\'epartition statistique des entiers sans grand facteur premier dans les progressions arithm\'etiques},   Proc. London Math. Soc.  \textbf{72} (1996),  481-514.

\bibitem[FT4]{FT22}  E. Fouvry, G. Tenenbaum,\emph{
Multiplicative functions in large arithmetic progressions and applications} 
Trans. Amer. Math. Soc. \textbf{375} (2022), no. 1, 245-299.

\bibitem[FI1]{FI1} J. B. Friedlander, H. Iwaniec, \emph{Opera de cribro}, Coll. Publ. \textbf{57} (2010),  AMS, Providence, RI.

\bibitem[FI2]{FI2} J. B. Friedlander and H. Iwaniec, \emph{Coordinate distribution of Gaussian primes}, J. Eur. Math. Soc.  \textbf{24} (2022),  737--772.

\bibitem[Ha]{Ha} A. Harper, \emph{Bombieri-Vinogradov and Barban-Davenport-Halberstam type theorems for smooth numbers},  {\tt arXiv:1208.5992}

 
\bibitem[HL]{HL} G. H. Hardy, J. E. Littlewood, \emph{Some problems of partitio numerorum; III: on the expression of a number as a sum of primes},  Acta Math. \textbf{44} (1922), 1-70.

\bibitem[HB]{HB} D. R. Heath-Brown, \emph{A new form of the circle method, and its application to quadratic forms}, 
J. Reine Angew. Math. \textbf{481} (1996), 149-206. 

\bibitem[HBT]{HBT} D. R. Heath-Brown, D. I. Tolev, \emph{Lagrange's four squares theorem with one prime and three almost-prime variables}, J. Reine Angew. Math. \textbf{558} (2003), 159--224. 


\bibitem[HT]{HT} A. Hildebrand, G. Tenenbaum, \emph{On integers free of large prime factors.} Trans. Amer. Math. Soc. \textbf{296} (1986),  265-290.

\bibitem[Ho1]{Ho} C. Hooley, \emph{On the representation of a number as the sum of two squares and a prime}, Acta Math. \textbf{97} (1957),
189-210. 

\bibitem[Ho2]{HoBook} C. Hooley, \emph{Applications of sieve methods to the theory of numbers}, Cambridge Tracts in Mathematics \textbf{70},  Cambridge University Press 1976

\bibitem[Ho3]{Ho2} C. Hooley, \emph{On the representation of a number as the sum of a prime and two squares of square-free numbers}, 
Acta Arith. \textbf{182} (2018),   201-229.

\bibitem[Iw1]{Iw} H. Iwaniec, \emph{Primes of the type $\phi(x,y)+A$ where $\phi$ is a quadratic form},  Acta Arith. \textbf{21} (1972), 203-234.

\bibitem[Iw2]{IwS}H. Iwaniec, \emph{Rosser's sieve}, Acta Arith. \textbf{36} (1980), 171-202.

\bibitem[Iw3]{IwS2}H. Iwaniec, \emph{A new form of the error term in the linear sieve}, Acta Arith. \textbf{37} (1980), 307-320.

\bibitem[Iw4]{Iw2} H. Iwaniec, \emph{Topics in classical automorphic forms}, Grad. Stud. Math. \textbf{17} (1997), AMS.

\bibitem[IK]{IK} H. Iwaniec, E. Kowalski, \emph{Analytic Number Theory}, AMS Colloquium Publications \textbf{53}, Providence 2004

\bibitem[Ju]{Ju} M. Jutila, \emph{A variant of the circle method}, in: Sieve methods, exponential sums, and their application in number theory (Cardiff 1995), 245-254.


\bibitem[KMS]{KMS} H. Ki, H. Maier, A. Sankaranarayanan, \emph{Additive problems with smooth integers}, 
Acta Arith. \textbf{175} (2016),  301-319.
    
\bibitem[KMV]{KMV} E. Kowalski, P. Michel, J. VanderKam, \emph{Rankin-Selberg $L$-functions in the level aspect}, Duke Math. J. \textbf{114}
(2002), 123-191.    
    
  \bibitem[XLi]{XLi} X. Li, \emph{Upper bounds on $L$-functions at the edge of the critical strip},  IMRN 2010, 727-755  
    
\bibitem[Li]{Li} Ju. V. Linnik, \emph{An asymptotic formula in an additive problem of Hardy-Littlewood}, Izv. Akad. Nauk SSSR
Ser. Mat. \textbf{24} (1960), 629-706.    


\bibitem[MT]{MT} K. Matom\"aki, J. Ter\"av\"ainen. \emph{On the Möbius function in all short intervals}, J. Eur. Math. Soc. (2022).




\bibitem[Mu]{Mu} R. Munshi, \emph{Shifted convolution of divisor function $d_3$ and Ramanujan $\tau$ function},   Ramanujan Math. Soc. Lect. Notes Ser.  \textbf{20} (2013), 251-260 

\bibitem[Pi]{Pi}  N. Pitt, \emph{On an analogue of Titchmarsh's divisor problem for holomorphic cusp forms}, J. Amer. Math. Soc.  \textbf{26} (2013),
735-776
\bibitem[S]{Sound} K. Soundararajan. \emph{The distribution of smooth numbers in arithmetic progressions}, Anatomy of
Integers, CRM Proc. and Lect. Notes, vol. 46, Amer. Math. Soc., Providence, RI, pp 115--128.
2008


\bibitem[To]{To} B. Topacogullari, \emph{The shifted convolution of divisor functions}, Q. J. Math. \textbf{67} (2016),  331-363. 

\bibitem[Vi]{Vi} A. I. Vinogradov, \emph{General Hardy-Littlewood equation}, 
Mat. Zametki \textbf{1} (1967), 189-197. 
    
\end{thebibliography}
\end{document}